\documentclass[a4paper,10pt]{article}
\usepackage{geometry}
\usepackage[utf8]{inputenc}

\title{Strongly contracting axes for fully irreducible automorphisms of Generalized Baumslag-Solitar groups}
\author{Chloé Papin}

\usepackage{graphicx}
\usepackage{amssymb}
\usepackage{tikz}
\usetikzlibrary{shapes}
\usetikzlibrary{calc}

\usepackage{pstricks}

\usepackage{enumitem}

\usepackage{amsmath}
\usepackage[english]{babel} 
\usepackage[T1]{fontenc}

\usepackage{amsthm}
\newtheorem{theo}{Theorem}[section]
\newtheorem{theointro}{Theorem}
\newtheorem{prop}[theo]{Proposition}

\newtheorem{lem}[theo]{Lemma}

\newtheorem{cor}[theo]{Corollary}

\theoremstyle{definition}
\newtheorem{defi}[theo]{Definition}
\theoremstyle{remark}
\newtheorem{rema}[theo]{Remark}
\newtheorem{remas}[theo]{Remarks}
\newtheorem{ex}[theo]{Examples}

\def \N{\mathbb N}
\def \Z{\mathbb Z}

\def \R{\mathbb R}
\def \D{\mathcal D}
\def \A{\mathcal A}

\def \L{\mathcal L}
\DeclareMathOperator{\Aut}{Aut}
\DeclareMathOperator{\Out}{Out}

\DeclareMathOperator{\diam}{diam}
\DeclareMathOperator{\axe}{Axe}

\DeclareMathOperator{\len}{len}

\DeclareMathOperator{\nbl}{nbl}

\DeclareMathOperator{\BS}{BS}
\DeclareMathOperator{\Inn}{Inn}

\DeclareMathOperator{\val}{val}
\DeclareMathOperator{\Lip}{Lip}
\DeclareMathOperator{\BBT}{BBT}
\DeclareMathOperator{\BCC}{BCC}

\DeclareMathOperator{\vol}{vol}
\DeclareMathOperator{\CV}{CV}
\DeclareMathOperator{\LEG}{Leg}
\DeclareMathOperator{\LR}{LR}

\DeclareMathOperator{\Wh}{Wh}
 
\begin{document}
\newcommand{\localhost}{.}

 \maketitle
 
 \begin{abstract}
  Similarly to the action of $Out(F_N)$ on Outer Space, the outer automorphism group of a Generalized Baumslag Solitar group acts on a \emph{deformation space} endowed with the Lipschitz metric and the action of any fully irreducible automorphism with a train track representative is hyperbolic. Inspired by previous work by Algom-Kfir, we prove that the axis of such an automorphism is strongly contracting.
 \end{abstract}


\section*{Introduction}

Let $G$ be a generalized Baumslag-Solitar (GBS) group, i.e. a group which is isomorphic to the fundamental group of a graph of groups where all vertex and edge groups are infinite cyclic. By Bass-Serre theory, $G$ admits a minimal action by isometries on a simplicial tree $T$ such that all vertex or edge stabilizers in $T$ are infinite cyclic.

In general $G$ admits infinitely many such actions. When $G$ is not isomorphic to $\Z, \Z^2$ or the fundamental group of a Klein bottle $\langle a, b | a^2 = b^2 \rangle = \langle a,t | tat^{-1} = a^{-1} \rangle$, then the \emph{cyclic deformation space} $\D$ associated to $G$ is defined as the projectivized set of minimal actions of $G$ by isometries on metric simplicial trees with edge and vertex stabilizers isomorphic to $\Z$, where actions $T, T'$ are identified if there exists a $G$-equivariant isometry $T \to T'$. The space $\D$ is analogous to Culler and Vogtmann's Outer Space $\CV_N$ for the free group $F_N$. The study of Outer Space is crucial for understanding the outer automorphism group $\Out(F_N)$. Just as $\Out(F_N)$ acts on $\CV_N$, the outer automorphism group $\Out(G)$ acts on $\D$ by pre-composition of the action: if $G$ acts on $T$, then define $T \cdot \phi$ as the action whose underlying space is $T$ and where for $t \in T$, $g \cdot _{T \cdot \phi} t = \phi(g) \cdot _T t$.

There is an important analogy between the study of the outer automorphism group of free groups $F_N$ and mapping class groups. In this analogy, Outer Space is the counterpart for Teichmüller space. Fully irreducible automorphisms of $\Out(F_N)$ are a special class of automorphisms which do not preserve the conjugacy class of any free factor. In some sense they act on Outer Space by translations along an axis whose points are actions which admit a \emph{train track representative}. Their equivalent in the mapping class groups context is pseudo-Anosov automorphisms. Such automorphisms can also be defined for GBS groups: we say that an automorphism $\phi \in \Out(G)$ is \emph{fully irreducible} if no conjugacy class of proper \emph{cyclic factor} is $\phi$-periodic. A \emph{cyclic factor} is the analogue of a free factor for the free group: it is a non-cyclic subgroup $G' < G$ such that there exists a graph of groups $\Gamma$ representing $G$, and a subgraph $\Gamma' \subset \Gamma$ such that the fundamental group $\pi_1(\Gamma')$ identifies to the subgroup $G'$. We say that $G'$ is a proper cyclic factor if $G' \neq G$.

One main property for a fully irreducible automorphism $\phi \in \Out(F_N)$ is that its action on $\CV_N$ is hyperbolic and admits an axis $\L_\phi$. In order for this to make sense, we endow $\CV_N$ with the non-symmetric \emph{Lipschitz metric} $d_{\Lip}$: it is defined by
\[
 d_{\Lip}(T,T'):= \Lip(T,T') \frac{\vol(T)}{\vol(T')}
\]
for $T, T' \in \CV_N$, where $\Lip(T,T') := \sup_{f: T \to T' \text{ Lipschitz }} \Lip(f)$.

One can define a closest point projection $\pi_f: \D \to \L_f$ such that for $T \in \D$, $d_{\Lip}(T, \pi_f(T))$ is minimal. Define the \emph{outward ball} $B_{\rightarrow}(Y,r):=\{T \in \D/ d_{\Lip}(Y,T) <r\}$.

In \cite{AlgomKfirStrongly} Algom-Kfir proves that the projections on axes have the \emph{strong contraction property}: 
there exists $B > 0$ depending only on $N$ and $\phi$ such that the diameter of the projection of any outward ball in $\D$ disjoint from $\L_\phi$ is bounded by $B$.

Likewise we can define these objects for the deformation space $\D$ of a GBS group. A few differences arise, for example because trees in $\D$ have non-trivial edge stabilizers: one consequence is that in general there exist $T, T' \in \D$ such that $d_{\Lip}(T, T')= 0$ and $d_{\Lip}(T',T)$ is arbitrarily large.

In this paper we prove the analogue of this property for fully irreducible automorphisms of a GBS group:
\begin{theointro}
Let $G$ be a GBS group with $b_1(G) \geq 3$.
 Let $\phi$ be an atoroidal fully irreducible automorphism such that $\phi, \phi^{-1}$ both admit train track representatives. 
 Let $\L_\phi$ be an axis for $\phi$ in $\D$ and let $\pi_\phi$ be a closest point projection to $\L_\phi$. Then there exists $D > 0$ such that for any $Y \in \D$ and $r>0$ such that $B_{\rightarrow}(Y, r) \cap \L_\phi = \varnothing$
 \[
  \diam(\pi_\phi(B_{\rightarrow}(Y, r))) \leq D
 \]
\end{theointro}

The proof for this result is similar to \cite{AlgomKfirStrongly}. 

We do not know whether train track representatives always exist for fully irreducible automorphism. Several cases are known:
\begin{itemize}
 \item The procedure in \cite{BestvinaHandelTrainTracks} may be adapted to construct a train track map by hand on a specific example. See Example \ref{ex:train-track}.
 \item If $G$ has no non-trivial integer modulus then \cite{ForesterSplittings} proves that $\D$ has finite dimension. Then \cite{MeinertTheLipschitzMetric} implies that all fully irreducible automorphisms admit train track representatives.
 \item If $G := \BS(p,pn)$ then Bouette proved in \cite{BouetteThese} that all automorphisms of $G$ are reducible and preserve the conjugacy class of a common cyclic factor $H$. She then introduces a new deformation space $\D_H$ consisting in all $G$-trees with cyclic edge stabilizers, and vertex stabilizers either cyclic or conjugate to $H$. In this deformation space, there exist fully irreducible automorphisms and they admit train track representatives.
\end{itemize}
The second case is not relevant here since we need the trees in the deformation space to be locally finite.

For technical reasons, we ask that the automorphism $\phi$ in the theorem be pseudo-atoroidal, which means that for all $g \in G$, $(\|\phi^n(g)\|_T)_{n \in \N}$ is unbounded.  

We also need the first Betti number $b_1(\Gamma)$ for any graph of cyclic groups $\Gamma$ with $\pi_1(\Gamma) \simeq G$ to be at least $3$. Actually $b_1(\Gamma)$ does not depend on the choice of $\Gamma$ when $G$ is not isomorphic to the fundamental group of a Klein bottle, which we exclude.

\bigskip

An important result in \cite{AlgomKfirStrongly} is the fact that axes of automorphisms have bounded projection on each other. This fact is of great interest since it enables the definition of a projection complex on which the quasi-tree construction of \cite{BBF} could be applied. However we do not know yet if this still holds in the GBS context, due to the fact that the Lipschitz metric is even less symmetric. For example it is not true that the Lipschitz metric is quasi-symmetric on the $\theta$-thick part of $\D$. Furthermore, there is no bound on the number of candidates in trees of $\D$.

\bigskip

In Section \ref{sec:axes-definitions} we give some background about GBS groups and their automorphisms. We develop the topic of laminations in Section \ref{sec:axes-laminations}. In Section \ref{sec:whitehead-et-laminations} we prove the analogue of 
results from \cite{AlgomKfirStrongly} which state that the axis of a simple element cannot follow both the stable and unstable lamination for a long distance; the method of the proof differs somehow from the original. Section \ref{sec:legalite} develops the behaviour of lines in $T$ such as axes of elements of $g$ when iterating a train track $f: T \to T$, and Section \ref{sec:projection} relies on it to define a projection of $\D$ on the axis of a fully irreducible element with a train track representative. The contents in this section are really close to \cite{AlgomKfirStrongly} and are there for completeness.

We prove negative curvature properties of the projection in Section \ref{sec:negative}, that is, inequalities about distances in $\D$. 
Although the former does not differ from the free group case, the latter needs some arguments which are specific to GBS groups. Finally we prove the strong contraction for balls of outward radius.

 \paragraph{Aknowledgement.}
I am especially grateful to my advisor Vincent Guirardel for his support and guidance. This article is a part of my thesis, which was written at Université de Rennes 1.

\section{Generalities} \label{sec:axes-definitions}
\subsection{Graphs and trees}  

A \emph{graph} $\Gamma$ is defined  by $(V(\Gamma), E(\Gamma), \bar \cdot, o, t)$ where
\begin{itemize}
 \item $V(\Gamma)$ is a set of \emph{vertices}
 \item $E(\Gamma)$ is a set of \emph{edges}
 \item the map $\bar \cdot$ is an involution $E(\Gamma) \to E(\Gamma)$ without fixed point; for $e \in E(\Gamma)$ the edge $\bar e$ is called the \emph{opposite edge}
 \item the maps $o, t : E(\Gamma) \to V(\Gamma)$ are the \emph{initial vertex} and \emph{terminal vertex} maps, with the property that every $e \in E(\Gamma)$ satisfies $o(e) = t(\bar e)$.
 \end{itemize}
It is finite if $V(\Gamma), E(\Gamma)$ are finite. See \cite{SerreArbresAmalgames} for more details on graphs.
 
 An edge path in $\Gamma$ is a sequence $e_1, \dots, e_k$ with $e_i \in E(\Gamma)$ for $i \in \{1, \dots, k\}$ and $t(e_i)= o(e_{i+1})$ for $i \leq k-1$. It is \emph{non-backtracking} if for all $i \in \{1, \dots, k-1\}$, $\bar e_i \neq e_{i+1}$. It is a \emph{loop} if $o(e_1) = t(e_k)$.
 
A \emph{tree} is a graph without non-backtracking loops.

Let $E^+(\Gamma)$ be an orientation of the edges, i.e. a subset of $E(\Gamma)$ such that $E^+(\Gamma) \sqcup \overline{E^+(\Gamma)}$ is a partition of $E(\Gamma)$.

A \emph{metric} on a graph $\Gamma$ is a map $\len_{\Gamma}: E(\Gamma) \to \R_+$ such that for all $e \in E(\Gamma)$, $\len_\Gamma(\bar e) = \len_\Gamma(e)$.

The \emph{geometric realization} of a graph $\Gamma$ is the union of points $(x_v)_{v \in V(\Gamma)}$ and segments $(\sigma_e)_{e \in E^+(\Gamma)}$ where $\sigma_e$ is isometric to $[0,\len_{\Gamma}(e)]$ for every $e \in E(\Gamma)$, where for every $e \in E^+(\Gamma)$ we identify the first point of $\sigma_e$ with $x_{o(e)}$ and its last point with $x_{t(e)}$. It is endowed with the associated path metric. It does not depend on the choice of $E^+(\Gamma)$. 

In the rest of the paper we will identify trees and other graphs with their geometric realizations. A \emph{path} in a tree $T$ is the image of a Lipschitz map from an interval to $T$. It is \emph{non-backtracking} if the map is an immersion, and equivalently if it is a the image of a geodesic. In the context of geometric realizations an \emph{edge path} is a path which is the image of an edge path in the graph. For two points $x, y \in T$, the segment $[x,y]$ is the unique geodesic from $x$ to $y$.

\bigskip

A \emph{graph of groups} is a graph $\Gamma$ together with collections of vertex groups $(G_v)_{v \in V(\Gamma)}$ and edge groups $(G_e)_{e \in E(\Gamma)}$ and monomorphisms $\iota_e: G_e \to G_{t(e)}$.  Let $\tau$ be a maximal subtree in the graph $\Gamma$. The fundamental group $\pi_1(\Gamma, \tau)$ of the graph of groups $\Gamma$ is defined as follows:
\[
 \pi_1(\Gamma, \tau) = \left \langle \bigcup_{v \in V(\Gamma)} G_v \cup (t_e)_{e \in E(\Gamma)} | \bigcup_{v \in V(\Gamma)} R_v,  \bigcup_{e \in E(\Gamma)} R_e, R_{\tau} \right \rangle 
\]
where
\begin{itemize}
 \item for $v \in V(\Gamma)$, $R_v$ is the set of relations of $G_v$
 \item for $e \in E(\Gamma)$, $R_e = \{ t_e \phi_e(h) t_{\bar e} \phi_{\bar e}(h)^{-1} / h \in G_e\}$ 
 \item $R_\tau := \{ t_e, e \in \tau\}$
\end{itemize}
Note that for $h=1$ we obtain the relation $t_{\bar e} = t_e^{-1}$.
Up to isomorphism the fundamental group does not depend on the choice of $\tau$. If the vertex groups and edge groups are finitely generated then $\pi_1(\Gamma)$ is finitely presented.

Let $G$ be a group. A \emph{marked graph of groups} for $G$ is a graph of groups $\Gamma$ together with a \emph{marking} (i.e. an identification) $\Psi: G \to \pi_1(\Gamma)$. The automorphism group of $G$ acts as follows on the set of marked graphs of groups: if $\phi \in \Aut(G)$ and $(\Gamma, \Psi)$ is a marked graph of groups then $(\Gamma, \Psi) \cdot \phi := (\Gamma, \Psi \circ \phi)$.

A \emph{$G$-tree} $T$ is a metric simplicial tree with an action of $G$ by isometries. The tree $T$ is \emph{minimal} if there is no proper $G$-invariant subtree.

The \emph{universal cover} of a graph of groups $\Gamma$ is a minimal $\pi_1(\Gamma)$-tree $T$ such that $T/G$ is isomorphic to $\Gamma$ as a graph and for every $v \in V(\Gamma)$ and every lift $\bar v \in T$, the stabilizer of $\bar v$ is isomorphic to $G_v$. By Bass-Serre theory in \cite{SerreArbresAmalgames}, universal covers exist and are unique up to $\pi_1(\Gamma)$-equivariant isomorphism. 

Moreover Bass-Serre theory gives a correspondance between marked graphs of groups for $G$ and $G$-trees.

If $\Gamma$ is a metric graph then the metric naturally lifts to its universal cover.

In a $G$-tree $T$ we denote the pointwise stabilizer of a vertex $v$ (resp. an edge $e$) by $G_v$ (resp. $G_e$). 

\bigskip

A \emph{generalized Baumslag-Solitar} (GBS) group is a group which is isomorphic to the fundamental group of a finite graph of groups where all vertex and edge groups are infinite cyclic. If a generator is chosen for every vertex and edge group then the monomorphisms $\phi_e$ are defined by the multiplication by an integer $\lambda(\bar e) \in \Z  \setminus \{0\}$.

Let $G$ be a GBS group. In the general case there exist infinitely many marked graphs of cyclic groups for $G$. Making $\Aut(G)$ act on a marked graph of groups often yields infinitely many other markings, besides in general cases there are infinitely many possible underlying graphs of groups.

Let $T$ be a $G$-tree. A subgroup $H < G$ is \emph{elliptic} in $T$ if it fixes a point in $T$.
Suppose all elliptic groups in $T$ are also elliptic in $S$. Then there exists a $G$-equivariant map $T \to S$ (see by example \cite{GuirardelLevitt07}). 

Let $T, T'$ be $G$-trees. We say that they lie in the same \emph{deformation space} if they have the same sets of elliptic subgroups.
Equivalently they are in the same deformation space if there exist $G$-equivariant maps $T \to T'$ and $T' \to T$.

Now let us define the \emph{cyclic deformation space} $\D$ associated to a group $G$, as the set of minimal simplicial $G$-trees with cyclic vertex and edge stabilizers, where we identify $T$ and $T'$ if there is a $G$-equivariant isometry or homothety $T \to T'$.

Equivalently we could define $\D$ with marked graphs of groups.

If $G$ is not isomorphic to $\Z^2$ or the fundamental group of a Klein bottle $\langle a,b| a^2 = b^2 \rangle \simeq \langle a,t| tat^{-1} = a^{-1}\rangle$, then $\D$ is a deformation space, i.e. all trees in $\D$ have the same elliptic subgroups. 

Let $T \in \D$ and let $e$ be an edge in $T$. Define the equivalence relation $\sim_e$ as the minimal $G$-invariant equivalence relation such that $x \sim_e y$ if $x, y \in e$. The \emph{collapse} of the edge $e$ is the quotient map $T \to T/ \sim_e$. The edge $e$ in $T$ is \emph{collapsible} if $T/ \sim_e \in \D$. Equivalently an edge is collapsible if its image  in the quotient is not a loop and one of its two labels is $\pm 1$. 

A tree of $\D$ is \emph{reduced} if none of its edges is collapsible.

GBS trees in the same deformation space share some properties. Let $\Gamma$ be a finite connected graph. Then the first Betti number $b_1(\Gamma)$ is defined by $b_1(\Gamma) = \# E(\Gamma) - \#V(\Gamma) + 1$. By \cite[Section 4]{GuirardelLevitt07} the first Betti number is an invariant of the deformation space.

Let $G$ be a GBS group with cyclic deformation space $\D$. We say that an elliptic subgroup $H<G$ is \emph{big} if there exists a tree $T \in \D$ such that $H$ fixes no edge in $T$. 

From \cite{GuirardelLevitt07} we deduce:
\begin{lem} \label{lem:big}
 Let $T \in \D$. 
 The number of vertices $v \in T$ such that for all edge $e$ with origin $v$, $G_e \neq G_v$ is bounded by the number of conjugacy classes of big subgroups of $G$.

 If $T$ is reduced then these numbers are equal.
\end{lem}

\begin{rema}
 The notion of big subgroups is defined in \cite{GuirardelLevitt07}, though it depends on a family $\A$ of subgroups of $G$: a subgroup $H<G$ is big if it is elliptic and is not conjugate into a subgroup of an element of $\A$. Here the corresponding choice for $\A$ is the family of subgroups which fix an edge in a reduced tree of $\D$, or equivalently in every reduced tree of $\D$. Thus an elliptic subgroup is big if fixes a single point in some (equivalently any) tree in $\D$.
\end{rema}

Solvable GBS groups are GBS groups isomorphic to $\Z$ and $\BS(1,n)$ for $n \in \N$ (which include $\Z^2$ and the fundamental group of a Klein bottle).

\subsection{Cyclic factors, irreducible automorphisms}

From now on we assume $G$ is a non-solvable GBS group. The automorphism group of $G$ is $\Aut(G)$. The outer automorphism group is $\Out(G):= \Aut(G)/\Inn(G)$ where $\Inn(G)$ is the subgroup of inner automorphisms $\{c_g : x \mapsto gxg^{-1}, g \in G\}$.

Cyclic factors are the GBS analogue of free factors for free groups.
\begin{defi}
 A \emph{cyclic factor} of $G$ is a subgroup $H$ such that there exists a graph of cyclic groups $\Gamma$ and an identification $G \simeq \pi_1(\Gamma)$, with a subgraph $\Gamma_H$ such that $H$ is conjugate to $\pi_1(\Gamma_H)$.
\end{defi}
The family of cyclic factors of $G$ is stable by conjugacy and by automorphisms.

\begin{ex}
 \begin{enumerate}
  \item If $G := \BS(2,4) = \langle a, t| ta^2t^{-1} = a^4\rangle$, the first graph of groups of Figure \ref{fig:ex-facteur-cyclique} represents $G$. The red subgraph represents the subgroup $H :=\langle a, t^{-1} a^2 t \rangle$ which is a cyclic factor.
  \item The second and third graphs of Figure \ref{fig:ex-facteur-cyclique} represents $G:= \langle u, r,s,t | tu^n t^{-1} = su^n s^{-1} = ru^n r^{-1} = u \rangle$. The subgroup $\langle u, r \rangle \simeq \BS(1,n)$ is a cyclic factor, it can be seen in the graph on the left. The subgroup $\langle u, ru r^{-1}, rs u s^{-1} r^{-1}, rst \rangle \simeq \BS(1, n^3)$ is a cyclic factor which can be seen in the graph on the right.
 \end{enumerate}
 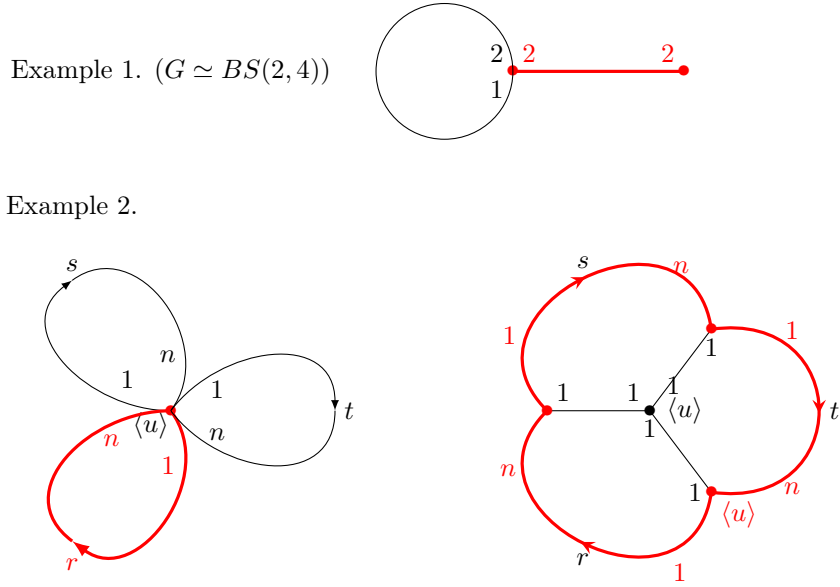
\begin{figure}
  \centering
  \begin{tikzpicture}[scale = 0.9]
 
\newcommand{\arrowIn}{
\tikz \draw[-stealth, very thick] (-1pt,0) -- (1pt,0);
}

  \begin{scope}[xshift= 6cm]
    \coordinate (w) at (0,0);
    \coordinate (v) at (2.5,0);

    \draw (v) node [color = red] {$\bullet$} node [above left, color = red] {$2$} ;
    \draw (w) node [color = red] {$\bullet$} node [above right, color = red]{$2$} node [above left] {$2$} node [below left] {$1$};
    \draw [very thick, color=red] (v)--(w);
    \draw (w) arc (0: 360:1);

    \coordinate (leg) at (-5,0);
    \draw (leg) node {Example 1. ($G \simeq BS(2,4)$)};
  \end{scope}

  \begin{scope}[xshift=1cm, yshift=-5cm, scale = 1.2]
   \coordinate (leg) at (-1.2,2.5);
   \draw (leg) node {Example 2.};
   
    \coordinate (o) at (0,0);
    \coordinate (r) at (3/5,4/5);
     \coordinate (s) at (3/5,-4/5);
     \coordinate (t) at (-1,0);
    \coordinate (r') at (-6/5,-8/5);
     \coordinate (s') at (-6/5,8/5);
     \coordinate (t') at (2,0);

    \coordinate (r+) at (-2/5,-11/5);
     \coordinate (s+) at (-2,1);
     \coordinate (t+) at (2,1);

    \coordinate (r-) at (-2,-1);
     \coordinate (s-) at (-2/5,11/5);
     \coordinate (t-) at (2,-1);

    \draw (o) node [below left = -0.1cm] {$\langle u \rangle$};

    \draw[color = red] (o) node {$\bullet$};
    \draw[-latex] (o) .. controls (r)  and (t+) .. (t') node [near start, below] {$1$} node [right] {$t$};
    \draw[-latex, color = red, very thick] (o) .. controls (s) and (r+) .. (r') node [near start, left] {$1$} node [below = 0.1cm] {$r$} ;
    \draw[-latex] (o) .. controls (t) and (s+) .. (s') node [near start, above right ] {$1$} node [above] {$s$};
    \draw (o) .. controls (s)  and (t-) .. (t') node [near start, above] {$n$};
    \draw [color = red, very thick] (o) .. controls (t) and (r-) .. (r')  node [near start, below] {$n$};
    \draw (o) .. controls (r) and (s-) .. (s')  node [near start, left] {$n$} ;
\end{scope}
\begin{scope}[xshift=8cm, yshift = -5cm, scale = 1.5]
    \coordinate (o) at (0,0);
    \coordinate (r) at (3/5,4/5);
     \coordinate (s) at (3/5,-4/5);
     \coordinate (t) at (-1,0);
    \coordinate (r') at (-6/5,-8/5);
     \coordinate (s') at (-6/5,8/5);
     \coordinate (t') at (2,0);

     \draw (o) node [right = 0.1cm] {$\langle u \rangle$};
     \draw[color = red] (s) node [below right] {$\langle u \rangle$};
    
    \coordinate (r+) at (2/5,11/5);
     \coordinate (s+) at (2,-1);
     \coordinate (t+) at (-2,-1);

    \coordinate (r-) at (2,1);
     \coordinate (s-) at (2/5,-11/5);
     \coordinate (t-) at (-2,1);

  \draw (o) node {$\bullet$} node [above right=0.1cm] {$1$} node [below] {$1$} node [above left] {$1$};
    \draw (o) -- (r) node [below ] {$1$}  ;
    \draw (o) -- (s) node [left] {$1$}; 
    \draw (o) -- (t) node [above right] {$1$};

   \draw[very thick, color = red]  (r) .. controls (r-) and (s+) .. (s) node [near start, above] {$1$} node [near end, below] {$n$} node {$\bullet$} node[
    sloped,
    pos=0.5,
    allow upside down]{\arrowIn} node [midway, right, color = black] {$t$};
   \draw[color = red, very thick] (s) .. controls (s-) and (t+) .. (t) node [near start, below right] {$1$} node [near end, above left] {$n$} node {$\bullet$} node[
    sloped,
    pos=0.5,
    allow upside down]{\arrowIn} node [midway, below, color =black] {$r$};
   \draw[color = red, very thick]  (t) .. controls (t-) and (r+) .. (r) node [near start,left] {$1$} node [near end, right] {$n$} node {$\bullet$} node[
    sloped,
    pos=0.5,
    allow upside down]{\arrowIn} node [midway, above, color = black] {$s$};

  \end{scope}
\end{tikzpicture} 
  \caption{Examples of cyclic factors} \label{fig:ex-facteur-cyclique}
 \end{figure}
\end{ex}

\begin{defi}
An automorphism $\phi \in \Aut(G)$ is \emph{fully irreducible} if no power of $\phi$ preserves the conjugacy class of a cyclic factor. Since inner automorphisms preserve conjugacy classes, the full irreducibility can be defined for outer automorphisms.

A \emph{representative} for $\phi$ is a map $f:T \to T$ with $T \in \D$ which is $\phi$-equivariant, i.e. $\forall t \in T$, $\forall g \in G$, $f(gt) = \phi(g) \cdot f(t)$.

A representative for an outer automorphism class $\psi \in \Out(G)$ is a representative for some automorphism in the class $\psi$.
\end{defi}

\begin{defi}
 Let $\phi \in \Out(G)$. A \emph{pseudo-periodic} conjugacy class for $\phi$ is the conjugacy class of an element $g \in G$ such that $\|\phi^n(g) \|_{n \in \N}$ is bounded.
 
 An automorphism $\phi \in \Out(G)$ is \emph{pseudo-atoroidal} if $\phi$ has no pseudo-periodic conjugacy class.
\end{defi}

Train track representatives for automorphisms of $\Out(F_N)$ were introduced in \cite{BestvinaHandelTrainTracks}. They are a main tool for studying these automorphisms. One can define likewise train tracks for other families of groups acting on trees.

\begin{defi}
 Let $T$ be a $G$-tree.
 A \emph{gate structure} on $T$ is a $G$-invariant family of equivalence relations $(\sim_v)_{v \in V(T)}$ on the sets $E_v$ of edges with origin $v$. Equivalence classes for these relations are called gates.
 
 Let $\tau: T \to T'$ be a $G$-equivariant map sending edges to non-degenerate non-backtracking paths. The gate structure \emph{induced by $\tau$} is the minimal gate structure such that for $v \in V(T)$, $e,e' \in E_v$, if $\tau(e)\cap\tau(e')$ has non-zero length then $e \sim_v e'$. 
\end{defi}

\begin{defi}
 Let $T$ be a $G$-tree with a gate structure. A \emph{turn} in $T$ is a pair of edges with same origin. The turn $\{e, e'\}$ is \emph{illegal} if $e$ and $e'$ belong to the same gate. Otherwise the turn is \emph{legal}.
\end{defi}

\begin{defi}
 Let $\phi \in \Out(G)$. Let $f: T \to T$ be a representative for $\phi$ sending vertex to vertex.. Then $f$ is \emph{train track} if for every $e \in E(T)$:
 \begin{itemize}
  \item $\len(f(e)) > 0$
  \item at every vertex $v \in V(T)$ there are at least two gates for the gate structure induced by $f$.
  \item for every $k \in \N$, $f^k(e)$ crosses only legal turns for the gate structure induced by $f$.
 \end{itemize}
 The \emph{train track structure} is the gate structure such that $e \sim e'$ if there exists $k \in \N$ such that $f^k(e) \cap f^k(e')$ is not a single point. 
\end{defi}

Suppose $f : T \to T$ is a train track representative for a fully irreducible automorphism $\phi \in \Out(G)$. Up to precomposing $f$ with a map $T \to T$ whose restriction to edges is a homeomorphism, we may assume that $f$ stretches the edges uniformly. Let $e_1, \dots, e_n$ be the edges of $T/G$. Let $A(f)$ be the transition matrix for $f$ where $A(f)_{ij}$ is the number of occurences of edges in the orbit $e_i$ in $f(e_j)$, for $1 \leq i,j \leq n$. Irreducibility of $\phi$ implies that up to collapsing edges in $T$ the matrix $A(f)$ is primitive, i.e. there exists $k \in \N$ such that $A(f)^k >0$ (see \cite[Lemma 1.9]{PapinDetection}). 

The theorem of Perron-Frobenius below then applies to $A(f)$:
\begin{theo}[Perron-Frobenius] 
 Let $A$ be a non-negative primitive matrix with size $n\times n$. There exists a real eigenvalue $\lambda>0$ (the \emph{Perron-Frobenius eigenvalue}) such that for every other eigenvalue $\mu \neq \lambda$ we have $|\mu|< \lambda$. The eigenvectors for $\lambda$ are unique up to scalar multiplication and there exists an eigenvector $v$ for $\lambda$ such that $v > 0$.
\end{theo}
A proof of the theorem can be found in \cite[Theorem 1.1]{Seneta}.

Let $\lambda$ be the Perron-Frobenius eigenvalue and $(l_1, \dots, l_n)$ be the left Perron-Frobenius eigenvector. Define a metric on $T$ by $\len(e_n):=l_n$. Then for every $e \in E(T)$ we have $\len(f(e)) = \lambda \len(e)$, and the Lipschitz constant $\Lip(f)$ is $\lambda$.

From now on we assume train track maps are linear on edges and have the same Lipschitz constant on all edges. 

\bigskip

Recall that it is not known whether all fully irreducible automorphisms admit train track maps. However some do exist: Example \ref{ex:train-track} is an example of fully irreducible automorphism with a train track representative.

\begin{ex}\label{ex:train-track}
 Let $G:= \langle r,s,t | ru^n r^{-1} = su^ns^{-1} = a
 tu^nt^{-1} = u\rangle$. Define $\phi \in \Aut(G)$ by
\[
 \phi :  \begin{cases}
           u \mapsto u \\
           r \mapsto s \\
           s \mapsto t \\
           t \mapsto rsts^{-1} t^{-1}
         \end{cases}
\]
Define the tree $T$ (whose quotient $T/G$ is represented on Figure \ref{fig:ex-train-track}) by a fundamental domain with vertices $v,x$ and edges
\begin{align*}
 e_a & = [ v, rv] \\
 e_b & = [v, t x ] \\
 e_e & = [v, x] \\
 e_f & = [v, s^{-1}x] \\
\end{align*}
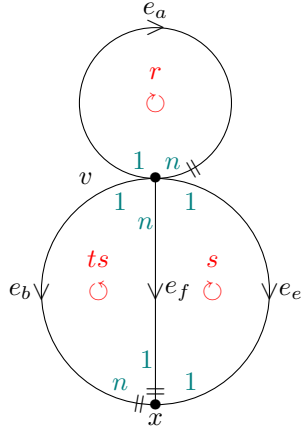
\begin{figure}
 \centering
   \begin{tikzpicture}
  \coordinate (v) at (0,2);
  \coordinate (x) at (0,-1);

 \draw (v) node {$\bullet$} node [left = 0.7cm] {$v$} arc (-90:270:1) node [pos=0.5, sloped] {$>$} node [pos=0.5, above] {$e_a$} node [pos=0.08, rotate= -60] {$=$};
 \draw (v) arc (90:-270:1.5) node [pos=0.05, below, teal] {$1$} node [pos = 0.55,above, teal] {$n$} node [pos = 0.45, above, teal] {$1$}  node [pos=0.95, below, teal] {$1$}  node [pos=0.25, sloped] {$>$} node [pos=0.25, right] {$e_e$} node [pos=0.75, sloped] {$<$} node [pos=0.75, left] {$e_b$} node [pos=0.52, rotate=80] {$=$};;
 \draw (v) -- (x) node {$\bullet$} node [below] {$x$} node [pos=0.2, left = -0.1cm, teal] {$n$} node [pos = 0.8, left = -0.1cm, teal] {$1$}  node [pos=0.5, sloped] {$>$} node [pos=0.5, right] {$e_f$} node [pos=0.95] {$=$};
\draw (v) node [above left, teal] {$1$} node [above right, teal] {$n$} ;

\draw[red] (0, 3) node {$\circlearrowright$} node [above=0.2cm] {$r$};
\draw[red] (-0.75, 0.5) node {$\circlearrowleft$} node [above=0.2cm] {$ts$};
\draw[red] (0.75, 0.5) node {$\circlearrowright$} node [above=0.2cm] {$s$};
\end{tikzpicture}
 \caption{The quotient $T/G$ in Example \ref{ex:train-track}} \label{fig:ex-train-track}
\end{figure}
Define $f$ on $T$ by
 \begin{align*}
f : & v \mapsto v, x \mapsto tx \\
    &e_a \mapsto e_e- s \cdot \bar e_f \\
    &e_b \mapsto e_a- r\cdot e_e - rs \cdot \bar e_f - rs \cdot e_b - rst \cdot \bar e_e - rst \cdot e_f \\
    & e_e \mapsto e_b  \\
    & e_f \mapsto e_e  
\end{align*}
The map $f$ is a train track representative for $\phi$. To see this we can compute the successive images of the turns taken by $f$. Consider the turn $\{\bar e_a, r e_e\}$ at vertex $r v$ which is taken by $f(e_b)$:

\begin{center}
\begin{tikzpicture}[scale =0.95]
 \node (0) at (0,0) {$r\{r^{-1}\bar e_a, e_e\}$};
  \node (1) at (3,0) {$\phi(r)\{e_f, e_b\}$};
  \node (2) at (6,0) {$\phi^2(r)\{e_e , e_a \}$};
  \node (3) at (9,0) {$\phi^3(r)\{e_b , e_e \}$};
  \node (4) at (12,0) {$\phi^4(r)\{e_a, e_b \}$};
  
  \draw [-latex] (0) --  (1) node [midway, above] {$f$};
  \draw [-latex] (1) -- (2) node [midway, above] {$f$};
  \draw [-latex] (2) -- (3) node [midway, above] {$f$};
  \draw [-latex] (3) -- (4) node [midway, above] {$f$};
  \draw [-latex, dashed] (4) .. controls (12, -2) and (6,-2) ..  (2)node [midway, above] {$f$} ;
\end{tikzpicture}
\end{center}

When applying $f$ to $\phi^4(r)\{e_a, e_b \}$ we find a turn in the orbit of $\phi^2(r)\{e_e , e_a \}$ which we had already found. Thus $\{\bar e_a, r e_e\}$ is a legal turn. The proof for the other turns goes the same way.

\bigskip

The inverse of $\phi$ is 
\[
 \phi^{-1} :  \begin{cases}
           u \mapsto u \\
           r \mapsto tsrs^{-1}r^{-1} \\
           s \mapsto r \\
           t \mapsto s
         \end{cases}
\]
which is the same automorphism as $\phi$, but with swapped roles for $r$ and $t$. Thus it admits a train track map.
\end{ex}

\bigskip

Let $f: T \to S$ be a $G$-equivariant map between trees of $\D$, sending vertex to vertex and edge to edge path. By \cite{ForesterDeformation} it is a quasi-isometry, i.e. a map such that for all $x,y \in T$, $K^{-1} d_{T}(x,y) - C \leq d_{T'}(h(x), h(y)) \leq K d_T(x, y) + C$. Thus the map has bounded backtracking property (see \cite{GJLL}) and there exists a constant $\BBT(f) \leq  K^2 C + C$ such that for every $x,y \in T$ the image $f([x,y])$ lies in a $\BBT(f)$-neighbourhood of $[f(x), f(y)]$.

In \cite{BestvinaFeighnHandelLaminations} and \cite{AlgomKfirStrongly}, a similar constant $\BCC(f)$ is used. Let us clarify the link between these constants.
Suppose $\alpha \cdot \beta$ is a geodesic concatenation of paths with $\alpha$ $f$-legal. When applying $f$ to $\alpha \cdot \beta$, there is a subsegment $\tau \subset \alpha$ such that $f(\tau) \subset f(\beta)$. Let $\tau' \subset \beta$ be a minimal prefix such that $f(\tau) \subset f(\tau')$. The first point of $\tau$ and last point of $\tau'$ are mapped to the same point. Then because $f$ is a quasi-isometry there is a constant $K$ depending only on $f$ such that the length of $\tau \cdot \tau'$ is bounded by $K$. This gives a bound $\BCC(f)$ (from \emph{bounded cancellation constant}, introduced by Cooper in \cite{CooperAutomorphisms}) on $f(\tau)$: the simplification which occurs at an illegal turn is bounded by $\BCC(f)$. By applying the previous paragraph to $\tau \cdot \tau'$ we have $\BCC(f) \leq \BBT(f)$. In order to keep a reduced number of constants, we will use the $\BBT$ constant instead of $\BCC$ where it could be used.

For a path $\eta$ in a tree $T$, we denote by $[\eta]$ the unique non-backtracking path which has the same endpoints as $\eta$. We can extend this notation to infinite paths which converge to a point in the boundary of $T$: for example, this is well-defined for a bi-infinite quasi-geodesic such as the image of a line by a quasi-isometry.

\begin{defi} \label{defi:critical-constant}
Let $\phi \in \Out(G)$ be a fully irreducible automorphism with a train track representative $f : T \to T$. The \emph{critical constant} is $C_f := \frac{2\BBT(f)}{\lambda -1}$. 
\end{defi}

The critical constant has the following property: for any geodesic concatenation $\alpha \cdot \beta \cdot \gamma \subset T$ such that $\beta$ is legal and $\len(\beta) \geq C_f$ then let $\alpha' \subset \alpha, \beta' \subset \beta, \gamma' \subset \gamma$ such that the path $[f(\alpha \cdot \beta \cdot \gamma)]$ can be written as the concatenation $[f(\alpha')] \cdot [f(\beta')] \cdot [f(\gamma')]$. Then $\len(f(\beta')) \geq C_f$. More specifically we have:

\begin{lem} \label{lem:thetas-disjoints}
 Let $\phi \in \Out(G)$ be fully irreducible. Let $f : T \to T$ be a train track representative.
 Let $\alpha$ be any path in $T$. Let $\beta \subset \alpha$ be a legal subpath with length at least $2C_f$. Define $\beta'$ as the legal subpath of $\beta$ obtained by cutting out the $\frac{C_f}{2}$-neighbourhood of the endpoints. Then $\beta'$ satisfies the following condition: for all $n \in \N$
 \[
  f^n(\alpha \setminus \beta') \cap f^n(\beta')= \varnothing \text{ and } f^n(\beta') \subset [f^n(\alpha)]
 \]
 In particular $\len([f^n(\alpha)]) \geq \lambda^n \len(\beta')$.
\end{lem}
\begin{proof}
 Let $\alpha, \beta$ be as above. 
 
 We actually prove a slightly stronger statement by induction on $n$.
 Define $\beta'_n$ as the segment obtained by cutting out a $\BBT(f) \displaystyle \sum_{k=1}^{n} \lambda^{-k}$-neighbourhood from the endpoints of $\beta$. Observe that $\beta' = \bigcap_{n \in \N} \beta'_n$ since $\BBT(f) \displaystyle \sum_{k=1}^{\infty} \lambda^{-k} = \frac{\BBT(f)}{\lambda - 1} = \frac{C_f}{2}$.
 
 We will prove that for all $n \in \N$ we have
 \[
  f^n(\alpha \setminus \beta'_n) \cap f^n(\beta'_n)= \varnothing \text{ and } f^n(\beta'_n) \subset [f^n(\alpha)]
 \]
 and since $\beta' \subset \beta'_n$ is a subsegment of the legal segment $\beta'_n$, the same holds for $\beta'$.
 
 The condition is true for $n=0$ since $\beta'_0 = \beta$.
 
 Suppose the lemma holds for some $n \in \N$. Since $\beta'_n$ is legal 
 we have $\len(f^n(\beta'_n)) = \lambda^n \len(\beta'_n)$. Besides $f^n(\beta'_n) \subset [f^n(\alpha)]$. Apply $f$ to the path $f^n(\alpha)$. There may be cancellation at the endpoints of the legal segment $f^n(\beta'_n)$ but this cancellation does not exceed $\BBT(f)$ when measured in $f^{n+1}(\beta'_n)$, since $\beta'_n$ is legal. 
 
 This neighbourhood in $f^{n+1}(\beta'_n)$ corresponds to a $\lambda^{-n-1}\BBT(f)$-neighbourhood of the endpoints of $\beta'_n$: as a result $f^{n+1}(\beta'_{n+1})$ does not intersect $f(f^n(\alpha) \setminus f^n(\beta'_n))$ since it is contained in $\beta'_n$. It does not intersect $f(f^n(\beta'_n) \setminus f^n(\beta'_{n+1}))$ either. Finally note that $f(f^n(\alpha) \setminus f^n(\beta'_n)) \cap f^{n+1}(\beta'_n)=f^{n+1}(\alpha \setminus \beta'_n) \cap f^{n+1}(\beta'_n)$. 
 
 Since $f^{n+1}(\beta'_{n+1})$ does not meet $f^{n+1}(\alpha \setminus \beta'_{n+1})$, is is contained in $[f^{n+1}(\alpha)]$.  
\end{proof} 

A consequence of Lemma \ref{lem:thetas-disjoints} is:
\begin{lem} \label{lem:croissance-segments-legaux} 
 Let $\phi \in \Out(G)$ be fully irreducible. Let $f : T \to T$ be a train track representative.
 For any path $\alpha$ in $T$, for any legal subpath $\beta \subset \alpha$ such that $\len(\beta) > 2C_f$, we have for every $n \in \N$:
 \[
  \len(f^n(\beta) \cap [f^n(\alpha)]) \geq \frac{1}{2} \lambda ^n \len(\beta)
 \]
\end{lem}
\begin{proof}
 Suppose $\alpha$ contains a legal subsegment $\beta$ of length greater than $2C_f$. Let $\beta'$ be the subsegment of $\beta$ at distance $\frac{C_f}{2}$ from the endpoints of $\beta$. By Lemma \ref{lem:thetas-disjoints}, for every $n \in \N$, $f^n(\beta') \subset [f^n(\alpha)]$ so we have:
 \begin{align*}
  \len(f^n(\beta')) &= \lambda^n \len(\beta') \\
     &= \lambda^n \left ( \len(\beta) - \frac{2 \BBT(f)}{\lambda - 1} \right ) \\
      & \geq  \lambda^n \len(\beta) \left(1- \frac{1}{2} \right)
 \end{align*}
 Therefore $\len(f^n(\beta) \cap [f^n(\alpha)]) \geq \frac{1}{2} \lambda^n \len(\beta)$.
\end{proof}

\begin{defi}
 Let $f: T \to T$ be a train track representative for $\phi \in \Out(G)$.
 A non-backtracking segment $\eta \subset T$ is a \emph{periodic Nielsen path} if there exists $g \in G$ and $n \geq 1$ such that $g [f^n(\eta)] = \eta$. We call it simply a Nielsen path if we can choose $n=1$.
 
 A periodic Nielsen path is \emph{indivisible} if it cannot be written as the non-backtracking concatenation $\alpha \cdot \beta$ of two periodic Nielsen paths.
\end{defi}
A result about periodic indivisible Nielsen paths (or \emph{pINPs}) from \cite{PapinDetection} is:
\begin{prop}
 Let $f: T \to T$ be a train track representative for an automorphism $\phi \in \Out(G)$. There are only finitely many orbits of periodic indivisible Nielsen paths.
\end{prop}

Periodic Nielsen paths give a characterization of pseudo-periodic conjugacy classes for $\phi \in \Out(G)$, proved in \cite[Section 3]{PapinDetection}: 
\begin{lem} \label{lem:carac-pseudo-per}
 Let $f: T \to T$ be a train track representative for $\phi \in \Out(G)$.
 The conjugacy class of an element $g \in G$ is pseudo-periodic for $\phi$ if and only if the axis of $g$ in $T$ is a geodesic concatenation of periodic indivisible Nielsen paths.
\end{lem}

\begin{lem}\label{lem:courteConcatenation}
 Suppose that $\phi$ is fully irreducible, pseudo-atoroidal and has a train track representative $f: T \to T$. Then there exists $m \in \N$ such that it is impossible to concatenate more than $m$ periodic Nielsen paths for $f$ together in $T$, and more than $m$ periodic Nielsen paths for $f_-$ in $T_-$.
\end{lem}
\begin{proof}
 There are finitely many orbits of periodic INPs; let $l$ be the number of orbits of periodic INPs. By contradiction, suppose $L$ is a path in $T$ which contains a concatenation of more than $2l$ pINPs. Then there exists $\eta, g\eta \subset L$ with $g \in G$ loxodromic such that there is a fundamental domain for $g$ which is a concatenation of pINPs. By Lemma \ref{lem:carac-pseudo-per}, this implies that $g$ is a loxodromic pseudo-periodic element for $\phi$, which is impossible since $\phi$ is pseudo-atoroidal, so we can set $m := 2l$.
\end{proof}


\subsection{The Lipschitz metric on $\D$}

The space $\D$ can be endowed with a pseudo-metric called the \emph{Lipschitz metric}. For  $T, T' \in \D$ define
\[
 \Lip(T,T')= \inf_{f: T \to T'} \Lip(f)
\]
where the infimum is taken over all $G$-equivariant Lipschitz functions $f : T \to T'$. In \cite{MeinertTheLipschitzMetric} the following is proved:
\begin{prop}
 For $T, T' \in \D$ there exists a $G$-equivariant map $f:T \to T'$, sending vertex to vertex and edge to edge path, linear on the edges, such that $\Lip(f) = \Lip(T, T')$.
\end{prop}

The Lipschitz metric is defined as follows: for $T, T' \in \D$
\[
 d_{\Lip}(T, T') = \log \left [ \Lip(T,T') \frac{\vol(T/G)}{\vol(T'/G)} \right ]
\]
The distance $d_{\Lip}(T, T')$ is unchanged by rescaling $T$ or $T'$: it only depends on their projective classes. If $T, T'$ are normalized so that $\vol(T)= \vol(T')=1$ then $d_{\Lip}(T, T') = \log \Lip(T,T')$. Sometimes it is more practical to work with $1$-Lipschitz maps, for example when $T \to T'$ is a collapse or a fold. When $\Lip(T,T')=1$ then $d_{\Lip}(T,T')= \log \frac{\vol(T)}{\vol(T')}$. 

The Lipschitz metric is not a metric in the actual sense.
\begin{lem}
 The Lipschitz metric has the following properties:
 \begin{enumerate} [label=(\roman*)]
  \item for $T, T' \in \D$, $d_{\Lip}(T,T') \geq 0$
  \item for $T, T', T'' \in \D$, $d_{\Lip}(T, T'') \leq d_{\Lip}(T, T') + d_{\Lip}(T', T'')$
 \end{enumerate}
\end{lem}
\begin{proof}
 (i) Let $T, T' \in \D$ scaled such that $\Lip(T,T') = 1$. Then $f$ induces a $1$-Lipschitz map on the quotients. By minimality $f$ is surjective so $\vol(T'/G) \leq \vol(T/G)$. Thus 
 \[
  \log \Lip(T,T') \frac{\vol(T/G)}{\vol(T'/G)} \geq 1
 \]
 
 (ii) Let $T, T', T'' \in \D$. Let $f: T \to T', f': T' \to T''$ be Lipschitz maps. Then $\Lip(T, T'') \leq \Lip(f' \circ f) \leq \Lip(f) \Lip(f')$. By taking the lower bound we get $\Lip(T, T'') \leq \Lip(T, T') \Lip(T', T'')$. By taking the logarithm we obtain what we want.
\end{proof}
\begin{rema}
 The other properties of metrics fail for $\D$:
 \begin{itemize}
  \item like in $\CV_N$ the Lipschitz metric is not symmetric. A common counter example is drawn on Figure \ref{fig:contrex-prop-metrique}: if $T, T'$ are the same tree with a different metric on edges such that $T'/G$ has a very short loop, $d_{\Lip}(T',T)$ is very big.
  \item unlike in $\CV_N$ there exist $T, T'$ in $\D$ such that $d_{\Lip}(T, T') = 0$ and $d_{\Lip}(T', T) \neq 0$. More precisely, $d_{\Lip}(T',T)$ can be chosen arbitrarily big. See Figure \ref{fig:contrex-prop-metrique}.
  \item If $\phi \in \Out(G)$ and $T \in \D$ then $d_{\Lip}(T, T \cdot \phi) = 0$ implies that $T = T\cdot \phi$. This is a consequence of the fact that the actions are minimal so maps are surjective. A $1$-Lipschitz $G$-equivariant surjective map from $T$ to itself is an isometry, unless $G$ is solvable. 
 \end{itemize}
  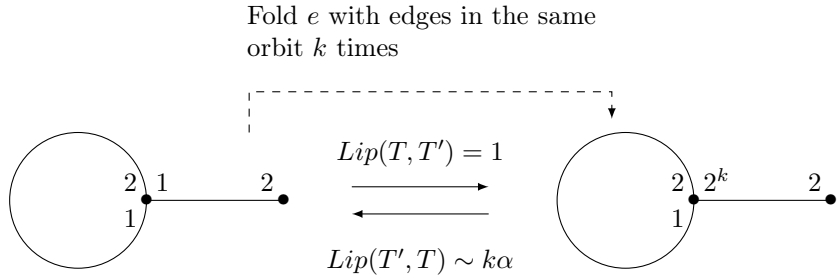
\begin{figure}
   \centering
   \begin{tikzpicture}[scale = 0.9]

\newcommand{\arrowIn}{
\tikz \draw[-stealth, very thick] (-1pt,0) -- (1pt,0);
}

\begin{scope}
    \draw (2, 4.5) node {1. Classical example for $F_2$};
     \coordinate (o) at (2,2);
     \coordinate (no) at (1.5, 3);
     \coordinate (so) at (1.5, 1);
     \coordinate (ne) at (2.5, 3);
     \coordinate (se) at (2.5, 1);

     \coordinate (a) at (0,2) ;
    \coordinate (b) at (4,2) ;
     \coordinate (a+) at (0,3);
    \coordinate (b+) at (4,3);
     \coordinate (a-) at (0,1);
    \coordinate (b-) at (4,1);

   \draw (o) node {$\bullet$} .. controls (no) and (a+).. (a) node [below left] {$a$} .. controls (a-) and (so) .. (o) node[near start, sloped]{\arrowIn};
   \draw (o) node {$\bullet$} .. controls (ne) and (b+).. (b) node [below right]{$b$} .. controls (b-) and (se) .. (o) node[near start, sloped]{\arrowIn};

    \draw (1,0.8) node {$1/2$};
   \draw  (3,0.8) node {$1/2$};

\draw (0.5,3.5) node {$T$};
  
 \draw[-latex] (5,2) -- (7,2);
\draw (6,2.5) node {$Lip(T,T') \sim 2$};

 \draw[latex-] (5,1.6) -- (7,1.6);
\draw (6,1) node {$Lip(T',T) \sim \frac{1}{2\varepsilon}$};

\end{scope}

\begin{scope}[xshift = 8cm] 
     \coordinate (o) at (3.5,2);
     \coordinate (no) at (3, 3.5);
     \coordinate (so) at (3, 0.5);
     \coordinate (ne) at (3.5, 2.5);
     \coordinate (se) at (3.5, 1.5);

     \coordinate (a) at (0,2) ;
    \coordinate (b) at (4,2) ;
     \coordinate (a+) at (0,3.5);
    \coordinate (b+) at (4,2.5);
     \coordinate (a-) at (0,0.5);
    \coordinate (b-) at (4,1.5);

   \draw (o) node {$\bullet$} .. controls (no) and (a+).. (a) node [below left] {$a$} .. controls (a-) and (so) .. (o) node[near start, sloped]{\arrowIn};
   \draw (o) node {$\bullet$} .. controls (ne) and (b+).. (b) node [below right]{$b$} .. controls (b-) and (se) .. (o) node[near start, sloped]{\arrowIn};

    \draw (1.5,0.5) node {$1- \varepsilon$};
   \draw  (3.8,1) node {$\varepsilon$};

\draw (0.5,3.5) node {$T'$};
\end{scope}

\begin{scope}[yshift=-5cm]
\draw (1.8,4) node {2. Example in $BS(2,4)$};

 \coordinate (w) at (2,0);
    \coordinate (v) at (4,0);

    \draw (v) node  {$\bullet$} node [above left] {$2$} ;
    \draw (w) node{$\bullet$} node [above right]{$1$} node [above left] {$2$} node [below left] {$1$};
    \draw  (v)--(w);
    \draw (w) arc (0: 360:1);

 \draw[-latex, dashed] (3.5,1)|- (6,1.6)  -| (8.8,1.2);

\draw (6.25, 2.5) node[text width=5cm] {Fold $e$ with edges in the same orbit $k$ times};

 \draw[-latex] (5,0.2) -- (7,0.2);
\draw (6,0.7) node {$Lip(T,T') =1$};

 \draw[latex-] (5,-0.2) -- (7,-0.2);
\draw (6,-0.9) node {$Lip(T',T) \sim k\alpha$};

\end{scope}

\begin{scope}[yshift= -5cm, xshift=8cm]
     \coordinate (w) at (2,0);
    \coordinate (v) at (4,0);

    \draw (v) node  {$\bullet$} node [above left] {$2$} ;
    \draw (w) node{$\bullet$} node [above right]{$2^k$} node [above left] {$2$} node [below left] {$1$};
    \draw  (v)--(w);
    \draw (w) arc (0: 360:1);
\end{scope}

\end{tikzpicture} 
   \caption{Counter examples for symmetry of the Lipschitz metric} \label{fig:contrex-prop-metrique}
  \end{figure}
\end{rema}

The Lipschitz metric and its computation have been explored before in \cite{BestvinaBersLike} for free groups, \cite{FrancavigliaMartino} for free products, and \cite{MeinertTheLipschitzMetric} for more general deformation spaces. The facts presented below can be found in these papers.

Let $T, T' \in \D$. Let $f: T \to T'$ be a piecewise linear $G$-equivariant map. The \emph{tension graph} $\Delta(f)$ is the subforest of $T$ spanned by edges $e \in E(T)$ such that the stretch factor on $e$ is $\Lip(f)$. The map $f$ is \emph{optimal} if it realizes the infimum of $\Lip(T, T')$ and if at every vertex $v \in \Delta(f)$, there are at least two gates at $v$ for the gate structure induced by $f$ which contain edges in $\Delta(f)$. In \cite{MeinertTheLipschitzMetric} Meinert proves that optimal maps exist. 

The distance between two points in $\D$ can be effectively computed by comparing translation lengths of some elements of $G$ in both trees.
Suppose $f : T \to T'$ is a $G$-equivariant map between trees of $\D$. For every $g \in G$ we have, by applying $f$ to a fundamental domain,
\[
 \frac{\|g\|_{T'}}{\|g\|_T} \leq \Lip(f)
\]
thus by taking the lower bound we have $\Lip(T,T') \geq \max_{g \in G} \frac{\|g\|_{T'}}{\|g\|_T} $. It is actually an equality, as this result from \cite{MeinertTheLipschitzMetric} states:
\begin{lem} \label{lem:existe-g-legal} 
 Let $T, T' \in \D$. Let $f : T \to T'$ be an optimal map. 
 There exists $g \in G$ such that $\axe_T(g)$ is $f$-legal and contained in the tension graph for $f$. In particular
 \[
  \Lip(T, T') = \frac{\|g\|_T'}{\|g\|_T} = \max_{h \in G} \frac{\|h\|_T'}{\|h\|_T}
 \]
\end{lem}

\begin{figure}
 \centering
 \begin{tikzpicture}[scale=1]

 \draw (-2,2) ellipse [x radius=1,y radius=1.2];
\draw (-2, 0) node {Loop};

     \coordinate (o) at (2,2);
     \coordinate (no) at (1.5, 3);
     \coordinate (so) at (1.5, 1);
     \coordinate (ne) at (2.5, 3);
     \coordinate (se) at (2.5, 1);

     \coordinate (a) at (0,2) ;
    \coordinate (b) at (4,2) ;
     \coordinate (a+) at (0,3);
    \coordinate (b+) at (4,3);
     \coordinate (a-) at (0,1);
    \coordinate (b-) at (4,1);

   \draw (o)  .. controls (no) and (a+).. (a)  .. controls (a-) and (so) .. (o) ;
   \draw (o)  .. controls (ne) and (b+).. (b)  .. controls (b-) and (se) .. (o);
 
 \draw (2,0) node {Figure eight};


\draw (6,2) arc (0:360:0.5) -- (7,2)  arc (180: 540:0.8);
\draw (6.6, 0) node {Barbell};

\draw (-2,-2) node {$\bullet$} -- (0.5,-2)  arc (180: 540:0.8);
\draw (0, -3.5) node {Singly degenerate barbell};

\draw (5,-2) node {$\bullet$} -- (7.5,-2)  node {$\bullet$};
\draw (6, -3.5) node {Doubly degenerate barbell};

\end{tikzpicture} 
 \caption{The five possible shapes for a candidate in the quotient graph} \label{fig:candidats}
\end{figure}
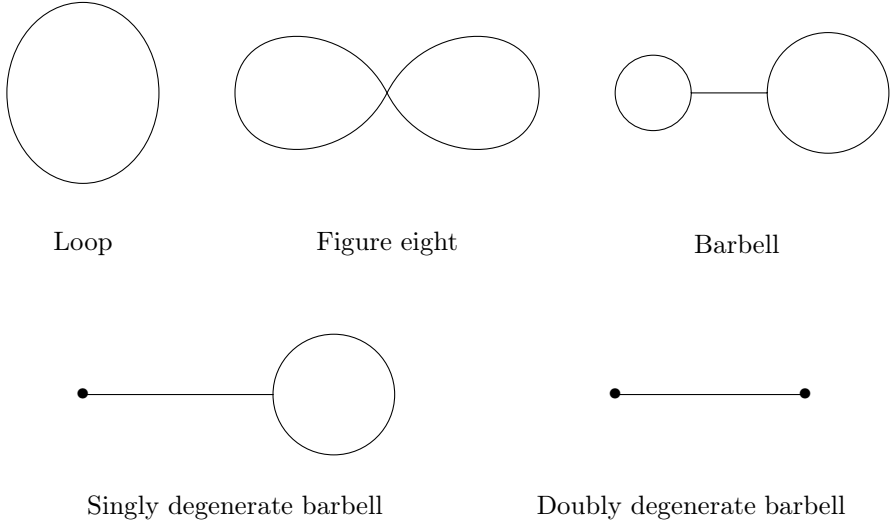

Let $T \in \D$. A \emph{candidate} of $T$ is an element $g \in G$ such that the map $\pi: \axe_T(g) / \langle g \rangle \to T/G$ has one of the following forms (see Figure \ref{fig:candidats}):
\begin{itemize}
 \item a loop: the map $\pi$ is an embedding
 \item a figure eight: there are two embedded circles $u,v$ in $T$ which intersect in exactly one point. The map $\pi$ maps the circle $\axe_T(g)/\langle g \rangle$ to the tight loop which crosses $u$ and $v$ successfully.
 \item a barbell: there are two disjoint embedded circles $u, v$ in $T/G$, and a segment $s$ which connects $u$ to $v$; $\pi$ maps the circle to the tight loop which crosses $u, s$, then $v$, and then $s$ backwards.
 \item an embedded singly degenerate barbell : it is the degenerate case of the barbell where $v$ is a single vertex. In that case, the vertex group at $v$ must be greater than the edge group of the last edge of $s$.
 \item an embedded doubly degenerate barbell : degenerate case of the barbell where both circles are single vertices. The vertex group at $u$ must also be greater than the edge group of the first edge of $g$.
\end{itemize} 

In particular a candidate crosses each orbit of edge at most twice.
In \cite[Theorem 9.10]{FrancavigliaMartino} the following theorem is proved for the case of free products: 
\begin{theo} \label{theo:distance-lipschitz-candidats} 
 Let $T \in \D$. For every $T' \in \D$, there exists a candidate $g$ of $T$ such that
 \[
  \Lip(T, T') =  \frac{\|g\|_T'}{\|g\|_T}
 \]
\end{theo}
For the proof we refer to \cite{FrancavigliaMartino}. Note that the context differs a little since the deformation space and the group are different. The relevant point is that contrary to elements of $\CV_N$, trees in the outer space for a free product may have non-free vertex stabilizers, which account for the degenerate barbells. In the case of GBS products we also have vertices whose stabilizer is greater than the stabilizers of incident edges, hence the presence of degenerate barbell candidates.

\subsection{The axis of an irreducible automorphism}

\begin{prop}
For an automorphism $\phi \in \Out(G)$ with a primitive train track map $f : T \to T$, then 
\[
\Lip(f) = \Lip(T, T \cdot \phi) = \min_{S \in \D} \Lip(S, S \cdot \phi)  
\]
\end{prop}
\begin{proof}
 Since $f$ is a train track map, at each vertex of $T$ there are at least two gates. Consequently there exists $h \in G$ such that $\axe_T(g)$ is legal.
 Let $n \in \N$. Since $h$ is $f^n$-legal for all $n \in \N$ we have $\frac{\|\phi^n(h)\|_{T}}{\|h\|_T} = \Lip(f^n) = \Lip(f)^n$.
 
 Then we have
 $\Lip(T, T \cdot \phi^n) = \sup_{g \in G} \frac{\|\phi^n(g)\|_{T}}{\|g\|_T}  = \frac{\|\phi^n(h)\|_{T}}{\|h\|_T} = \Lip(f)^n$.  Let $\lambda := \Lip(f)$.

 Let $S \in \D$. Let $f': S \to S \cdot \phi$ be such that $\Lip(S, S \cdot \phi)= \Lip(f')=: \lambda'$. By triangular inequality, for $n \in \N$ we have $\Lip(S,S \cdot \phi^n) \leq \lambda'^n$. 
 
  By triangular inequality we have
 \[
  \lambda^n = \Lip(T, T \cdot \phi^n) \leq \Lip(T, S) \Lip(S, S \cdot \phi^n) \Lip(S,T) \leq \lambda'^n \Lip(T,S) \Lip(S,T)
 \]
 If $\lambda' < \lambda$ this inequality becomes false when $n$ is big enough, hence the minimality of $\lambda= \Lip(T, T \cdot \phi)$.
\end{proof}

A map $\gamma:\R \to \D$ is a \emph{geodesic} for the Lipschitz distance if for any $t,s \in \R$, $t<s \Rightarrow d_{\Lip}(\gamma(t), \gamma(s))= s-t$. Since the metric is not symmetric, the distance $d_{\Lip}(\gamma(s), \gamma(t))$ needs not be $|s-t|$: in fact it can even be zero.

\begin{prop} \label{prop:existe-geod-continue}
 Let $f:T \to T$ be a train track representative for $\phi$. There exists a geodesic $\L_f: \R \to \D$ such that for every $n \in \Z$, $T \cdot \phi^n \in \L_f$, and such that the map $\L_f$ is continuous for the axes topology.
\end{prop}
For a construction see \cite[Section 4.4]{MeinertTheLipschitzMetric}. The \emph{axes topology} on $\D$ is the coarsest topology such that the functions $T \mapsto \|g\|_T$ for $g \in G$ are continuous; for more information on the topologies of $\D$ see \cite{GuirardelLevitt07}.

 If $f : T \to T$ is a train track representative we choose an arbitrary axis $\L_f$ crossing $T$ and denote by $T_t$ the unique point of the axis such that $d_{\Lip}(T, T_t) = t$ if $t \geq 0$ and $d_{\Lip}(T_t, T) = -t$ if $t \leq 0$.

 Similarly there is an axis $\L_{f_-}$ for $\phi^{-1}$ defined from $f_-: T_- \to T_-$. Since $d_{\Lip}$ is $\Out(G)$-invariant the axes stay within a bidirectional bounded neighbourhood of each other.

\bigskip

We already stated that between any two trees in $\D$, there exists a $G$-equivariant quasi-isometry. The quasi-isometry constants can be chosen uniformly if the trees lie in a subsegment of $\L_f$:
\begin{lem} \label{lem:borne-kc}
 For any $T, S \in \D$ and $D \geq 1$ such that $\max \{\Lip(T,S), \Lip(S,T) \} \leq D$, any optimal map $T \to S$ is an equivariant $(D, 4D)$-quasi-isometry.
\end{lem}
\begin{proof}
 Let $T , S \in \D$. 
 Scale $T,S$ such that $\vol(T/G)=\vol(S/G)=1$. There exists $X \subset T$ such that $G \cdot X = T$ and $\diam(X) \leq 1$. 
 For every oriented edges $\overrightarrow{e},   \overrightarrow{e'} \in E(T)$, there exists an edge path with length at most $2$ with first edge $\overrightarrow{e}$ and last edge $h \overrightarrow{e'}$ for some $h \in G$. This fails if $G$ is solvable, but we assumed that it is not the case. 
 
 Let us prove this fact. First we will prove that for any edge $e$ there exists a path $\overrightarrow{e} \dots \overleftarrow{he}$ for some $h \in G$.
 
 Since $G$ is not solvable, then the action of $G$ on $T$ is irreducible.
 Since $T$ is not a line, there exist translates $h_1e, h_2e$ such that for any line containing both edges, the orientations of the edges along the line differ.

 There exists $h \in G$ such that $\axe_T(h)$ crosses both $h_1e$ and $h_2 e$. Either $h_1 e$ and $h_2 e$ point towards each other, or they may point away from each other, but in that case there exists $k \in \Z$ such that $h^k h_1e$ and $h_2e$ point towards each other. Then we get a path $\overrightarrow{h_1 e} \dots \overleftarrow{h_2e}$. 
 
 Now let $e, e'$ be edges in the tree. There is a path which connects both edges, but the path may fail to contain a translate of $\overrightarrow{e} \dots \overrightarrow{he'}$ for some  $h \in G$. By concatenating paths which reverse the orientation on one or both sides we obtain a path satisfying the condition.
 
 As for the bound on the length of the path, observe that if an edge appears in the path twice with same orientation, then subpath can be deleted to obtain a shorter path. 
 
 \bigskip
 
 Let $\tau: T \to S$ be an optimal map.
 Let $x, y \in T$. 
 Let $e$ be the first edge of $[x, y]$. Let $e'$ be an edge with origin $y$ and not in $[x,y]$: such an edge exists since $T$ is minimal. By the fact above, there exists $g \in G$ and a path containing $\overrightarrow{e'}, g\overrightarrow{e'}$ with these orientations and length at most $2$. Thus $d_T(gx, y) \leq 2$ and $[x,y] \subset \axe_T(g)$. We have $d_T(x, y) + 2 \leq \|g\|_T$.

 Then by Lemma \ref{lem:distance-entre-axes} we have $\Lip(S,T) \leq D, \Lip(T,S) \leq D$ so
 \[
  d_S(\tau(x), \tau(y)) \leq D d_T(x, y)
 \]
 and
 \begin{align*}
  d_S(\tau(x), \tau(y)) &\geq d_S(\tau(x), \tau(gx)) - d_S(\tau(gx), \tau(y)) \\
     & \geq \|g\|_S - 2\Lip(\tau)\\
     & \geq \frac{\|g\|_T}{D} -2\Lip(\tau) \\
     &\geq \frac{d_T(x,y)}{D} - \frac{2}{D} -2D
 \end{align*}
Then the optimal map $\tau : T \to T_{t'}$ is a 
$(D, 4D)$-quasi-isometry. 
\end{proof}

 We have the following result about the axes $\L_f$ and $\L_{f_-}$:
\begin{lem} \label{lem:distance-entre-axes}
  Let $\phi \in \Out(G)$ be a fully irreducible automorphism. Let $f: T \to T$ be a train track representative for $\phi$ and $f_-: T_- \to T_-$ be a train track representative for $\phi^{-1}$. Let $\L_f, \L_{f_-}$ be axes in $\D$ for $\phi$ and $\phi^{-1}$.
  
  Let $a< b$, $c,d$ be real numbers. There exists a constant $D_{a,b,c,d} >1$ such that for every $X, Y \in \{T_t/ a \leq t \leq b\} \cup \{ (T_-)_t /c \leq t \leq d\}$, $\Lip(X,Y) \leq D_{c,d}$.
\end{lem}
\begin{proof}
 Let $n,m \in \N$ be such that 
 \[
  n \log(\lambda) \leq a, m \log(\lambda) \geq b, n\log(\lambda_-) \leq c, m \log (\lambda_-) \geq d
 \]

 There is a quadrilateron which crosses $T \cdot \phi^n, T \cdot \phi^m, T_- \cdot \phi^m, T_- \cdot \phi^n$, which contains $\{T_t/ a \leq t \leq b\}$ and $\{ (T_-)_t /c \leq t \leq d\}$. Its length is
 \[
  d:= d_{\Lip}(T \cdot \phi^n, T \cdot \phi^m) + d_{\Lip}(T, T_-)+ d_{\Lip}(T_- \cdot \phi^m, T_- \cdot \phi^n) +  d_{\Lip}(T_-,T)
 \]
 Therefore, for every $X, Y$ as above we have $L(X,Y) \leq D:= e^d$.
\end{proof}

\begin{rema} \label{rema:borne-kc}
 Lemmas \ref{lem:borne-kc} and \ref{lem:distance-entre-axes} imply that for every $S \in \D$, there exist $(K, C)$ such that for every $t \in [0, \log(\lambda)]$, there exist equivariant $(K,C)$-quasi-isometries $T_t \to S$ and $S \to T_t$.
\end{rema}


\section{The stable and unstable laminations} \label{sec:axes-laminations}

Let $T,T' \in \D$. There exists a $G$-invariant quasi-isometry $T \to T'$. In fact all equivariant quasi-isometries $T \to T'$ are close:
\begin{lem} \label{lem:applications-proches} 
 Let $T,T'$ be metric $G$-trees such that $T$ is co-compact. Let $u,v$ be continuous $G$-equivariant maps $T \to T'$. There exists a constant $C$ depending on $u$ and $v$ such that for every $x \in T$
 \[
  d_{T'}(u(x), v(x)) \leq C
 \] 
\end{lem}
\begin{proof}
 Let $K \subset T$ be a compact subset such that $G \cdot K = T$. Let $C:= \max_{x \in K} d_{T'}(u(x), v(x))$. For every $y \in T$ there exists $g \in G$ and $x \in K$ such that $y = gx$ so by $G$-equivariance 
 \[
  d_{T'} (u(y),v(y)) = d_{T'} (gu(x),gv(x)) =d_{T'}(u(x), v(x)) \leq C.
 \]
\end{proof}
Recall that a $G$-invariant quasi-isometry $f$ induces a $G$-equivariant homeomorphism $\partial T \to \partial T'$. Because of Lemma \ref{lem:applications-proches} the homeomorphism does not depend on $f$ so there is a canonical $G$-invariant identification of the boundaries of all trees of $\D$.

\medskip

A \emph{lamination} $\Lambda$ is a $G$-invariant, symmetric, closed subset of $\partial T \times \partial T \setminus \Delta$ where $\Delta$ is the diagonal, for some $T \in \D$. The discussion above implies that 
for any $S \in \D$, $\Lambda$ can be canonically identified with a subset of $\partial S \times \partial S \setminus \Delta$ so we may drop the reference to $T$.

When we fix a tree $T$, $\Lambda$ identifies with a $G$-invariant set of unoriented bi-infinite geodesics of $T$ which we call the \emph{realization} of $\Lambda$ in $T$ and which we denote by $\Lambda_T$.
Its elements are called \emph{leaves}. A  \emph{leaf segment} is a subsegment of a leaf of $\Lambda_T$. The assumption that $\Lambda$ is a closed subset of $\partial T \times \partial T \setminus \Delta$ translates into the following fact: if $(\sigma_n)_{n \in \N}$ is an increasing sequence of leaf segments in $T$ whose union is a bi-infinite geodesic $\ell \subset T$, then $\ell$ is a leaf of $\Lambda_T$.

If $T' \in \D$ is another tree, there exists a $G$-invariant quasi-isometry $h: T \to T'$. For any leaf $\lambda \in \Lambda_T$, the line of $T'$ obtained by tightening $h(\lambda)$ is a leaf of $\Lambda_{T'}$, and conversely all leaves of $\Lambda_{T'}$ are tightened images of leaves of $\Lambda_T$.

\bigskip

Let $\phi \in \Out(G)$ be a fully irreducible automorphism. Let $f: T \to T$ be an primitive train track representative for $\phi$ with Lipschitz constant $\lambda > 1$. Define the \emph{stable lamination} $\Lambda^+_f$ by its realization in the train track tree $T$, as the set of bi-infinite geodesics whose subsegments belong to 
\[
 \{\sigma \subset T / \exists e \in E(T), \exists n \in \N, \: \sigma \subset f^n(e)\}
\]

We call $\Lambda^+_{f_-}$ the \emph{unstable lamination}.

\begin{rema}
 Since $\lambda > 1$, if $e$ is an edge of $T$ such that $e \subset  \mathring{f(e)}$, then the limit of $f^n(e)$ when $n \to \infty$ is a leaf of the lamination. Since $f$ is primitive, one can check that the set of leaves which can be obtained by this process by replacing $f$ with $g f$ for $g \in G$ is a $G$-invariant subset of the stable lamination and its closure is the stable lamination.
\end{rema}

\begin{lem}\label{lem:longue-feuille}
 For every $l > 0$ there exists $n_\Lambda > 0$ such that if $\alpha \subset T$ contains a legal subsegment with length at least $2C_f$, then $[f^n(\alpha)]$ contains a leaf segment of $\Lambda^+_f$ with length at least $l$ for all $n \geq n_{\Lambda}$.
\end{lem}
\begin{proof}
Let $\beta$ be a legal subsegment of $\alpha$ with length at least $2C_f$. By Lemma \ref{lem:thetas-disjoints} the subpath $\theta \subset \beta$ obtained by truncating the $C_f/2$-neighbourhood of the endpoints has the following property: for any $n \in \N$, $f^n(\theta) \subset [f^n(\alpha)]$.

We have $\len(\theta) \geq C_f$. There exists $n_1 \in \N$ depending on $T$ and $\lambda$ such that $\lambda^{n_1} C_f \geq 2 \max_{e \in E(T)} \len(e)$ so $f^{n_1}(\theta)$ contains an edge.

There exists $n_2 \in \N$ depending on $l$ such that $\lambda^{n_2} \min_{e \in E(T)} \len(e) \geq l$. For $n \geq n_2$, for any $e \in E(T)$, $f^n(e)$ contains a leaf segment of $\Lambda^+_f$ with length $l$.

Then for any $n \geq n_1 + n_2$, $[f^n(\alpha)]$ contains a leaf segment of $\Lambda^+_f$ with length $l$.
\end{proof}

\begin{defi}
 A lamination $\Lambda$ is \emph{minimal} in $T$ if all leaves of $\Lambda_T$ have the same leaf segments up to the action of $G$, i.e. for every leaf segment $\sigma \subset T$ of $\Lambda_T$, for every leaf $\ell$ in $T$, there exists $g \in G$ such that $g \sigma \subset \ell$.
\end{defi}

\begin{defi}
 Let $\Lambda$ be a lamination. Let $T \in \D$. A leaf $\ell \in \Lambda_T$ is \emph{quasi-periodic} if for every $C > 0$ there exists $L >0$ such that for every subsegment $\sigma \subset \ell$ with $\len(\sigma) = C$, for every subsegment $\gamma \subset \ell$ with $\len(\gamma) > L$, there exists $g \in G$ such that $g \sigma \subset \gamma$.
\end{defi}

\begin{rema}
 If $T \to S$ is a quasi-isometry and $\Lambda_T$ is minimal, then $\Lambda$ is minimal in $S$. Similarly, if a leaf in $T$ is quasi-periodic, then the realization of this leaf in another tree $S \in \D$ is also quasi-periodic. A proof is given in \cite[Remark 1.17]{PapinDetection}.
\end{rema}

The following is proved in \cite{BestvinaFeighnHandelLaminations}, although for a slightly different definition of the stable lamination.
\begin{lem}
 Let $f: T \to T$ be an irreducible train track representative for an automorphism $\phi \in \Out(G)$ with Lipschitz constant $\lambda$. Then the stable lamination $\Lambda^+_f$ is minimal and its leaves in $T$ are quasi-periodic.
\end{lem}
\begin{proof}
 First let us prove the minimality.
 There exists $N \in \N$ such that for all $e,e' \in E(T/G)$,  $f^N(e)$ contains an edge in the orbit of $e'$. 
 
 Let $\ell\subset T$ be a leaf of $\Lambda^+_f$.  Let $k \in \N$. Let us prove that there exists a constant $L_k$ such that every segment of $\ell$ longer than $L_k$ contains a translate of $f^k(e_0)$ for every $e_0 \in E(T/G)$.
 
 Let $\sigma$ be a subsegment of $\ell$ with length at least $L_k :=2\lambda^k \max_{e \in E(T/G)} \len (f^N(e))$. By definition of the stable lamination, there exists $e \in E(T)$ such that $\sigma \subset f^{N+n}(e)$ with $n \geq k$. The segment $\sigma$ is contained in the concatenation of segments $f^k(f^N(e'))$ for edges $e' \subset f^{n-k}(e)$. Each of these segments contains a translate of $f^k(e_0)$ and is shorter than $\len(\sigma)/2$ so one of them is contained in $\sigma$, hence in $\ell$. This proves that $\ell$ contains a translate of every leaf segment of $\Lambda^+_f$ contained in $f^k(e_0)$ for any $e_0 \in E(T)$.
 
 This also proves the quasi-periodicity: for all $C > 0$ there exists $k \in \N$ such that every segment $I$  longer than $C$ is contained in $f^k(e)$ for some edge $e$, and every leaf segment of the leaf $\ell$ longer than $L_k$ contains a copy of $f^k(e)$ and thus a copy of $I$.
\end{proof}

\begin{lem} \label{lem:laminations-memes-segments} 
 Let $\Lambda, \Lambda'$ be two distinct minimal closed $G$-invariant laminations with quasi-periodic leaves. Then for any $G$-tree $T$, there exists a bound $C_T$ on the length of leaf segments which are common to both laminations.
\end{lem}
\begin{proof}
 Let $T \in \D$.
 By contraposition we prove that if the bound $C_T$ does not exist then $\Lambda$ and $\Lambda'$ have the same sets of leaf segments in $T$. Since they are closed this implies $\Lambda= \Lambda'$.
 
 Suppose that there exist leaf segments $(\eta_n)_{n \in \N}$ with $\eta_n \subset \Lambda_T \cap \Lambda'_T$ for all $n \in \N$, and $\len(\eta_n) \to \infty$. Then for any leaf segment $\sigma \subset \Lambda_T$, there exists $n \in \N$ and $g \in G$ such that $g \sigma \subset \eta_n \subset \Lambda'_T$ so $\sigma$ is also a leaf segment of $\Lambda'_T$. By symmetry $\Lambda_T$ and $\Lambda'_T$ have the same leaf segments, so the laminations $\Lambda$ and $\Lambda'$ are equal.
\end{proof}

We need the following:
\begin{lem}\label{lem:qi-et-feuille2} 
 Let $h : T \to T'$ be a $(K,C)$-quasi-isometry where $K \geq 1, C \geq 0$.
 
 For every $l > 0$ there exists $L > 0$ depending on $l, K, C$ such that if $\eta$ is a bi-infinite geodesic, if $\sigma$ is a subsegment of $\eta$ with length at least $L$, then $[h(\sigma)]$ contains a subsegment of $[h(\eta)]$ with length at least $\ell$.
\end{lem}
\begin{proof}
 Let $l >0$. Let $\eta$ be a bi-infinite geodesic.
 Let $\sigma$ be a subsegment of $\eta$. Then the length of $[h(\sigma)]$ is at least $ K^{-1} \len(\sigma)-C$.
 
 The image $h(\eta)$ lies in a $\BBT(h)$-neighbourhood of $[h(\eta)]$, where $\BBT(h) \leq K^2C + C$. The endpoints of $[h(\sigma)]$ are in this neighbourhood so the length of $[h(\sigma)] \cap [h(\eta)]$ is at least $K^{-1} \len(\sigma) - C -2 \BBT(h)$. Thus by taking $L \geq K(l + C + 2\BBT(h))$, if $\len(\sigma) \geq L$, then $[h(\sigma)]$ contains an $l$-segment of $[h(\eta)]$.
\end{proof}

For a lamination $\Lambda$, $T \in \D$ and $C >0$, a \emph{$C$-piece} of $\Lambda_T$ is a leaf segment of $\Lambda_T$ with length $C$. 

\begin{lem} \label{lem:laminations-egales} 
 Let $f: T \to T$ and $f': T' \to T'$ be two train track representatives for a fully irreducible automorphism $\phi \in \Out(G)$. Then the stable laminations $\Lambda^+_f$ and $\Lambda^+_{f'}$ are equal.
\end{lem}
\begin{proof}
 We will prove that every leaf segment of $(\Lambda^+_f)_{T'}$ is also a leaf segment of $(\Lambda^+_{f'})_{T'}$: by symmetry we will get the result.
 
 Let $C>0$. By quasi-periodicity of the leaves of $\Lambda^+_f$, there exists $L>0$ such that every leaf segment of $(\Lambda^+_f)_{T'}$  longer than $L$ contains every orbit of leaf segment of $(\Lambda^+_f)_{T'}$ with length at most $C$.

 Let $h: T \to T'$ be a $G$-equivariant quasi-isometry. 

 By Lemma \ref{lem:qi-et-feuille2} there exists $L_0 > 0$ such that for every bi-infinite geodesic $\eta \subset T$, for every segment $\sigma \subset \eta$ with length at least $L_0$, the segment $[h(\sigma)] \cap [h(\eta)]$ has length at least $2L + 2\BBT(h)$. Without loss of generality, we may assume $L_0 \geq C_f$.

 There exists $g \in G$ be such that $\axe_{T'}(g)$ is legal for $f'$. The conjugacy class of $g$ is not pseudo-periodic since $\|\phi^n(g)\|_{T'} \to \infty$ when $n \to \infty$.

 The axis of $g$ in $T$ does not have to be legal, however the number of orbits of $f$-illegal turns under the action of $\langle \phi^n(g) \rangle$ cannot increase when $n \to \infty$: $f$ sends $f$-legal subsegments to $f$-legal subsegments and $\axe_{T}(\phi^n(g)) = [f^n(\axe_T(g))]$. Since $\{\|\phi^n(g)\|_{T}, n \in \N\}$ is unbounded, this implies that there exists $n_0 \in \N$ such that $\axe_T(\phi^{n_0}(g))$ contains an $f$-legal subsegment with length $L_0$. Since $\axe_{T'}(\phi^{n_0}(g))$ is also $f'$-legal, up to replacing $g \in \phi^{n_0}(g)$ we can assume $n_0 = 0$.
 
 As in Lemma \ref{lem:longue-feuille} there exists $N > 0$ such that for every $e \in E(T')$, $f'^N(e)$ is a leaf segment of $(\Lambda^+_{f'})_{T'}$ with length at least $2L$. The axis of $\phi^N(g)$ in $T$ still contains a legal subsegment with length at least $L_0$ since $L_0 \geq C_f$. Once again, up to replacing $g$ by $\phi^N(g)$, we may assume that the axis of $g$ in $T'$ can be cut into pieces of $(\Lambda^+_{f'})_{T'}$ with length at least $2L$.

 The map $h$ maps $\axe_T(g)$ to a $\BBT(h)$-neighbourhood of $\axe_{T'}(g)$.  
 By Lemma \ref{lem:qi-et-feuille2} and by definition of $L_0$, $h(\axe_T(g)) \subset T'$ contains a leaf segment of $[h(\Lambda^+_f)_{T})]=(\Lambda^+_f)_{T'}$ with length at least $2L + 2\BBT(h)$ and a subsegment with length at least $2L$ is contained in $\axe_{T'}(g)$.
 
 Then there exists a segment $\gamma' \subset T'$ of length greater than $L$ that is both a leaf segment $(\Lambda^+_f)_{T'}$ and $(\Lambda^+_{f'})_T$. Thus every leaf segment of $(\Lambda_f)_{T'}$ of length $C$ is a leaf segment of $(\Lambda_{f'})_{T'}$.
\end{proof}

From now on, we can simply refer to the stable lamination as $\Lambda^+_\phi$ or simply $\Lambda^+$ when the automorphism is obvious. The notation $\Lambda^-_\phi$ denotes the unstable lamination $\Lambda^+_{\phi^{-1}}$.

\begin{lem}\label{lem:laminations-distinctes} 
 For a fully irreducible automorphism $\phi \in \Out(G)$, the stable lamination and unstable lamination are distinct.
\end{lem}
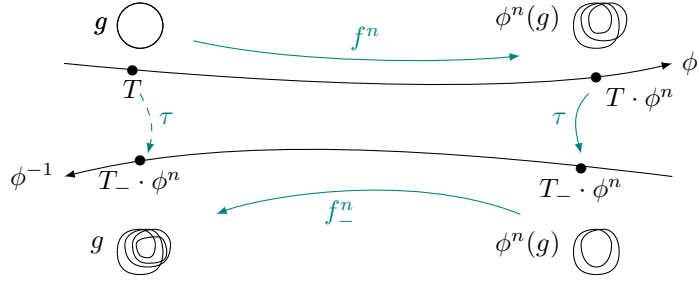
\begin{figure}
 \centering
   
\begin{tikzpicture}
  \draw[-latex] (0,0.5) .. controls (2, 0.3) and (6, 0) .. (8,0.5) node [right] {$\phi$}; 
  \draw[latex-] (0,-1) node [left] {$\phi^{-1}$} .. controls (1.8, -0.55) and (4, -0.5) .. (8,-1); 

 \draw  (0.9, 0.4) node {$\bullet$} node [below] {$T$};
\draw  (7, 0.3) node {$\bullet$} node [below right] {$T \cdot \phi^n$};
\draw  (1, -0.8) node {$\bullet$} node [below] {$T_- \cdot \phi^n$};  
\draw  (6.8, -0.9) node {$\bullet$} node [below] {$T_- \cdot \phi^n$};  

\draw [-latex, teal, dashed] (1, 0.1) .. controls (1.1, -0.1) and (1.2, -0.3) .. (1.1, -0.7) node [midway, right] {$\tau$};
\draw [-latex, teal] (6.9, 0.1) .. controls (6.7, -0.1) and (6.7, -0.3) .. (6.8, -0.7) node [midway, left] {$\tau$};

\draw [-latex, teal] (1.7, 0.8) .. controls (3, 0.5) and (5, 0.5) .. (6, 0.6) node [midway, above] {$f^n$};
\draw [latex-, teal] (2, -1.5) .. controls (2.5, -1.3) and (4.5, -0.9) .. (6, -1.5) node [midway, below] {$f_-^n$};

 \draw (1,1) circle [radius=0.3cm] node [left = 0.3cm] {$g$};
\draw (7, 0.7) .. controls (6.8, 0.7) and (6.7, 0.8) .. (6.7, 1);
\draw (7, 1.3) .. controls (6.8, 1.3) and (6.7, 1.2) .. (6.7, 1);
\draw (7, 0.7) .. controls (7.2, 0.7) and (7.3, 0.8) .. (7.3, 1);
\draw (7.1, 1.3) .. controls (7.2, 1.3) and (7.3, 1.2) .. (7.3, 1);

\draw (7.1, 1.3) .. controls (6.9, 1.3) and (6.8, 1.2) .. (6.8, 1.1) node [left = 0.15cm] {$\phi^n(g)$};
\draw (7.1, 0.8) .. controls (6.9, 0.8) and (6.8, 0.9) .. (6.8, 1.1);
\draw (7.1, 0.8) .. controls (7.3, 0.8) and (7.4, 0.9) .. (7.4, 1.1);
\draw (7.2, 1.3) .. controls (7.3, 1.3) and (7.4, 1.2) .. (7.4, 1.1);

\draw (7.2, 1.3) .. controls (7, 1.3) and (6.9, 1.2) .. (6.9, 1.1);
\draw (7, 0.9) .. controls (7, 0.9) and (6.9, 0.9) .. (6.9, 1.1);
\draw (7, 0.9) .. controls (7.1, 0.9) and (7.2, 0.9) .. (7.2, 1.1);
\draw (7, 1.3) .. controls (7.1, 1.3) and (7.2, 1.2) .. (7.2, 1.1);

 \draw (1,1) circle [radius=0.3cm] node [left = 0.3cm] {$g$};
\draw (7, -2.3) .. controls (6.8, -2.3) and (6.7, -2.2) .. (6.7, -2);
\draw (7, -1.7) .. controls (6.8, -1.7) and (6.7, -1.8) .. (6.7, -2);
\draw (7, -2.3) .. controls (7.2, -2.3) and (7.3, -2.2) .. (7.3, -2);
\draw (7.1, -1.7) .. controls (7.2, -1.7) and (7.3, -1.8) .. (7.3, -2);

\draw (7.1, -1.7) .. controls (6.9, -1.7) and (6.8, -1.8) .. (6.8, -1.9) node [left = 0.15cm] {$\phi^n(g)$};
\draw (7, -2.2) .. controls (6.9, -2.2) and (6.8, -2.1) .. (6.8, -1.9);
\draw (7, -2.2) .. controls (7.1, -2.2) and (7.2, -2.1) .. (7.2, -1.9);
\draw (7, -1.7) .. controls (7.1, -1.7) and (7.2, -1.8) .. (7.2, -1.9);


\draw (1,-2.3) .. controls (0.8, -2.3) and (0.7, -2.2) .. (0.7, -2);
\draw (1, -1.7) .. controls (0.8, -1.7) and (0.7, -1.8) .. (0.7, -2);
\draw (1, -2.3) .. controls (1.2, -2.3) and (1.3, -2.2) .. (1.3, -2);
\draw (1.1, -1.7) .. controls (1.2, -1.7) and (1.3, -1.8) .. (1.3, -2);

\draw (1.1, -1.7) .. controls (0.9, -1.7) and (0.8, -1.8) .. (0.8, -1.9) node [left = 0.15cm] {$g$};
\draw (1.1, -2.2) .. controls (0.9, -2.2) and (0.8, -2.1) .. (0.8, -1.9);
\draw (1.1, -2.2) .. controls (1.3, -2.2) and (1.4, -2.1) .. (1.4, -1.9);
\draw (1.2, -1.7) .. controls (1.3, -1.7) and (1.4, -1.8) .. (1.4, -1.9);

\draw (1.2, -1.7) .. controls (1, -1.7) and (0.9, -1.8) .. (0.9, -1.9);
\draw (1.1, -2.15) .. controls (1, -2.15) and (0.9, -2.1) .. (0.9, -1.9);
\draw (1.1, -2.15) .. controls (1.3, -2.15) and (1.35, -2.1) .. (1.35, -2);
\draw (1.25, -1.8) .. controls (1.3, -1.8) and (1.35, -1.9) .. (1.35, -2);

 \draw (1.25, -1.8) .. controls (1.2, -1.8) and (0.95, -1.8) .. (0.95, -1.9);
\draw (1, -2.05) .. controls (1, -2.05) and (0.95, -2) .. (0.95, -1.9);
\draw (1, -2.05) .. controls (1.1, -2.1) and (1.2, -2.1) .. (1.2, -1.9);
 \draw (1, -1.7) .. controls (1.1, -1.7) and (1.2, -1.8) .. (1.2, -1.9);
\end{tikzpicture}
 \caption{When applying $f^n$ and then $f_-^n$, the translation length of $g$ increases if $n$ is big enough} \label{fig:laminations-distinctes}
\end{figure}

\begin{proof}
 Assume by contradiction that $\Lambda^-_{\phi} = \Lambda^+_\phi$. See Figure \ref{fig:laminations-distinctes}.
 
 Let $f: T \to T$ be a train track representative for $\phi$ and $f_- : T_- \to T_-$ be a train track representative for $\phi^{-1}$. Let $\tau : T \to T_-$  be a $G$-equivariant quasi-isometry. 
 
 There exists $g \in G$ loxodromic whose axis in $T$ is $f$-legal.
 By Lemma \ref{lem:qi-et-feuille2} there exists $L> 0$ such that for every leaf segment $\sigma$ of $(\Lambda^+)_T$ longer than $L$, $[\tau(\sigma)]$ contains a leaf segment of $(\Lambda^+)_{T_-}$ longer than $2C_{f_-}$. 
 
 Let $n_0 \in \N$ be such that for any edge $e \in T$, $\len(f^n(e)) \geq 2C_f$. Let $n \geq n_0$. Then every $e \in \axe_{T}(\phi^n(g))$ contributes a leaf segment of $(\Lambda^+)_{T_-} = (\Lambda^-)_{T_-}$ longer than $2C_f$ in $\axe_{T_-}(\phi^n(g))$. By Lemma \ref{lem:thetas-disjoints}, the images $\beta_e :=[\tau(f^n(e))]$ for different edges $e$ contain subsegments $\beta'_e$ obtained from $\beta_e$ by cutting out the $\frac{C_{f_-}}{2}$-neighbourhood of the endpoints. The subsegments $\beta'_e$ satisfy the following: for any edge $e$, for any $m \in \N$, $f_-^m(\beta'_e) \subset \axe_{T_-}(\phi^{n-m}(g))$ and for any other edge $e'$, $f_-^m(\beta'_e)  \cap f_-^m(\beta'_{e'}) = \varnothing$. 
 
 This implies that $[f_-^{n}(\axe_{T_-}(\phi^n(g)))]= \axe_{T_-}(g)$ contains two disjoint leaf segments $f_-^n(\beta'_e)$,  $gf_-^n(\beta'_e) =f_-^n(\beta'_{\phi^n(g)e})$ longer than $\lambda^n C_f$. 
 
 Then we must have $\|g\|_{T_-} \geq \lambda^n C_f$, which is a contradiction for $n$ sufficiently big.
\end{proof}

\section{Laminations and simple elements of $G$} \label{sec:whitehead-et-laminations} 

\subsection{Simple elements, simple pairs and Whitehead graphs} 

\begin{defi}
 A loxodromic element $g \in G$ is \emph{simple} if it is contained in a proper cyclic factor of $G$.
 
\medskip
  
 A pair of elements $g, h \in G$ is \emph{simple} if there exists cyclic factors $H_g, H_h$ such that $g \in H_g, h \in H_h$ and a graph $\Gamma$ of cyclic groups with $\pi_1(\Gamma) \simeq G$, with disjoint subgraphs $\Gamma_g, \Gamma_h$ such that $\pi_1(\Gamma_g) \simeq H_g, \pi_1(\Gamma_h) \simeq H_h$.
\end{defi}
In that case $H_g, H_h$ belong to a proper \emph{system of cyclic factors}, i.e. collection of conjugacy class of cyclic factors which can be simultaneously seen in some graph of groups as the fundamental groups of disjoint subgraphs.

\begin{lem}
 Suppose $b_1(G) \geq 3$. 
 For any tree $T \in \D$, the candidates in $T$ are simple.
\end{lem}
\begin{proof}
 If $b_1(G) \geq 3$ then every candidate $g$ for $T$ avoids at least an orbit of edges $G \cdot e$. Then $T \setminus G \cdot e$ is a proper subforest of $G$ which contains the axis of $g$. It defines a cyclic factor containing $g$, so $g$ is simple. 
\end{proof}

\begin{lem} \label{lem:troisieme-candidat-compatible}
 Suppose $b_1(G) \geq 3$.
 Suppose $g,h\in G$ are candidates in $T \in \D$. There exists a candidate $k \in G$ such that $\{k,g\}$ and $\{k,h\}$ are both simple. 
\end{lem}
\begin{proof}
 Since $b_1(\Gamma) \geq 3$ and $g, h$ are candidates, neither of their axes crosses every orbit of edges in $T$.
 
 Let $e_g, e_h$ be edges in $\Gamma$ such that $\pi(\axe_T(g))$ avoids $e_g$ and $\pi(\axe_T(h))$ avoids $e_h$, where $\pi: T \to \Gamma$ is the quotient map. Since $b_1(g) \geq 3$ the graph $\Gamma' :=\Gamma \setminus \{e_g, e_h\}$, which may be disconnected, has a connected component with first Betti number $b_1(\Gamma') \geq 1$. There exists an element $k \in G$ whose axis in $T$ is in a lift of $\Gamma'$. Then $\axe_T(k)$ crosses neither $e_g$ nor $e_h$ so:
 \begin{itemize}
  \item the axes of $g,k$ are in $\Gamma \setminus e_g$, so $\{k,g\}$ is simple
  \item the axes of $h,k$ are in $\Gamma \setminus e_h$, so $\{k,h\}$ is simple
 \end{itemize}
 Thus $k$ is the element of $G$ that we were looking for.
\end{proof}

\begin{defi}
 Let $\ell$ be a bi-infinite geodesic in $T \in \D$. A \emph{turn} in $\ell$ is a pair $\{e, e'\} \subset E(T)$ of distinct edges such that $o(e) = o(e')$ and $e \cup e' \subset \ell$. 
\end{defi}

\begin{defi}
 Let $\mathcal G$ be a collection of bi-infinite geodesics of some $T \in \D$. Let $v \in V(T)$. The Whitehead graph $W:=\Wh_T(\mathcal G, v)$ is the following graph:
 \begin{itemize}
  \item vertices of $W$ are edges of $T$ with origin $v$
  \item there is an edge $e-e'$ in $W$ if there exists $\ell \in \mathcal G$ and $g \in G$ such that $g \cdot\ell$ contains both $e$ and $e'$, i.e. if $\{e, e'\}$ is a turn crossed by $\ell$
 \end{itemize}
\end{defi}
\begin{rema}
 For $\mathcal G, T, v$ as in the definition, we have $\Wh_T(\mathcal G, v) = \Wh_T(G \cdot \mathcal G,v)$.
\end{rema}

\begin{ex}
 The two main examples, which we will both use in this paper, are the following.
 \begin{enumerate}
  \item $\mathcal G = \{ \axe_T(hgh^{-1}), h \in G\}$ is the collection of axes of all conjugates of some $g \in G$. In that case we write $\Wh_T(\mathcal G, v) = \Wh_T(g,v)$.
  \item let $f : T \to T$ be a train track representative for a fully irreducible automorphism and let $S \in \D$, let $\mathcal G = (\Lambda^+)_S$. Since all leaves in $\Lambda^+$ have the same subsegments in $S$, for any leaf $\ell \in (\Lambda^+)_S$, for any $v \in V(S)$, $\Wh_S(\Lambda^+, v)= \Wh_S(\ell, v)$. 
 \end{enumerate}
\end{ex}

The interest of Whitehead graphs is that they help understanding cyclic factors.
In \cite{PapinWhitehead} we prove the following theorem (Theorem 2.14):
\begin{theo} \label{theo:whitehead-axes}
 \begin{itemize}
  \item Let $g \in G$ be a loxodromic element. Then $g$ is simple if and only if for every $T \in \D$ there exists $v \in T$ such that $\Wh_T(g,v)$ is disconnected or has a cut vertex, i.e. a vertex $p \in \Wh_T(g,v)$ such that $\Wh_T(g,v) \setminus \{p\}$ is disconnected.
  \item Let $g, h \in G$ be loxodromic elements. Then $\{g,h\}$ is simple if and only if for every $T \in \D$ there exists $v \in T$ such that $\Wh_T(\{g,h\}, v)$ is disconnected or has a cut vertex.
 \end{itemize}
\end{theo}

The stable lamination of an automorphism $\phi \in \Out(G)$ is \emph{carried} by a cyclic factor $H$ if and only if there exists $T \in \D$ such that for every leaf $\ell \in (\Lambda^+)_T$ there exists a translate of the minimal subtree $T_H$ which contains $\ell$. 

Minimality of these laminations imply that the stable lamination is carried by $H$ if and only if there exists a leaf $\ell \in (\Lambda^+)_T$ such that $\ell$ is contained in $T_H$. For all $S \in T$, there exists a quasi-isometry $T \to S$ and it implies that the realization of $\ell$ in $S$ is in the subtree $S_H$, so the fact that $\Lambda^+$ is carried by $H$ can be seen in every $S \in \D$. These facts are proved in \cite{PapinDetection}.

Lemma 2.5 of \cite{PapinDetection} implies:
\begin{prop} \label{prop:gphe-lam-connex} 
 Suppose $\phi \in \Out(G)$ is a fully irreducible automorphism with irreducible train track representative $f : T \to T$. Then no leaf of the stable lamination $(\Lambda^+_\phi)_T$ is carried by a cyclic factor.
\end{prop} 

\subsection{Long segments of laminations in axes of elements of $G$} 

In this section, we assume that $\phi \in \Aut(G)$ is a fully irreducible automorphism such that both $\phi$ and $\phi^{-1}$ have train track representatives.

 The following lemma is a transposition of Lemma 2.17 from \cite{PapinWhitehead}. The original lemma gives a link between the Whitehead graph of the axis of a loxodromic element $g \in G$ in a tree $S$ and the existence of a tree $\hat S$ where $\hat S \to S$ is either a fold or a collapse which induces an isometry $\axe_{\hat S}(g) \to \axe_S(g)$. The proof of the lemma does not use the specific fact that $\axe_S(g)$ is an axis of an element and could actually work with any bi-infinite geodesic. In particular, it can be transposed to laminations:
\begin{lem} \label{lem:whitehead-laminations}
 Let $S \in \D$ such that no edge in $S / G$ is a loop. Let $\ell^+$ be a leaf of the stable lamination $\Lambda^+$ and $\ell^-$ be a leaf of the unstable lamination $\Lambda^-$. The following are equivalent:
 \begin{itemize}
  \item There exists a vertex $v \in V(S)$ such that $\Wh_S(\{ \ell^+, \ell^-\}, v)$ is disconnected or has a cut vertex
  \item There exists a tree $\hat S \in \D$ and a non-injective map $\pi:=\hat S \to S$ such that if $\hat \ell^+, \hat \ell^-$ are the leaves in $\hat S$ corresponding to $\ell^+, \ell^-$ then $\pi$ induces isometries $\hat \ell^+ \to \ell^+$ and $\hat \ell^- \to \ell^-$.
 \end{itemize}
\end{lem}
\begin{remas} \label{rem:pliage-ecrase}
 \begin{enumerate}
  \item The assumption that $S$ has no loop is not especially restrictive: in fact, up to subdividing all loops before applying the lemma, we may assume that $S$ has no loop. The tree $\hat S$ produced by the lemma does not have any loop either. 
  \item A non-injective $G$ map $\hat S \to S$ sending vertex to vertex and edge to edge is a composition of collapses and folds (see \cite{BestvinaFeighnBounding}). 
 In particular we can assume that the map $\pi$ given by the lemma is either a collapse or a fold.
 \end{enumerate}
\end{remas}

\begin{prop} \label{prop:existence-arbre-wh-connexe}
 There exists $S \in \D$ such that for every $v \in V(S)$, the Whitehead graph $\Wh_S(\Lambda^+ \cup \Lambda^-, v)$  is connected without cut vertex.
\end{prop}
\begin{proof}
 Let $T$ be the initial train track representative for $\phi$. We will change $T$ gradually using Lemma \ref{lem:whitehead-laminations}.  
 
 By Proposition \ref{prop:gphe-lam-connex}, the Whitehead graph $\Wh_T(\Lambda^+, v)$ is connected for every $v \in S$. The graph $\Wh_T(\Lambda^+ \cup \Lambda^-, v)$ has even more edges so it is also connected.
 
 If an edge of $T$ is a loop then subdivide it so that Lemma \ref{lem:whitehead-laminations} applies. Endow $T$ with the combinatorial metric, i.e. give each edge the length 1. This does not change the Whitehead graphs.
 
 Suppose there exists a Whitehead graph in $T$ which has a cut point. By applying Lemma \ref{lem:whitehead-laminations} and Remark \ref{rem:pliage-ecrase}, we can construct a sequence
 \[
  \dots T_n \to T_{n-1} \to \dots \to T_0 = T
 \]
 where each map $T_i \to T_{i-1}$ is either a collapse or a fold whose restriction to the leaves of both laminations are isometric. The construction of the sequence stops when we find $n$ such that every Whitehead in $T_n$ has no cut vertex.
 
 The maps $T_i \to T_{i-1}$ are in fact not collapses, since a collapsed edge in $T_i$ would not be crossed by any leaf of the lamination, contradicting Proposition \ref{prop:gphe-lam-connex}. They are folds.
 
 \bigskip
 
 We want to prove that the sequence above cannot be infinite. By contradiction, assume it is infinite.
 In \cite[Lemma 1.25]{PapinWhitehead} we proved that the number of orbits of edges of the trees of the sequence built by iterating Lemma \ref{lem:whitehead-laminations} has to go to infinity in that case.
 
 The first Betti number $b_1(T_n)$ is constant. Recall that for a connected graph $\Gamma$ with $V$ vertices and $A$ edges we have $b_1(\Gamma) = A - V + 1$.
 
 Recall that a \emph{big} vertex stabilizer is a vertex stabilizer which does not fix any edge in some (equivalently any) reduced tree.
 By Lemma \ref{lem:big} there is a bound on the number of vertices of valence 1 in $T_n / G$. In fact, the associated vertex groups are big since trees in $\D$ are minimal. Since there exist finitely many conjugacy classes of big vertex stabilizers, this gives a bound on the number of vertices of valence 1 in $T_n/G$. 
 
  Let $A_n, V_n$ be the number of edges and vertices in $T_n / G$. For every $v \in T_n$ denote by $\val(v)$ the valence of $v$. Then we have
 \[
  2V_n + 2 b_1(\Gamma) - 2 = 2A_n = \sum_{v \in T_n / G} \val(v)
 \]
 Therefore
 \[
  2 b_1(\Gamma) - 2 = \sum_{v \in T_n/G} (\val(v) - 2)
 \]
 The only negative terms in the sum correspond to valence 1 vertices so there is a lower bound on their sum. This implies that there is a bound on the number of vertices with valence $\geq 3$.
 
 As a result, since the number of edges in $T_n$ when $n$ goes to infinity is unbounded, the number of vertices of valence 2 in $T_n/G$ is unbounded.
 
 Let $v \in T_n/G$ be a vertex of valence 2. Let $l_1, l_2$ be the labels at $v$. If $|l_1| > 1$ and $|l_2| > 1$ then the stabilizers of vertices in the orbit of $v$ are big. Thus the number of vertices of valence 2 with both labels distinct from $\pm 1$ is bounded by $m(G)$. For every other vertex of valence 2, one of the labels is $1$ or $-1$.
 
 \bigskip
 
 A \emph{topological edge} in $T_n/G$ is a connected component of 
 \[
   \Gamma \setminus \{v \in V(\Gamma) / \val(v) \neq 2 \text{ or no label at } v \text{ is } \pm 1\}
 \] 
 There is a bound $B$, independent of $n$, on the number of topological edges in $T_n/G$. Since $\vol(T_n/G)$ is unbounded there is no bound on the length of topological edges.
 
 A subsegment $\sigma:= e_0, \dots, e_k$ of a topological edge is \emph{increasing} if for any $i \in \{1, \dots, k\}$, $G_{o(e_i)} = G_{e_i}$. It is \emph{decreasing} if for any $i \in \{0, \dots, k-1\}$, $G_{t(e_i)} = G_{e_i}$. A subsegment with a single edge is both increasing and decreasing.
 This can be understood efficiently with labels: $\sigma$ is increasing (resp. decreasing) if $\lambda(e_i) = \pm 1$ (resp. $\lambda(\bar e_i) = \pm 1$) for all $i \in \{1, \dots, k\}$.

 We will now prove that any topological edge of $T_n/G$ can be cut into at most $2m(G)+1$ subsegments which are either increasing or decreasing. See Figure \ref{fig:arete-topo} for an example.
 
 Let $\sigma:= e_0, \dots, e_k$ be a topological edge. Cut it into subsegments by the following process. Let $\sigma_1 := e_0, \dots, e_{i_1}$ be the maximal decreasing prefix of $\sigma$: it has at least one edge. Let $\sigma_2$ be the maximal increasing prefix of $\sigma \setminus \sigma_1$. The label $\lambda(\bar e_{i_1})$ must be nonzero unless $\sigma_1= \sigma$, so since $\sigma$ is a topological edge, $\lambda(e_{i_1 + 1}) = \pm 1$ so $\sigma_2$ also has at least one edge. Continue this procedure to construct an alternating sequence of disjoint decreasing and increasing subsegments.
 
 Write $\sigma$ as the concatenation $\sigma_1, \dots, \sigma_N$ of subsegments. We claim that whenever an increasing subsegment is followed by a decreasing subsegment, the last edge of the former has a big stabilizer, and no edge in the latter does. Thus the number of increasing subsegments in $\sigma$ is bounded by $m(G)+ 1$.
 
 Suppose $\sigma_j$ is an increasing subsegment followed by a decreasing subsegment $\sigma_{j+1}$. Let $e_i$ be the last edge of $\sigma_j$. Then by maximality of $\sigma_j$, $\lambda(e_{i+1}) \neq \pm 1$. Besides, $\sigma_1$ is a decreasing subsegment so there is a decreasing subsegment before $\sigma_j$. Let $e_{i'}$ be its last edge. We have $i' < i$ and $\lambda(\bar e_{i'}) \neq \pm 1$. For all $p \in \{i' + 1, i\}$, $\lambda(e_p) = \pm 1$ so $\sigma_j$ is collapsible and collapses to a vertex $v$ with labels $\lambda(e_{i+1})$ on the right, $\prod_{i' \leq p \leq i} \lambda(e_p)$ on the left. Both labels are not $\pm 1$ so the vertex group associated to $v$ is big. It is also the edge group associated to $e_i$.

 For a topological edge of length $k$, at least one the maximal topological edges is longer than $\frac{k}{2m(G) + 2}$. Thus there is no bound on the maximal length of half topological edges when $n$ increases.
 
  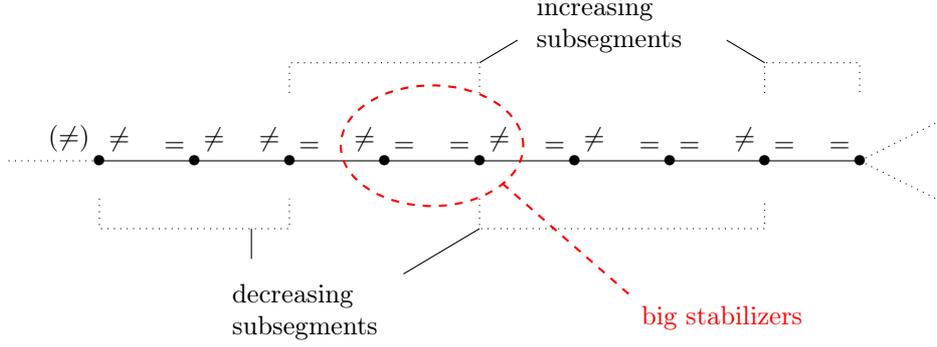
\begin{figure}
  \centering
   \begin{tikzpicture}

 \draw (0,0) node {$\bullet$} node [above right] {$\neq$}  node [above left] {$(\neq)$}-- (10, 0) node {$\bullet$} node [above left] {$=$};
\draw[dotted] (10, 0) -- +(1,0.5);
\draw[dotted] (10, 0) -- +(1,-0.5);
\draw[dotted] (-0, 0) -- +(-1.2,0) ;

\draw (1.25,0) node {$\bullet$} node [above left] {$=$} node [above right] {$\neq$};
\draw (2.5,0) node {$\bullet$} node [above left] {$\neq$} node [above right] {$=$};
\draw (3.75,0) node {$\bullet$} node [above left] {$\neq$} node [above right] {$=$};
\draw (5,0) node {$\bullet$} node [above left] {$=$} node [above right] {$\neq$};
\draw (6.25,0) node {$\bullet$} node [above left] {$=$} node [above right] {$\neq$};
\draw (7.5,0) node {$\bullet$} node [above left] {$=$} node [above right] {$=$};
\draw (8.75,0) node {$\bullet$} node [above left] {$\neq$} node [above right] {$=$};

\draw [dotted] (0,-0.5) |- (2,-0.9) -| (2.5,-0.5);
\draw [dotted] (5,-0.5) |- (8,-0.9) -| (8.75,-0.5);

\draw [dotted] (2.5,0.9) |- (4,1.3) -| (5,0.9);
\draw [dotted] (8.75,0.9) |- (9,1.3) -| (10,0.9);

\draw (3,-2) node [text width=2.5cm] {decreasing \\ subsegments};
\draw (2, -1.3) -- (2, -0.9);
\draw (4, -1.5) -- (5, -0.9);
\draw (7, 1.8) node [text width=2.5cm] {increasing \\ subsegments};
\draw (5.5, 1.6) -- (5, 1.3);
\draw (8.25, 1.6) -- (8.75, 1.3);

\draw [dashed, thick, red] (4.375,0.2) ellipse [x radius=1.2,y radius=0.8];
\draw [dashed, thick, red] (5.3,-0.3) -- (7, -1.8) node [below right]  {big stabilizers};
\end{tikzpicture} 
  \caption{General form of a topological edge and monotonous subsegments} \label{fig:arete-topo}
 \end{figure}

 Suppose $\sigma$ is an increasing subsegment of a topological edge with length $k$ in $T_n/G$. Let us prove that
 there exists a leaf of $\Lambda^+$ and a leaf of $\Lambda^-$ 
 which overlap along a segment with length $k$. Let $\pi_n : T_n \to T_n/G$ be the quotient map.
 
 Write $\sigma$ as the concatenation $e_0, \dots, e_k$ with $\lambda(e_i) = \pm 1$ for all $i \in \{1, \dots, k\}$ (see Figure \ref{fig:pre-image}). Let $w= t(e_k), v=o(e_0)$. 
 
 The subsegment $\sigma$ lifts in $T$ to a subforest $Y_\sigma$. Let $\tilde \sigma$ be  a connected component of $Y_\sigma \setminus \pi_n^{-1}(v)$. It is a finite rooted tree with root $\tilde w$ which is a lift of $w$. For every $i \in \{1, \dots, k\}$, any lift of $\tilde e_i$ is an edge of $T$ which points towards $\tilde w$. The terminal vertices of $\tilde \sigma$, other than possibly $\tilde w$, are the lifts of $v$.
 
 The vertex group $G_{\tilde w}$ acts transitively on the set of lifts of $v$.
 
 Let $\ell$ be a leaf of $\Lambda^+$. There is a translate of $\ell$ which crosses a lift of $e_k$. Thus it contains $[\tilde v, \tilde w]$ where $\tilde v$ is a lift of $v$. By transitivity of the action of $G_{\tilde w}$, for every $\tilde v' \in \pi_n^{-1}(v)$, the segment $[\tilde v', \tilde w]$ is contained in a leaf of $\Lambda^+$.

 Similarly, if $\ell_-$ is a leaf of $\Lambda^-$, there exists a translate of $\ell_-$ crossing $\tilde e_k$ and by translating further by an element of $G_{\tilde w}$ we can make sure that it crosses $[\tilde v', \tilde w]$ for any arbitrary $\tilde v' \in \pi_n^{-1}(v)$. Then the translates of $\ell$ and $\ell_-$ overlap on a length at least $k$.  
 
 \begin{figure}
  \centering
   \begin{tikzpicture}
   \draw (0, 4) node {$\bullet$} node [above right] {$1$}  node [below right = 0.2cm] {$e_0$} -- (1.5, 4) node {$\bullet$} node [above right] {$1$}  node [above left] {$n_0$}   node [below right = 0.2cm] {$e_1$} -- (3, 4) node {$\bullet$} node [above right] {$1$} node [above left] {$n_1$} --  (7.5, 4) node {$\bullet$} node [above right] {$1$} node [above left] {$n_{k-1}$}  -- (9,4)node {$\bullet$} node [above left] {$n_k$}  node [below left = 0.2cm] {$e_k$};

\draw (9,0)  node {$\bullet$} node [above] {$\tilde w$}-- (7.5, 1) -- (6,1.5);
\draw (9,0) -- (7.5, -1) -- (6,-1.5);
\draw (7.5, 1) -- (6, 0.5);
\draw (7.5, -1) -- (6, -0.5);

\draw[dotted] (6, 1.5) -- (5, 1.7);
\draw[dotted] (6, 1.5) -- (5, 1.3);

\draw[dotted] (6, -1.5) -- (5, -1.7);
\draw[dotted] (6, -1.5) -- (5, -1.3);

\draw[dotted] (6, 0.5) -- (5, 0.7);
\draw[dotted] (6, 0.5) -- (5, 0.3);

\draw[dotted] (6, -0.5) -- (5, -0.7);
\draw[dotted] (6, -0.5) -- (5, -0.3);

\draw [dotted] (3, 1) -- (3.5,0.9);
\draw (3, 1) -- (1.5, 1.4);
\draw (3, 1) -- (1.5, 0.6);
\draw (1.5, 1.4)-- (0, 1.6);
\draw (1.5, 1.4)-- (0, 1.2);
\draw (1.5, 0.6)-- (0, 0.8); 
\draw (1.5, 0.6)-- (0, 0.4);

\draw [dotted] (3, -1) -- (3.5,-0.9);
\draw (3, -1) -- (1.5, -1.4);
\draw (3, -1) -- (1.5, -0.6);
\draw (1.5, -1.4)-- (0, -1.6);
\draw (1.5, -1.4)-- (0, -1.2) node[left] {$\tilde v$};
\draw (1.5, -0.6)-- (0, -0.8); 
\draw (1.5, -0.6)-- (0, -0.4); 

\draw (0.5,2.2) node {$\vdots$};
\draw (0.5,-2.2) node {$\vdots$};
\draw (0.5,0) node {$\vdots$};

\draw (4.5,-0.8) node {$\vdots$};
\draw (4.5,0.8) node {$\vdots$};

\draw [thick, green!70!black] (9.2, 0) -- (9,0) -- (7.5, -1)-- (6, -0.5) -- (5, -0.7)  -- (3.5, -0.9) node [midway, above]  {$\ell$} -- (3, -1) -- (1.5, -1.4) -- (0, -1.2)  -- (-0.2, -1.2) ;

\draw [-latex](0,1.8) .. controls (-1, 1) and (-1,-1) .. (0, -1.8) node [midway, right] {$G_{\tilde w}$};

\draw (-0.6, 4.5) node {$T_n/G$};
\draw (-0.6, 2) node {$T_n$};
\end{tikzpicture} 
  \caption{Pre-image of an increasing subsegment} \label{fig:pre-image}
 \end{figure}
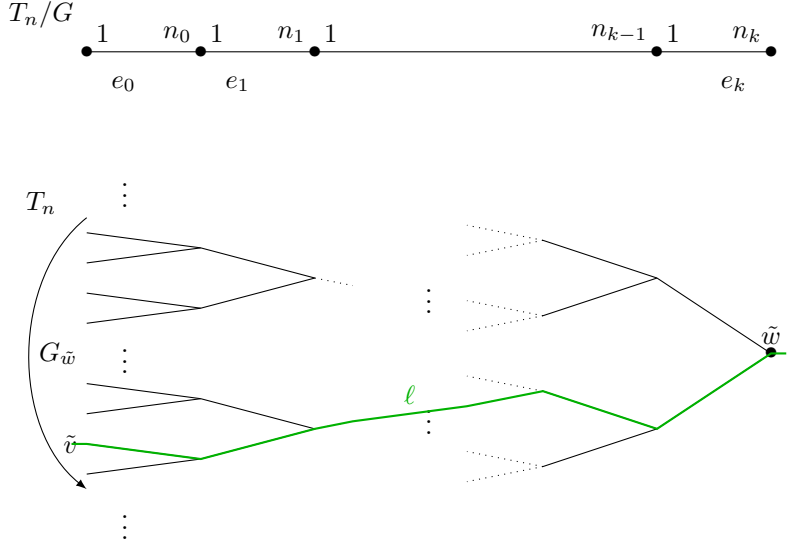

 The same proof can be transposed to the case of a decreasing subsegment.
 
 \bigskip
 
 Since there is no bound on the length of monotonous subsegments, for every $L > 0$ there exists $n \in \N$ such that $T_n$ contains an increasing or decreasing subsegment $\sigma$ longer than $L$. Let $\tilde \sigma$ be a lift for $\sigma$ in $T_n$. There exist leaves $\ell^+, \ell^-$ of the stable and unstable laminations which both cross $\tilde \sigma$. Therefore the leaves $\ell^+$ and $\ell^-$ overlap on a length bigger than $L$.
 
 The maps $T_n \to T$ for $n \in \N$ are isometric in restriction to the leaves of the laminations. Thus there is no bound on the length of common subsegments of both laminations in $T$, so by Lemma \ref{lem:laminations-memes-segments} the laminations are equal. This is a contradiction to Lemma \ref{lem:laminations-distinctes}.
\end{proof}

\begin{prop}  \label{prop:candidats-lamination}
 Let $\phi \in \Out(G)$ be a train track automorphism. Let $f: T \to T$ be a train track representative for $\phi$. Let $\mathcal L _f$ be an axis for $\phi$ in $\D$ passing through $T$. 
 
 There exists $L > 0$ such that for any $T_t \in \mathcal L_f$:
 \begin{enumerate}[label=(\roman*)]
  \item If $g$ is a simple loxodromic element in $G$, then $\axe_{T_t}(g)$ cannot simultaneously contain an $L$-piece of $\Lambda^+$ and an $L$-piece of $\Lambda^-$.
  \item If $g, h$ are simple loxodromic elements such that $\axe_{T_t}(g)$ contains an $L$-piece of $\Lambda^+$ and $\axe_{T_t}(h)$ contains an $L$-piece of $\Lambda^-$, then the pair $\{g,h\}$ is not simple.
 \end{enumerate}
\end{prop}
\begin{proof}
 Let $S$ be a tree obtained with Proposition \ref{prop:existence-arbre-wh-connexe}, i.e. such that for every $v \in V(S)$ the Whitehead graph $\Wh_S(\Lambda^+ \cup \Lambda^-, v)$ is connected without cut vertex. Note that for any $n \in \Z$, $S \cdot \phi^n$ has the same property since $\Lambda^+, \Lambda^-$ are $\phi$-invariant.
 
 By quasi-periodicity of leaves, there exists a constant $L_1 > 0$ such that any leaf segment of $\Lambda^+$ in $S$ (resp. $\Lambda^-$) longer than $L_1$ contains all turns in $\Lambda^+$ (resp. $\Lambda^-$). Suppose $\axe_S(g)$ contains an $L_1$-piece of $\Lambda^+$ and an $L_1$-piece of $\Lambda^-$, then $\Wh_S(g,v)$ is connected without cut vertex for all $v \in V(S)$. By theorem \ref{theo:whitehead-axes} $g$ is not simple. This proves assertion (i) in the specific case where the tree is $S$: now we would like to prove it for $T_t$ in the axis of $\phi$.

 By Remark \ref{rema:borne-kc}
 there exist constants $K>1, C>0$ such that for every $t \in [0, \log(\lambda)]$ there exists a $G$-equivariant $(K,C)$-quasi-isometry $h_t:T_t \to S$. There is an upper bound $B \geq 0$, depending only on $(K,C)$, on $\BBT(h_t)$. 
 
 Let $L_0$ be the constant of Lemma \ref{lem:qi-et-feuille2} for the quasi-isometry constants $K,C$ and $\ell=L_1 + 2B$. If a fundamental domain  $\gamma$ of $g$ in $T_t$ contains a leaf segment $\sigma$ with $\len(\sigma) \geq L_0$ then $[h_t(\gamma)]$ contains a leaf segment with length at least $L_1 +2B$, and a fundamental domain of $g$ in $S$ contains an $L_1$-piece of the lamination. This works for both laminations. 
 
 Let $L > \max\{L_0, 2KL_1\}$. Suppose $\axe_{T_t}(g)$ contains leaf segments $ \sigma_+, \sigma_-$ of the stable and unstable laminations, both longer than $L$. Up to replacing $g$ with $g^l$ for some $l \in \N$ we can suppose a fundamental domain for $g$ contains both $\sigma_+$ and $\sigma_-$. Lemma \ref{lem:qi-et-feuille2} ensures that $\axe_S(g)$ contains $L_1$-pieces of both laminations, therefore implying that $g^l$ hence $g$ is not simple.
 
 Finally suppose $t \notin [0, \log(\lambda)]$. There exists $n \in \Z$ such that $T_{t'} := T_t \cdot \phi^n$ with $t' = t + n \log(\lambda) \in [0, \log(\lambda)]$. Then $h_{t'} : T_{t'} \to S$ induces a $G$-equivariant $(K,C)$-quasi-isometry $T_{t} \to S \cdot \phi^{-n}$. With the same arguments as above we come to the same result, with the same constant $L$. This proves (i).
 
 \bigskip
 
 The proof of (ii) is analogous. We just proved that there exists $L > 0$, such that for $T_t \in \mathcal L_f$ there exists $n$ such that for any $g \in G$, if $\axe_{T_t}(g)$ contains an $L$-piece of any lamination then $\axe_{S \cdot \phi^{n}}(g)$ contains an $L_0$-piece of the same lamination. Applying this to $g$ with the stable lamination and $h$ with the unstable lamination, we get that $\Wh_{S \cdot \phi^n}(\{g, h\}, v)$ is connected without cut vertex for any $v \in V(S)$ and therefore $\{g,h\}$ is not simple.
\end{proof}


\section{Legality} \label{sec:legalite}

Let $G$ be a GBS group with first Betti number $b_1(G) \geq 3$.

In this section we fix a pseudo-atoroidal fully irreducible automorphism $\phi \in \Out(G)$ with a train track representative $f : T \to T$ and a train track representative $f_- : T_- \to T_-$ for $\phi^{-1}$. The goal is to study the evolution of $\|\phi^n(g)\|_T$, $\|\phi^n(g)\|_{T_-}$ when $n \to \pm \infty$ for $g \in G$.

The following three lemmas prove an analogue of Lemmas 2.9 and 2.10 in \cite{BestvinaFeighnHandelLaminations}. The point of view differs a little since we state the results in the trees and not in the quotient graphs.  A notable difference which is caused by non-trivial edge stabilizers is the fact that a concatenation of pINPs is not always a Nielsen path, since it might only be pre-periodic. The statements also differ a little for technical reasons.

\begin{lem} \label{lem:evolutionBoucle}
 Let $\phi \in \Aut(G)$ be an automorphism with a train track representative $f : T \to T$. For every $C > 0$ there exists $M \in \N$ such that for any edge path $\sigma \subset T$, one of the following holds:
 \begin{enumerate}[label=(\alph*)]
  \item $[f^M(\sigma)]$ is legal
  \item $[f^M(\sigma)]$ contains a legal segment of length $\geq C$ between two illegal turns 
  \item $[f^M(\sigma)]$ has fewer illegal turns than $\sigma$
  \item $\sigma$ is a concatenation $\eta_0\cdot \eta_1 \cdot \dots \cdot \eta_{k+1}$ for some $k \geq 1$ where $\eta_0, \eta_{k+1}$ are legal subpaths, and for $1 \leq i \leq k$ the path $f^M(\eta_i)$ is a periodic indivisible Nielsen path, and turns at the concatenation points are legal. 
 \end{enumerate}
\end{lem}

\begin{proof}
 Without loss of generality we may assume $C > C_f$. Then if $\sigma:=[y, y']$ contains a legal path with length $C$ then for any $n \in \N$, $[f^n(\sigma)]$ also does, by definition of the critical constant. 
 
 Let $M \in \N$ be a big enough integer, to be determined later. Suppose there exists a path $\sigma \subset T$ such that both (a), (b) and (c) fail. Since $[f^M(\sigma)]$ cannot have more illegal turns than $\sigma$, it has exactly the same number of illegal turns $k \geq 1$. There exist maximal legal subsegments $\gamma_0, \dots, \gamma_{k}$ such that $\sigma$ is the concatenation $\gamma_0 \cdot \dots \cdot \gamma_k$. Since (a) fails we have $\len(\gamma_i) \leq C$ for every $i \in \{1, \dots, k-1\}$. 
 
 The map $f$ maps legal segments to legal segments 
 and since the number of illegal turns is constant, for every $n \leq M$, there is a unique decomposition $[f^n(\sigma)] = \gamma_0^n \cdot \gamma_1^n \dots \cdot \gamma_k^n$ into maximal legal subsegments. We have $\len(\gamma_i^n) \leq C_f$ for all $i \in \{1, \dots, k-1\}$. 
 
 There are finitely many orbits of edge paths of $T$ with length at most $2C_f$. Let $N$ be the number of orbits of such subpaths.

 Let $i \in \{2, \dots, k-1\}$. There exists $p_i \leq N$ and $g_i \in G$ such that $\gamma_{i-1}^N \cdot \gamma_{i}^N =  g_i \left (\gamma_{i-1}^{N+p_i} \cdot \gamma_{i}^{N+p_i}  \right ) $.
 
 There also exists $p_1 \leq N$ and $g_1 \in G$ such that the restrictions of $\gamma_0^N \cdot \gamma_1^N$ and $g_1 \left (\gamma_0^{N+p_1} \cdot \gamma_1^{N+p_1} \right )$ to a $2C_f$-neighbourhood of the illegal turn are equal. Similarly define $p_k \in \N, g_k \in G$ such that restrictions of $\gamma_{k-1}^N \cdot \gamma_k^N$ and $g_k \left (\gamma_{k-1}^{N+p_k} \cdot \gamma_k^{N+p_k} \right )$ coincide.
 
 By taking the smallest common multiple of all $p_i$ for $i \in \{1, \dots, k\}$ we may replace $p_i$ by some $P \in \N$ which does not depend on $i$.
 
 Thus for every $i \in \{2, \dots, k-1\}$, we have $\gamma_{i-1}^N \cdot \gamma_{i}^N \subset g_i [f^P(\gamma_{i-1}^N \cdot \gamma_{i}^N)]$. The same holds for restrictions to a $2C_f$-neighbourhood of the illegal turn for $i \in \{1, k\}$.
 Note that this implies that there exists no $n \in \N$ such that $\gamma_0^n$ or $\gamma_k^n$ vanish. 
 
 For any $i \in \{1, \dots, k\}$ the path $\gamma_i^N$ contains a unique point $x_i$ such that $x_i =  g_i f^P(x_i)$, and $\gamma_{0}^N$ contains a unique point $x_{0}$ such that $x_{0} = g_{1} f^P(x_{0})$. The point $x_0$ might be equal to $y$, $x_k$ might be equal to $y'$. Note that $g_i \gamma_i^{N+P} = g_{i+1} \gamma_i^{N+P} = \gamma_i^N$ so for $i \in \{0, \dots, k-1\}$ actually $g_{i+1} f^P(x_i) =g_i f^P(x_i) = x_i$. Then for any $i \in \{1, \dots, k\}$ the path $[x_{i-1}, x_{i}]$ is a periodic indivisible Nielsen path. 
 
 For any $i \in \{0, \dots, k\}$ there exists a unique point $y_i \in \sigma$ such that $f^N(y_i)= x_i$. The points $y_0, \dots, y_{k}$ subdivide $\sigma$ into $k+2$ subsegments. Define $\eta_i= [y_i, y_{i+1}]$ for $i\notin \{0,k+1\}$, $\eta_0 = [y, y_0]$ and $\eta_{k+1} = [y_{k}, y']$: we just proved that $\eta_i$ is a pre-Nielsen path for $1 \leq i \leq k$. The other subpaths $\eta_0$ and $\eta_{k+1}$ are legal, and the other subpaths are pre-Nielsen paths, so $\sigma$ satisfies (d). 
 
 \bigskip
 
 The integers $N$ and $P$ only depend on $T$. If $M \geq N + P$ we proved that for any path $\sigma$ such that (a), (b) and (c) fail, (d) holds.
\end{proof}

The following result, which is a key for Lemma \ref{lem:aideLegalite}, implies that when neither Case (a), Case (b) nor Case (d) of Lemma \ref{lem:evolutionBoucle} occur, the decrease of the number of illegal turns is a definite proportion of the length of the segment.
\begin{lem} \label{lem:decroissance-segment}
 Let $\phi \in \Aut(G)$ with a train track representative $f : T \to T$. Let $C > C_f$. Let $M_0 \geq 1$ be the corresponding integer given by Lemma \ref{lem:evolutionBoucle}. There exists $p \in \N$ with the following property. Let $M:=pM_0$. There exists $K < 1$ and $K' \geq 0$ such that for any loxodromic $g \in G$, for any $\sigma \subset \axe_T(g)$, there exists segments $\alpha, \beta, \sigma'$ such that $[f^M(\sigma)] = \alpha \cdot \sigma' \cdot \beta$
 with
 \begin{itemize}
  \item $\len(\alpha) \leq K'/2 , \len(\beta) \leq K' / 2$
  \item $\sigma' \subset \axe_T(\phi^{M}(g))$ or $\len(\sigma') = 0$
  \item either $\sigma'$ contains a legal subsegment with length greater than $C$, or 
    \[
     \len(\sigma') \leq K \len(\sigma) + K'
    \] 
 \end{itemize}
 Moreover, if $\sigma$ contains a legal subsegment longer than $C$, then so does $\sigma'$.
\end{lem}
\begin{proof}
 Let $C > C_f$. Let $g \in G$ and $\sigma \subset \axe_T(g)$. 
 Let $M_0$ be the constant of Lemma \ref{lem:evolutionBoucle} for $C$. 
 
 Let $m \in \N$ be the maximal number of pINPs which may be concatenated in $T$ or $T_-$, from Lemma \ref{lem:courteConcatenation}.
 
 If $\len([f^{M_0}(\sigma)]) \leq K'_0:= 2 \BBT(f^{M_0}) + 2mC$, then define $\alpha:=[f^{M_0}(\sigma)]$, and $\sigma'$ and $\beta$ as single points such that $[f^{M_0}(\sigma)] = \alpha \cdot \sigma' \cdot \beta$. These subsegments satisfy the statement for any choice of $p$, with $\len(\sigma') = 0$.
 
 Suppose $\len([f^{M_0}(\sigma)]) > K'_0$.
 Define $\sigma_0:= [f^{M_0}(\sigma)] \cap \axe_T(g)$, which is not empty, and let $\alpha_0, \beta_0$ be the remaining subsegments.
 Note that if $\sigma$ contains a legal subsegment $\theta$ with length greater than $C > C_f$, then $[f^{iM_0}(\sigma) \cap \axe_T(\phi^{iM_0}(g)) \cap f^{iM_0}(\theta)$  is longer than $C$ for all $i \geq 1$, hence the last statement. 

 Let $\nbl(\theta)$ denote the number of maximal legal subsegments in a segment $\theta$.
 We will now prove that if $\sigma_0$ contains no legal subsegment longer than $C$, then its number of maximal legal subsegments $\nbl(\sigma)$ decreases.
 
 Write $\sigma := \theta_0 \cdot \theta_1 \cdot \dots \cdot \theta_n$ where each subsegment $\theta_i$ except $\theta_0$ has $m+2$ maximal legal subsegments. Let $i \in \{1, \dots, n\}$. Apply Lemma \ref{lem:evolutionBoucle} to $\theta_i$. Case (d) cannot happen since $\theta_i$ has $m+1$ illegal turns.
 If Case (a) happens then Case (c) also happens: $\nbl([f^{M_0} (\theta_i)]) < \nbl(\theta_i) = m+1$. Suppose Case (b) happens: either $[f^{M_0}(\theta_i)] \cap \axe_T(\phi^{M_0}(g))$ contains a legal subsegment with length $C$, or an illegal turn of $\theta_i$ is sent outside $\axe_T(\phi^{M_0}(g))$, in which case $\nbl([f^{M_0}(\theta_i)] \cap \axe_T(\phi^{M_0}(g))) < \nbl(\theta_i)$.
 
 Thus if $\sigma_0$ does not contain any legal subsegment longer than $C$, then 
 \begin{align*}
  \nbl(\sigma_0) &\leq \nbl(\theta_0) + \nbl([f^{M_0}(\theta_1)] \cap \axe_T(\phi^{M_0}(g))) + \dots + 
      \nbl([f^{M_0}(\theta_n)] \cap \axe_T(\phi^{M_0}(g))) \\
   & \leq \nbl(\theta_0) + (\nbl(\theta_1) -1 )+ \dots + (\nbl(\theta_n)-1) \\
   & \leq \nbl(\sigma) - n
 \end{align*}
where $n = \left \lfloor \frac{\nbl(\sigma)}{m+1} \right \rfloor$.

Therefore, with $k := (1 - \frac{1}{m+1})$, we obtain
\[
 \nbl(\sigma_0) \leq k \nbl(\sigma)
\]
This can be iterated as long as $[f^{iM_0}(\sigma)] \cap \axe_T(\phi^{iM_0}(g))$ contains no legal subsegment of length $C$
by applying the same argument to $\sigma_1$ instead of $\sigma$, creating a decreasing subsequence $\sigma_1, \sigma_2, \dots$ of $\sigma$ and increasing sequences $\alpha_1, \alpha_2, \dots$ and $\beta_1, \beta_2, \dots$. Then for $p \geq 1$, $[f^{pM_0}(\sigma)]$ can be cut into subsegments $\alpha_p \cdot \sigma_p \cdot \beta_p$ with 
\begin{itemize}
 \item $\len(\alpha_p), \len(\beta_p) \leq K'_0 \sum_{i=0}^{p-1} \lambda^i =: K'_p$
 \item $\sigma_p \subset \axe_T(\phi^{pM_0}(g))$
 \item $\nbl(\sigma_p) \leq k^p \nbl(\sigma)$
\end{itemize}

Now convert this result into lengths: we obtain
\[
 \len(\sigma_p) \leq \frac{Ck^p}{l_{\min }} \len(\sigma)
\]
where $l_{\min } := \min_{e \in E(t)} \len(e)$.

Choose $p$ such that $\frac{Ck^p}{l_{\min }} <1$. Let $K := \frac{Ck^p}{l_{\min }} $ and $K' := K'_p$.

Finally we get
\[
 \len(\sigma) \leq K \len(\sigma) + K'
\]
thus we obtain the lemma with $\alpha:= \alpha_p$, $\beta:= \beta_p$, $\sigma' := \sigma_p$.
\end{proof}

\begin{lem} \label{lem:aideLegalite} 
 Recall that $\phi$ is pseudo-atoroidal.
 Let $h : T \to T_-$ be a Lipschitz $G$-equivariant map, sending vertex to vertex, and edge to edge path.
 
 For every $C > 0$ there exists $N \in \N$ and $L > 0$ such that for any $g \in G$, 
 for any geodesic $\sigma \subset \axe_T(g)$ with length greater than $L$ and possibly infinite, if $\sigma':= [h(\sigma)]$, then one of the followings holds: 
 \begin{enumerate}[label=(\Alph*)]
  \item $[f^N(\sigma)] \cap \axe_T(\phi^N(g))$ contains a legal segment of length $> C$
  \item $[f_-^N(\sigma')]\cap \axe_{T_-}(\phi^{-N}(g))$ contains a legal segment of length $>C$.
 \end{enumerate}
\end{lem}
\begin{proof}
 Let $M$ be the constant from Lemma \ref{lem:decroissance-segment}. Assume 
 \[
  C \geq \max \left \{C_f, C_{f_-} \right \}
 \]

 We will suppose by contradiction that the lemma fails for $N:=Mi$ with $i$ sufficiently big. We will take a segment $\sigma$ in $T$ for which both (A) and (B) fail, and show that this assumption leads to a contradiction in the following sense. For $j \in \{0, \dots, i\}$ the segment $[f^{Mj}(\sigma)]$ can be cut into three segments: one in the axis of $\phi^{Mj}$ and two ``error'' parts outside of the axis. Using the fact that (A) fails, we see that the part in the axis must not contain any long legal subsegment, thus its length can be estimated by counting the number of maximal legal subsegments in it. Lemma \ref{lem:decroissance-segment} controls the decrease of the number of maximal legal subsegments in $\sigma$ up to error parts. 
 
 However the error parts may grow, for two reasons: they are stretched by $f^M$ and the inner part produces small errors too, which add to the previous error. The aim of the proof is to take $\sigma$ long enough to keep the growth of these error parts small in comparison with the decrease of the inner part, so that the overall effect of $f^{Mi}$ on $\sigma$ is a decrease.
 
 Then we apply the reverse: we look at the evolution of $f_-^{Mj} \circ h \circ f^{Mi}(\sigma)$ for $j \in \{0, \dots, i\}$. Now the argument for the absence of long legal subsegments in the inner part is the fact that (B) fails, and the conclusion is similar: the overall length of the segment decreases. As a result $[f_-^{Mi} \circ h \circ f^{Mi}(\sigma)]$ is a lot shorter than $\sigma$ in proportion.
 
 The contradiction comes from Lemma \ref{lem:applications-proches}: the maps $h$ and $f_-^{Mi} \circ h \circ f^{Mi}$ are equal up to a bounded error, so when $\sigma$ is long enough, it cannot decrease much in proportion when applying $f_-^{Mi} \circ h \circ f^{Mi}$.
 
 We now write a formal argument along these lines.

 \bigskip
 
 Let $i \in \N$; set $N:=Mi$.
 Let $\sigma$ be a segment in $\axe_T(g)$. 
 
 By Lemma \ref{lem:applications-proches}, there exists a constant $B_{i}$ such that for any $x \in T$,  $d_{T_-}(h(x), f_-^{Mi} \circ h \circ f^{Mi}(x)) \leq B_i$. There exists a constant $L_{1,i}$ depending on $B_i$ and $h$ such that if $\len(\sigma) > L_{1,i}$ then $\len([h(\sigma)])> 4B_i$ so
  \begin{align}
  \frac{\len([f_-^{Mj} \circ h \circ f^{Mi}(\sigma)])}{\len(\sigma)} > \frac{1}{2}
  \label{eq:minoration}
 \end{align}
 
 Moreover there exists $D_i \geq 0$ depending on $f, f_-, h, i$ such that for any $\sigma \subset \axe_T(g)$, for any segment $\theta \subset \sigma$ at distance greater than $D_i$ from the endpoints of $\sigma$,  $[f_-^{Mi}\circ h \circ f^{Mi} (\theta)]\cap \axe_T(g) \subset [h(\sigma)] \cap \axe_T(g)$. Let $\tilde \sigma$ be the subsegment of $\sigma$ obtained by cutting out a $D_i$-neighbourhood of the endpoints.
 
 \bigskip
 
 Suppose that both (A) and (B) fail for $\sigma$ and for $N=Mi$.
 
 Since (A) fails for $\sigma$, no segment $[f^{Mj}(\sigma)] \cap \axe_T(\phi^{Mj}(g))$ for $j \in \{0, \dots, i\}$ can contain a legal subsegment longer than $C$. By Lemma \ref{lem:decroissance-segment} there exist constants $K <1, K' \geq 0$ and a decomposition $[f^{M}(\sigma)] = \alpha_1 \cdot \sigma_1 \cdot \beta_1$ such that
  \[
  \len(\sigma_1) \leq K\len(\sigma) + K'
 \]
 Define by induction $\alpha_j \cdot \sigma_{j} \cdot \beta_j = [f^M(\sigma_{j-1})]$ using Lemma \ref{lem:decroissance-segment}. For all $j \in \{1, \dots ,i\}$ we have
 \[
  \len(\sigma_j) \leq K\len(\sigma_{j-1}) + K'
 \]
 and $\len(\alpha_j), \len(\beta_j) \leq K'/2$ so
 \[
  \len([f^M(\sigma_{j-1})]) \leq K\len(\sigma_{j-1}) + K' + 2 \lambda^M \frac{K'}{2}
 \]
 where $\lambda:= \Lip(f)$. Thus we have
 \[
  \len([f^{Mi}(\sigma)]) \leq K^i \len(\sigma) + K'(1+ \lambda^M) \sum_{j = 0}^{i-1} K^j \leq K^i \len(\sigma) + \frac{K'(1 + \lambda^M)}{1-K}
 \]

 Remember that $\tilde \sigma$ be the subpath obtained from $\sigma$ by cutting out the $D_i$-neighbourhood of the endpoints.
 Let $\sigma'_j:=[f_-^{Mj}\circ h \circ f^{Mi}(\tilde \sigma)] \cap \axe_{T_-}(\phi^{M(i-j)}(g))$ for $j \in \{0, \dots, i\}$: it cannot contain legal subsegments longer than $C$. Indeed $[f_-^{Mi} \circ h \circ f^{Mi}(\tilde \sigma)] \cap \axe_T(g) \subset [h(\sigma)] \cap \axe_T(g)$ and by assumption that (B) fails, the latter does not contain any legal subsegments longer than $C$. 
 
 Applying the same argument as above using Lemma \ref{lem:decroissance-segment} with $f_-$, $C$, we obtain again $K_-, K_-'$ such that
 \[
  \len([f_-^{Mi} \circ h \circ f^{Mi}(\tilde \sigma)]) \leq K_-^i \len([h \circ f^{Mi}(\tilde \sigma)]) + \frac{K_-'(1 + \lambda_-^M)}{1-K_-}
 \]
 with $\lambda_-:= \Lip(f_-)$.
 By combining both inequalities, using the fact that $h$ is Lipschitz: 
 \begin{align*}
  \len([f_-^{Mi} \circ h \circ f^{Mi}(\tilde \sigma)]) & \leq \Lip(h) K_-^i    \len([f^{Mi}(\tilde \sigma)])  + \frac{K_-'(1 + \lambda_-^M)}{ 1 - K_-} \\  
   & \leq \Lip(h) K_-^i \len([f^{Mi}(\sigma)])  + 2 \Lip(h) K_-^i D_i \Lip(f^M)^i + \frac{K_-'(1 + \lambda_-^M)}{ 1 - K_-} \\
   &\leq \Lip(h) (K_- K)^i\len(\sigma) + \Lip(h)K_-^i \left ( \frac{K'(1 + \lambda^M)}{1-K} + 2D_i \Lip(f^{M})^i \right ) \\
          & \hspace{2cm} + \frac{K_-'}{1-K_-}
 \end{align*}
 so
 \[
  \len([f_-^{Mi} \circ h \circ f^{Mi}(\sigma)]) \leq \Lip(h) (K_- K)^i\len(\sigma) + S_i
 \]
 where $S_i$ is an additive constant depending on $i$.
 
 Now assume that $i$ is big enough  so that $\Lip(h) (K_- K)^{i} <1/4$. Note that the choice of $i$ does not depend on $\sigma$ but only on the maps $f^M, f_-^M, h$ and on $C$. Then there exists a constant $L_{2,i} \geq 0$ such that if $\sigma$ is longer than $L_{2,i}$ then
 \begin{align}
  \frac{\len([f_-^{Mi} \circ h \circ f^{Mi}(\sigma)])}{\len(\sigma)} \leq \frac{1}{2}
  \label{eq:majoration}
 \end{align}

 \bigskip
 
 If $\sigma$ is longer than $L:=\max \{L_{1,i}, L_{2,i}\}$ then inequations \ref{eq:minoration} and \ref{eq:majoration} contradict each other. This achieves the proof.
\end{proof}

\begin{cor}\label{coro:aideLegalite}
 For every $C > 0$ there exists $N \in \N$ such that  
 for every $g \in G$  one of the followings holds: 
 \begin{enumerate}[label=(\Alph*)]
  \item $\axe_T(\phi^N(g))$ contains a legal segment of length $> C$
  \item $\axe_{T_-}(\phi^{-N}(g))$ contains a legal segment of length $>C$.
 \end{enumerate}
\end{cor}
\begin{proof}
 It suffices to apply Lemma \ref{lem:aideLegalite} to a long enough subsegment $\sigma \subset \axe_T(g)$.
\end{proof}

In the rest of the section, our aim will be to prove that if $g\in G$ then $\axe_{T \cdot \phi^n(g)}(g)$ has increasingly long legal subsegments. This is a key step in the definition of the projection $\D \to \L_f$.

\bigskip

As in \cite{AlgomKfirStrongly} we define the legality threshold and the legality of a path in $T$:
\begin{defi} \label{def:legalite}
 Let $\kappa:= \frac{4\BBT(f)}{\lambda - 1}= 2 C_f$ be the \emph{legality threshold}.
 
 For every finite path $\alpha \subset T$ we define the \emph{legality ratio} of $\alpha$ with respect to the train track structure as follows. Let $\alpha_1, \dots, \alpha_k$ be the maximal legal subsegments of $\alpha$. Then the legality of $\alpha$ is
 \[
  \LEG_f(\alpha,T) := \frac{\displaystyle \sum_{\len(\alpha_i) \geq \kappa} \len(\alpha_i)}{\len(\alpha)}
 \]
 which is the proportion of $\alpha$ which belongs to a legal subpath longer than $\kappa$.
 
 If $g$ is a loxodromic element of $G$, then we distinguish two cases:
 \begin{itemize}
  \item $g$ is legal and we define $\LEG_f(g,T)=1$
  \item there exists a fundamental domain $\alpha$ for $g$ which starts and ends at an illegal turn of the axis. Then $\LEG_f(g,T) = \LEG_f(\alpha,T)$.
 \end{itemize}
\end{defi}
\begin{remas} \label{rem:legaliteDomaineFonda!} 
 \begin{itemize}
  \item If $\alpha$ is a fundamental domain of $g$ which does not start and end at an illegal turn while $\axe_T(g)$ contains one, then $\LEG_f(\alpha) \leq \LEG_f(g)$.
  \item If $(\alpha_n)_{n \in \N}$ is a sequence of nested subsegments of $\axe_T(g)$ whose length goes to infinity then $\LEG_f(g,T) = \lim_{n \to \infty} \LEG_f(\alpha,T)$.
  \item For $l \in \Z \setminus \{0\}$ we have $\LEG_f(g^l, T) = \LEG_f(g,T)$.
  \item To define $\LEG_{f_-}(g, T_-)$ we use the threshold $\kappa_- := 2 C_{f_-}$.
 \end{itemize}
\end{remas}

The following result states that if $\alpha$ contains sufficiently many long legal subsegments, then the length of $f^n(\alpha)$ grows exponentially, as though $\alpha$ were legal.
\begin{lem}\label{lem:legalite-et-f}
 Let $\varepsilon > 0$. There exists a constant $C(\varepsilon)$ such that for every finite path $\alpha$ in $T$ such that $\LEG_f(\alpha, T) \geq \varepsilon$, for every $n \in \N$ we have $\len([f^n(\alpha)]) \geq C(\varepsilon) \lambda^n \len(\alpha)$.
\end{lem}
\begin{proof}
  Since $\LEG_f(\alpha, T) \geq \varepsilon$, $\alpha$ contains at least one legal subsegment of length greater than $\kappa$. Let $\beta_1, \dots, \beta_k$ be the maximal legal subsegments of $\alpha$ longer than $\kappa$. For every $i \in \{1, \dots, k\}$ let $\theta_i$ be the subsegment of $\beta_i$ obtained by cutting out the $\frac{C_f}{2}$-neighbourhood of the endpoints. By Lemma \ref{lem:thetas-disjoints} the images of $\theta_1, \dots, \theta_k$ by $f^n$ are disjoint for any $n \in \N$ and contained in $[f^n(\alpha)]$. Moreover for every $i \in \{1, \dots, k \}$, $\len(\theta_i) \geq \frac{1}{2} \len(\beta_i)$. Since $\kappa = 2 C_f$ we have
 \begin{align*}
  \len([f^n(\alpha)]) &\geq \sum_{i= 1}^k  \len(f^n(\theta_i)) \\
      &\geq  \lambda^n \sum_{i= 1}^k \len(\theta_i)  \\
      &\geq  \frac{1}{2} \lambda^n \sum_{i= 1}^k \len(\beta_i)  \\
      & \geq  \frac{1}{2} \lambda^n \varepsilon \len(\alpha)
 \end{align*}
 Therefore we obtain the desired result, with $C(\varepsilon) =  \varepsilon / 2$.
\end{proof}

\begin{cor} \label{coro:legalite-et-f}
 Let $\varepsilon > 0$. There exists a constant $C(\varepsilon)$ such that for every loxodromic $g \in G$ such that $\LEG_f(g, T) \geq \varepsilon$, for every $n \in \N$ we have $\|\phi^n(g)\|_T \geq C(\varepsilon) \lambda^n \|\phi^n(g)\|$.
\end{cor}
\begin{proof}
 Let $\varepsilon > 0$. Let $g \in G$. There exists $x \in \axe_T(g)$ such that $\LEG_f([x,gx], T)=\LEG_f(g, T)\geq \varepsilon$. By Lemma \ref{lem:legalite-et-f}, for any $k \in \N$ we have $d_T(f^n(x), \phi^n(g^k) f^n(x)) \geq C(\varepsilon) \lambda^n d_T(x, g^kx)$. 
 
 Thus 
 \begin{align*}
   \|\phi^n(g)\|_T &= \inf_{k \in \N} \frac{d_T(f^n(x),  \phi^n(g^k)f^n(x)}{k}\\
    &\geq C(\varepsilon) \lambda^n \frac{d_T(x, g^kx)}{k}\\
    &= C(\varepsilon) \lambda^n \|g\|_T
 \end{align*}
\end{proof}

\begin{figure}
 \centering
   \begin{tikzpicture}[scale=1]

\draw [-latex] (-5, 0) -- (6, 0) node [below] {$n$};

\draw[-latex] (0,0) -- (0,5) node [left] {$\|\phi^n(g)\|_T$};

\draw (0,1) node {$+$} node [above left] {$\|g\|_T$};

\draw [thick, color=red] (0,0.5) .. controls (1,0.6) and (4, 2) .. (5,5); 

\draw [dashed] (0, 1) -| (1.4,0);
\draw [dotted] (1.4, -0.2) |- (6, -0.4);
\draw (3.7, -1) node {$\|\phi^{n}(g) \|_T\geq \|g\|_T$};

\draw (-3,0) node {$+$} node [below] {$m$};
\draw (-3,0.4) node {$+$} node [above left] {$\|\phi^{-m}(g)\|_T$};

\draw [thick, color=green!80!black] (-3,0.2) .. controls (-2,0.24) and (1, 0.8) .. (2,2) .. controls (2.5, 2.6) and (3,3).. (3.5,5); 

\end{tikzpicture}
 \caption{If $\LEG_f(g, T) \geq \varepsilon_0$ then $\|\phi^n(g)\|_T$ is above the red graph. If $\LEG_f(\phi^{-m}(g), T) \geq \varepsilon_0$ then $\|g\|_T$ must be above the green graph, so if $m$ is big then $\|\phi^{-m}(g)\|_T$ is small.} \label{fig:tout-bete}
\end{figure}
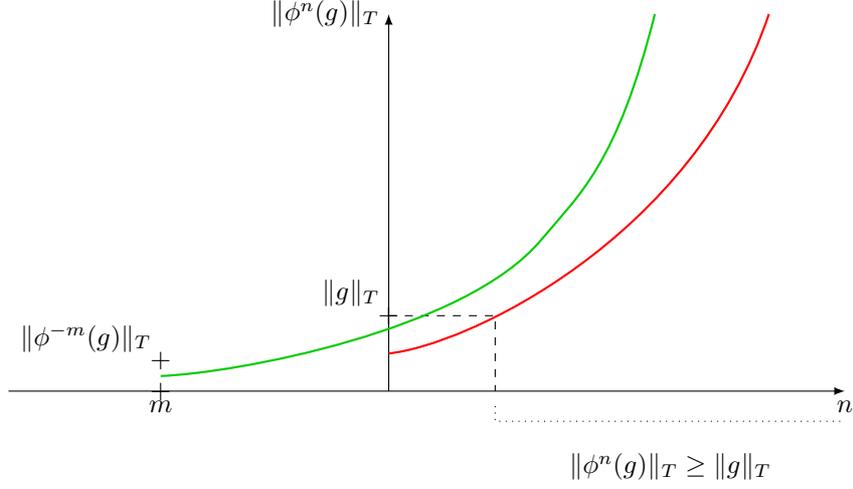

Lemmas \ref{lem:tout-bete} and \ref{lem:legalite-basse-existe} aim to prove basic properties which can be deduced from Lemma \ref{lem:legalite-et-f}. Together they prove that the legality function $n \mapsto \LEG_f(\phi^n(g), T)$ cannot be greater than $\varepsilon$ in a neighbourhood of $-\infty$.
They are illustrated by Figure \ref{fig:tout-bete}.

\begin{lem} \label{lem:tout-bete}
 For any $\varepsilon > 0$ there exists $M \in \N$ such that for any $g \in G$ such that $\LEG_f(g)\geq \varepsilon$, for any $m \geq M$, $\|\phi^m(g)\|_T > \|g\|_T$.
\end{lem}
\begin{proof}
 By applying Lemma \ref{lem:legalite-et-f} to a fundamental domain for $g$ starting at a legal turn, there exists $C(\varepsilon)$ such that for all $m \in \N$
 \[
  \|\phi^m(g)\|_T \geq C(\varepsilon) \lambda^m \|g\|_T
 \]
 so with $M \geq - \frac{\log C}{\log \lambda}$ we get the lemma.
\end{proof}
\begin{rema}
 Similarly there exists $M_-$ such that if $\LEG_{f_-}(g) \geq \varepsilon$ then for any $m \geq M_-$,  $\|\phi^{-m}(g)\|_{T_-} > \|g\|_{T_-}$.
\end{rema}

\begin{lem} \label{lem:legalite-basse-existe}
 For any loxodromic $g \in G$ there exists $m_g \in \N$ such that for any $m \geq m_g$, $\LEG_f(\phi^{-m}(g)) < \varepsilon$. 
\end{lem}
\begin{proof}
 Suppose $\LEG_f(\phi^{-m}(g)) \geq \varepsilon$ for some $m \in \N$ and $\varepsilon > 0$. By applying Corollary \ref{coro:legalite-et-f} to $\phi^{-m}(g)$, 
 there exists $C(\varepsilon)$ such that $\|g\|_T \geq C(\varepsilon) \lambda^m \|\phi^{-m}(g)\|_T$. Let $l_e$ be the length of the shortest edge in $T$, then $\|\phi^{-m}(g)\|_T \geq l_e$ so 
 \[ 
  m \leq \frac{\log(\|g\|_T) - \log C - \log l_e}{\log \lambda}.
 \]
\end{proof}

Corollary \ref{coro:aideLegalite} proves that for any $g$ in $G$, either $\phi^N(g)$ has a $f$-legal segment of length $C$, either $\phi^{-N}(g)$ has an $f_-$-legal segment of length $C$, where the integer $N$ does not depend on $g$ at all. A crucial point is Lemma \ref{lem:legalite-vrai}, i.e. that such a result also works with the legality ratio, i.e. up to choosing a greater $N$, either the legality ratio of $\phi^N(g)$ in $T$ or the legality ratio of $\phi^{-N}(g)$ in $T_-$ is greater than a definite $\varepsilon_0$. 
Combined with Lemma \ref{lem:legalite-et-f} we will then be able to prove that $\len(g)$ grows exponentially when $n \rightarrow \pm \infty$, and has a minimum in a bounded subset of $\L_f$.  

The following lemma needs the fact that $\phi$ is pseudo-atoroidal since it relies on Lemma \ref{lem:courteConcatenation} through Corollary \ref{coro:aideLegalite}. It is proved for the free group case in \cite{BestvinaFeighnHandelLaminations}.

\begin{lem} \label{lem:legalite-vrai}
 There exists $\varepsilon_0 > 0$ and $N\in \N$ such that for every loxodromic element $g \in G$, one of the followings holds:
 \begin{itemize}
  \item $\LEG_f(\phi^N(g), T) > \varepsilon_0$
  \item $\LEG_{f_-}(\phi^{-N}(g), T_-) > \varepsilon_0$
 \end{itemize}
\end{lem}
\begin{proof} 
 Fix a $G$-equivariant quasi-isometry $h: T \to T_-$. 
  Recall that for a $G$-equivariant quasi-isometry between $G$-trees $u : T_1 \to T_2$, such as $h, f, f_-$ and their products, if a segment $\sigma$ is contained in the axis of an element $g$ in $T_1$, then $[u(\sigma)]$ is contained in the axis of $g$ in $T_2$ apart from a $\BBT(u)$-neighbourhood of its endpoints.
  
  \bigskip

 Let $C:= \max \{2C_f, 2C_{f_-}\}$. Let $L, N$ be the constants given by Lemma \ref{lem:aideLegalite}.
 
 There exists a constant $K_N$ depending on $L$ and the quasi-isometry constants for $f, f_-$ such that for any points $x,y \in T$, $f^N(x)=f^N(y) \Rightarrow d_T(x,y) \leq K_N$ and $f_-^N\circ h(x)=f_-^N \circ h(y) \Rightarrow d_{T}(x,y) \leq K_N$).
 
 Observe that for any subsegment $\sigma \subset \axe_T(g)$ of length at least $K_N$, there exists $x \in \sigma$ such that $f^N(x) \in \axe_T(\phi^N(g))$ and $x_- \in \sigma$ such that $f_-^N \circ h (x_-) \in \axe_{T_-}(\phi^{-N}(g))$.

 The axis of $g$ in $T$ can be cut into subsegments $\theta_i, i \in \Z$ of length $L$ separated by other subsegments of length $K_N$ (see Figure \ref{fig:decoupage-axe}).
 
 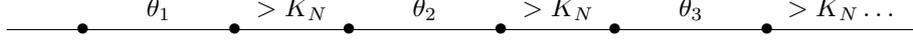
\begin{figure}
  \centering
  \begin{tikzpicture}[]
 \draw (0,0) -- (1,0) node {$\bullet$} -- (3,0) node [midway, above]{$\theta_1$} node {$\bullet$} -- (4.5, 0) node [midway, above] {$> K_N$} node {$\bullet$}
 -- (6.5,0) node [midway, above]{$\theta_2$} node {$\bullet$} -- (8, 0) node [midway, above] {$> K_N$} node {$\bullet$}
 -- (10,0) node [midway, above]{$\theta_3$} node {$\bullet$} -- (12, 0) node [midway, above] {$> K_N \dots$} ;
\end{tikzpicture}
  \caption{The axis of $g$ in $T$ is cut into subsegments $\theta_i$ separated by subsegments longer than $K_N$.} \label{fig:decoupage-axe}
 \end{figure}

 By the choice of $K_N$, for any $i \in \Z$,  $[f^N(\theta_i)] \cap [f^N(\theta_{i+1})] = \varnothing$ and $[f_-^N \circ h(\theta_i)] \cap [f_-^N \circ h (\theta_{i+1})] = \varnothing$.

 There exists a power $g^l$ such that $2L + 3K_N < \|g^l\|_T$. 
 Let $k:= \lfloor \frac{\|g^l\|_T-K_N}{L + K_N} \rfloor$.
 There exists a fundamental domain $\alpha$ for $g^l$ in $\axe_T(g)$ which contains at least $k$ consecutive $\theta_i, \theta_{i+1}, \dots, \theta_{i+k-1}$ of the segments defined above, and at distance at least $K_N$ from its endpoints. 
 
 Let $\alpha_N \subset T$ be a fundamental domain for $\phi^N(g^l)$ contained in $[f^N(\alpha)]$, and let $\alpha_{-N} \subset T_-$ be a fundamental domain for $\phi^{-N}(g^l)$ contained in $[f_-^N \circ h (\alpha)]$.
 Since there is a $K_N$-margin between $\theta_i$, $\theta_{i+k-1}$ and the endpoints of $\alpha$, for every $j \in \{i, \dots, i+k-1\}$ we have $[f^N(\theta_j)] \cap \axe_T(\phi^N(g)) \subset \alpha_N$ and $[f_-^N(\theta_j)] \cap \axe_{T_-} (\phi^{-N}(g)) \subset \alpha_{-N}$.

 By Lemma \ref{lem:aideLegalite}, for each $j \in \Z$, either $[f^N(\theta_j)] \cap \axe_T(\phi^N(g))$ contains an $f$-legal segment with length $C$, or $[f_-^N \circ h (\theta_j)] \cap \axe_{T_-}(\phi^{-N}(g))$ contains an $f_-$-legal segment with length $C$.
 Suppose the first case happens for at least half of the indices in $\{i, \dots, i+k-1\}$. Then since the images of the segments $\theta_j$ do not overlap, there are at least $k/2$ legal segments with length $C$ in the fundamental domain $\alpha_N$. Thus $\LEG_f(\phi^N(g)) \geq \frac{Ck}{2\len(\alpha_N)}$.
 
 Since  $k > \frac{\len(\alpha) - L - 2 K_N}{L + K_N}$ and the fact that $\|g^l\|_T = \len(\alpha) \geq 2L + 3K_N$ we obtain
 \begin{align*}
  \LEG_f(\phi^N(g)) &\geq \frac{C}{2\len(\alpha_N)} \frac{\len(\alpha) - L - 2 K_N}{L + K_N} \\
   &\geq \frac{C \len(\alpha)}{2\len(\alpha_N)} \frac{1 - \frac{L - 2 K_N}{2L+3K_N}}{L + K_N} \\
   & \geq \frac{C \len(\alpha)}{2\len(\alpha_N)} \frac{L + 5K_N}{(L + K_N)(2L + 3K_N)} \\
   & \geq \frac{C}{2 \Lip(f^N)} \frac{L + 5K_N}{(L + K_N)(2L + 3K_N)} > 0
 \end{align*}
The bound does not depend on $g$ nor on $\alpha$.

Similarly, if the second case happens, i.e. if there are more long legal segments in $\alpha_{-N} \subset \axe_{T_-}(\phi^{-N}(g))$, then we get
 \[
  \LEG_{f_-}(\phi^{-N}(g)) \geq  \frac{C}{2 \Lip(f_-^N \circ h)} \frac{1 - \frac{L - 2 K_N}{2L+3K_N}}{L + K_N}
 \]
 Since at least one of these two cases occurs, we can define $\varepsilon_0$ as the smallest of both bounds and we obtain the lemma.
\end{proof}

For a geodesic in $T$ (resp. $T_-$) and a constant $L>0$ we define the \emph{lamination ratio} $\LR(g, T, \Lambda^+, L)$ (resp. $\LR(g, T_-, \Lambda^-, L)$ as the upper bound of the proportion of $\axe_T(g)$ (resp. $\axe_{T_-}(g)$) which can be covered by pairwise disjoint leaf segments of $\Lambda^+$ (resp. $\Lambda^-$) with length at least $L$.

\begin{lem} \label{lem:lamination-ratio}
 Let $\varepsilon_0 > 0$. For any $L > 0$, there exists $N_1 \geq 0$ such that for any $g \in G$, for any $n \geq N_1$, if $\LEG_f(g, T) > \varepsilon_0$ then
 \[
  \LR(\phi^{n}(g), T, \Lambda^+, L) >\varepsilon_0/4.
 \]
\end{lem}
\begin{proof}
 Let $g \in G$ be such that $\LEG_f(g, T) > \varepsilon_0$.
 Let $\beta \subset \axe_T(g)$ be a maximal legal subsegment with length at least $\kappa$. Let $\beta' \subset \beta$ be the subsegment obtained by cutting out the $\frac{C_f}{2}\leq \frac{\kappa}{4}$-neighbourhood of the endpoints. Its length is at least $\len(\beta)-C_f$ and for all $n \in \N$, $f^n(\beta') \subset \axe_T(\phi^{n}(g))$.
 
 There exists $n_1 \in \N$ such that $\lambda^{n_1} C_f \geq 4 l_{\max}$ where $l_{\max} :=\max_{e \in E(T)} \len(e)$. Thus $f^{n_1}(\beta')$ contains at least one edge of $T$. In fact, the number of edges of $T$ contained in $f^{n_1}(\beta')$ is at least $k_{\beta'}:=\left \lfloor \frac{\lambda^{n_1} (\len(\beta) - C_f)}{l_{\max}} \right \rfloor$ and their total length  is at least $c_{\beta'}:=\lambda^{n_1}(\len(\beta) - C_f) - 2l_{\max}$.
 
 There exists $n_2 \in \N$ such that for every edge $e \in T$, $f^{n_2}(e)$ is a leaf segment with length greater than $L$.
 
 Thus $f^{n_1 + n_2} (\beta')$ is contained in $\axe_{T}(\phi^{n_1+n_2}(g))$ and contains at least $k_{\beta'}$ disjoint open leaf segment with length at least $L$ whose total length is at least $\lambda^{n_2}c_{\beta'}$.
 
 \bigskip
 
 Let $\mathcal B$ be the set of $\langle g \rangle$-orbits of maximal legal subsegments of $\axe_T(g)$.
 
 Let $n \in \N$. 
 The proportion of $\axe_T(\phi^{n_1+n_2+ n}(g))$ covered by the leaf segments is at least 
 \[
  \frac{\lambda^{n_2+ n} \displaystyle \sum_{\beta \in \mathcal B}   \lambda^{n_1}(\len(\beta) - C_f) - 2l_{\max}}{\lambda^{n_1+ n_2 + n} \|g\|_T }
 \]
 Thus
 \begin{align*}
  \LR(\phi^{n_1 + n_2 + n}(g), T, \Lambda^+, L) 
     &\geq \frac{\displaystyle \sum_{\beta \in \mathcal B}   \lambda^{n_1}(\len(\beta) - C_f) - 2l_{\max}}{\lambda^{n_1} \|g\|_T } \\
     &\geq \frac{\displaystyle \sum_{\beta \in \mathcal B}   \lambda^{n_1}(\len(\beta) - C_f) - \lambda^{n_1} \frac{C_f}{2}}{\lambda^{n_1} \|g\|_T }  \\
     &\geq \frac{\displaystyle \sum_{\beta \in \mathcal B}   \frac{\len(\beta)}{4}}{\|g\|_T } \\
     &\geq \frac{\varepsilon_0}{4}
 \end{align*}
 This proves the lemma with $N_1= n_1 + n_2$.
\end{proof}

 Lemma \ref{lem:legalite-vrai} yields a constant $\varepsilon_0$. For the rest of the paper we fix such an $\varepsilon_0$. Define $k(g), k_-(g)$ as follows:
 \begin{itemize}
  \item $k(g) = \min\{k \in \Z / \LEG_f(\phi^k(g)) \geq \varepsilon_0 \}$
  \item $k_-(g) = \max\{k \in \Z / \LEG_{f_-}(\phi^k(g)) \geq \varepsilon_0 \}$
 \end{itemize}
 By Lemma \ref{lem:legalite-basse-existe} these integers are well-defined.

\begin{lem} \label{lem:legalitePuissances}
 There exists $N \in \N$ such that for any loxodromic $g \in G$, $| k(g)- k_-(g) | \leq N$.
\end{lem}
\begin{proof}
 Let $g \in G$ be a loxodromic element. Let $N_0$ be the constant given by Lemma \ref{lem:legalite-vrai}. Lemma \ref{lem:legalite-vrai} implies that either $\LEG_f(\phi^{N_0}(\phi^{k(g)-N_0 -1}(g)))$ or $\LEG_{f_-}(\phi^{-N_0} \circ \phi^{k(g)-N_0-1}(g))$ is greater than $\varepsilon_0$. By definition of $k(g)$, the former does not hold so $\LEG_{f_-}(\phi^{-2N_0 + k(g) - 1}(g))>\varepsilon_0$. Therefore we have 
 \[
  k(g) - k_-(g) \leq 2N_0 + 1.
 \]
 
 Lemma \ref{lem:tout-bete}
 gives $M$ such that for all $m \geq M$ we have 
 \[ 
  \|\phi^m(\phi^{k(g)}(g))\|_T > \|\phi^{k(g)}(g) \|_T.
 \]
 It also gives a similar constant $M_-$ for $f_-$. If we had $k_-(g) - k(g)\geq \max\{M, M_-\}$ then by applying Lemma \ref{lem:tout-bete} twice with $m = k_-(g) - k(g)$ we would get a contradiction:
 \[
  \| \phi^k(g) \|_T = \|\phi^{-m} \circ \phi^{m} \circ \phi^k(g) \|_T > \|\phi^{m} \circ \phi^k(g) \|_T  > \| \phi^k(g) \|_T 
 \]
 This gives an upper bound for $k_-(g) - k(g)$. 
\end{proof}


\section{Defining the projection} \label{sec:projection}

Let $g \in G$ be a loxodromic element. Like in \cite{AlgomKfirStrongly} we define $t_0(g):= k(g) \log (\lambda_+)$. The following lemma is the same as \cite[Lemma XX]{AlgomKfirStrongly}. The fact that $G$ is a GBS group instead of $F_N$ has no influence.

\begin{lem} \label{lem:estimation-longueur-g}
 There exists a constant $C>0$ such that
 for every loxodromic element $g \in G$ we have for $t \geq t_0$:
 \[
  C^{-1} \lambda^{\left \lfloor \frac{t - t_0}{\log(\lambda)} \right \rfloor}\|g\|_{T_0} \leq  \|g\|_{T_t} \leq C \lambda^{\left \lfloor \frac{t - t_0}{\log(\lambda)} \right \rfloor}\|g\|_{T_0}
 \]
 and for $t \leq t_0$:
 \[
  C^{-1} \lambda_-^{\left \lfloor \frac{t_0- t}{\log(\lambda)} \right \rfloor} \|g\|_{T_0} \leq  \|g\|_{T_t} \leq C \lambda_-^{\left \lfloor \frac{t_0 - t}{\log(\lambda)} \right \rfloor} \|g\|_{T_0}
 \]
\end{lem}
\begin{proof}
 We will prove the inequalities in the case where $t$ is a multiple of $\log(\lambda)$. The result for other values of $t$ can be obtained by applying Lemma \ref{lem:distance-entre-axes} to a translate of the subsegment $\{T_t/ 0 \leq t \leq \log(\lambda)\}$, and it will only result in increasing the multiplicative constants by a controlled amount.

 Write ${t_0} = t_0(g)$.
 First let us deal with the case $t \geq t_0$. Let $n \in \N$ and let $t_n = t_0 + n \log \lambda$.
 
 We have $T_{t_n} = T_{t_0} \cdot \phi^n$.
 
 Since $f$ is $\lambda$-Lipschitz we have
 \[
  \|g\|_{T_{t_n}} \leq \lambda^{n} \|g\|_{T_{t_0}}
 \]

 Let us prove the other side of the inequality. Lemma \ref{lem:legalite-et-f} can be  applied to a well-chosen fundamental domain for $g$ and
 gives a constant $C(\varepsilon_0)$ independant of $n$ and $g$ such that
 \[
  \|g\|_{T_{t_n}} \geq C(\varepsilon_0) \lambda^{n} \|g\|_{T_{t_0}}
 \]
 which gives the first part of the Lemma.
 
 \bigskip
 
Now let us deal with the case $t \leq t_0$. Let $n \in \N$.
Let $t_n = t_0 - n \log(\lambda)$. In that case we have $n= \frac{t_0 - t}{\log(\lambda)}$. 

We have
\[
 \Lip(T, T_-)^{-1} \leq \frac{\|g\|_{T}}{\|g\|_{T_-}} \leq \Lip(T_-, T)
\]
and by applying this to $\phi^n(g)$ instead of $g$ for any $n \in \Z$, we have 

\[
 \Lip(T, T_-)^{-1} \leq \frac{\|g\|_{T\cdot \phi^n}}{\|g\|_{T_- \cdot \phi^n}} \leq \Lip(T_-, T)
\]
In particular, since $T_{t_0}= T \cdot \phi^{k(g)}$ this also works when replacing $T, T_-$ with $T_{t_0}, {T_-}_{t_0}$.

We now have ${T_-}_{t_n} = {T_-}_{t_0} \cdot \phi^{-n}$ so
 $\|g\|_{{T_-}_{t_n}} \leq \lambda_- ^{n} \|g\|_{{T_-}_{t_0}}$. 
We deduce the right inequality:
\begin{align*}
 \|g\|_{T_{t_n}} &\leq  \Lip(T_-, T) \|g\|_{{T_-}_{t_n}} \\
        & \leq  \Lip(T_-, T) \lambda_- ^{n} \|g\|_{{T_-}_{t_0}} \\
        & \leq  \Lip(T_-, T)\Lip(T, T_-) \lambda_- ^{n} \|g\|_{{T}_{t_0}}
\end{align*}

Now let us prove the left inequality.
Lemma \ref{lem:legalite-vrai} gives an integer constant $N$ such that $k_-(g) \geq k(g)- N$. Thus we have by the same arguments as above, and for $n \geq N$ and $T_{t_N} = T_{t_0} \cdot \phi^{-N}$ we obtain
\[
 \|g\|_{{T_-}_{t_n}}  \geq C(\varepsilon_0) \lambda_-^{n - N} \|g\|_{{T_-}_{t_N}}
\]
Since $\|g\|_{{T_-}_{t_N}} \geq \Lip(T_-\cdot \phi^{-N}, T_- )^{-1}\|g\|_{{T_-}_{t_0} }$ we have
\[
 \|g\|_{{T_-}_{t_n}}  \geq C(\varepsilon_0) \lambda_-^{- N} \Lip(T_-\cdot \phi^{-N}, T_- )^{-1} \lambda_-^{n}\|g\|_{{T_-}_{t_0}}
\]
When $0 \leq n < N$ then we have $\frac{\|g\|_{{T_-}_{t_0}}} {\|g\|_{{T_-}_{t_n}}} \leq \Lip(T_-, T_- \cdot \phi^{N})$ so in any case there is a constant $C>1$ depending only on $T, T_-, N$ such that
\[
 \|g\|_{T_{t_n}} \geq C^{-1} \lambda_-{n(t)}\|g\|_{T_{t_0}}
\]
This proves the lemma.
\end{proof}

Define $\Theta(g):=\{ t \in \R/ \|g\|_{T_t} \text{ minimal } \}$. 

\begin{lem}\label{lem:projection-g}
 There exists $s >0$ such that for every loxodromic $g \in G$ and $t \in \Theta(g)$, then $|t - t_0(g)| <s$.
\end{lem}
\begin{proof}
 Let $C$ be the constant from Lemma \ref{lem:estimation-longueur-g}. 
 Suppose $t > t_0(g)$. Then we have $\|g\|_{T_t}\geq C^{-1} \lambda^{\left \lfloor \frac{t - t_0(g)}{\log(\lambda)} \right \rfloor} \|g\|_{T_{t_0}}$ so since $\|g\|_{T_t} \leq \|g\|_{T_{t_0}}$ this implies $t - t_0 \leq \log( C \lambda)$.
 
 We get a similar inequality for $t < t_0$, hence the result.
\end{proof}
\begin{rema}
 The diameter of $\Theta(g)$ is bounded by $2s$. 
\end{rema}

An important property of the projection is that projections of simple pairs are close:
\begin{lem} \label{lem:distance-gh}
 There exists $s' > 0$ with the following property.
 Let $\{g, h\}$ be a simple pair of loxodromic elements of $G$. Then $|t_0(g) - t_0(h) | < s'$.
\end{lem}
\begin{figure}
 \centering
 \begin{tikzpicture}

  \draw[-latex] (0,0.5) .. controls (2, 0.3) and (6, 0) .. (10,0.5) node [right] {$\phi$}; 
  \draw[latex-] (0,-2) node [left] {$\phi^{-1}$} .. controls (1.8, -1) and (7, -1) .. (10,-2); 

 \draw  (0.9, 0.4) node {$\bullet$} node [below] {$t_0(g)$};
\draw  (9.6, 0.45) node {$\bullet$} node [above] {$t_0(h)$};
\draw  (9.5, -1.85) node {$\bullet$} node [above] {$t_0(h)$};

 \draw[-latex, blue] (1, 0.5) .. controls (1.5, 0.6) and (2, 0.6) .. (2.7,0.4) node [midway, above] {$n \log \lambda$}; 
 \draw  (2.8, 0.25) node {$\bullet$} node [below] {$t_1$};

\draw  (7.6, -1.5) node {$\bullet$} ;
 \draw[-latex, blue] (9.4,-2) .. controls (9, -2) and (8.5, -1.9) .. (7.8,-1.7) node [midway, below] {$N \log \lambda_-$}; 

\draw (6,-1.3) node {$\bullet$}; 
\draw[latex-, blue] (6.2,-1.5) .. controls (6.4, -1.6) and (7, -1.7) .. (7.4,-1.7) node [midway, below] {$n \log \lambda_-$}; 

\draw (4.3,-1.2) node {$\bullet$}; 
\draw[-latex, blue] (5.8,-1.45) .. controls (5.2, -1.5) and (5, -1.5) .. (4.5,-1.4) node [midway, below] {$m \log \lambda_-$}; 

\draw [dotted] (4.3, -1.2) -- (4.3, 0.22) node {$\bullet$};

\draw [latex-, very thin] (7.6, -2.5) |- (8, -3.5) node [right, text width=3cm] {$LEG_{f_-}(h)$ becomes $\geq \varepsilon_0$};
\draw [latex-, very thin] (6.0, -2.5) |- (6, -3.5) node [below, text width=3.3cm] {$\axe(h)$ contains long leaf segments of $\Lambda^-$};
\draw[latex-, very thin] (4.3, 0.5) |- (6, 1.8) node [right, text width=3.9cm] {$\axe(h)$ contains long leaf segments of $\Lambda^- $};
\draw [latex-, very thin] (2.8, 0.5) |- (2.8, 1.8) node [above, text width=3.9cm] {$\axe(g)$ contains long leaf segments of $\Lambda^- $};

\end{tikzpicture}
 \caption{If $t_0(g), t_0(h)$ are sufficiently far apart, then there exists $t$ such that the axes of $g$ and $h$ in $T_t$ contain long leaf segments of the opposite laminations} \label{fig:distance-gh}
\end{figure}
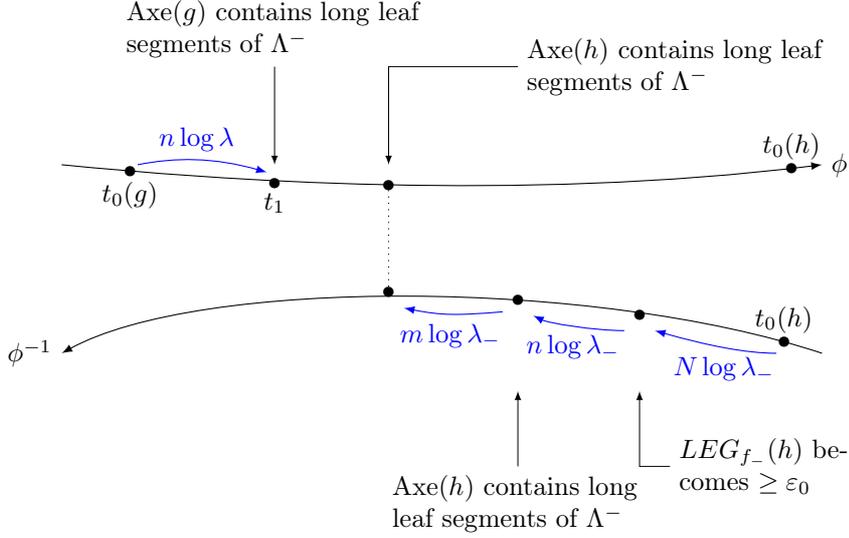

\begin{proof}
  We prove this by contraposition: we will show that if $t_0(g)$ and $t_0(h)$ are too far apart then we can find $t$ in between such that the axes of $g$ and $h$ in $T_t$ contain long segments of the stable and unstable lamination (see Figure \ref{fig:distance-gh}).
  
  Let $L_0$ be the constant from Proposition \ref{prop:candidats-lamination} such that elements in a simple pair cannot contain leaf segments longer than $L_0$ of both lamination in their axes. Without loss of generality we may assume $L_0 > \kappa$.
  
  By Lemma \ref{lem:qi-et-feuille2} there exists a constant $L_1$  such that if a path $\beta \subset T_-$ contains a leaf segment of $\Lambda^-$ longer than $L_1$, then $[h_- (\beta)]$ contains a leaf segment of $\Lambda^-$ of length greater than $L_0$. The choice of $L_1$ depends only on $h_-, f_-$, and $L$.
 
  \bigskip 
  
  Without loss of generality suppose $t_0(g) < t_0(h)$.
  
  Let $N_1 > 0$ be the integer given by Lemma \ref{lem:lamination-ratio} for $f$, $\varepsilon_0$ and $L_0$. Similarly define $N_{1,-} >0$ as the integer given for $f_-$, $\varepsilon_0$, $L_1$.
  
  By definition of $t_0$ we have $\LEG_f(g, T_{t_0(g)}) \geq \varepsilon_0$. Thus for all $n \geq N_1$ we have 
  \[
   \LR(g, T_{t_0(g) + n \log(\lambda)}, \Lambda^+, L_0) \geq \varepsilon_0 / 4.
  \]

  Similarly for all $n \geq N_{1,-}$ we have $\LEG_{f_-}(h, T_{t_0(h) - n \log(\lambda_-)}, \Lambda^-, L_1) \geq \varepsilon_0/4$.
  
  \bigskip
  
  Suppose $t_0(h) - t_0(g) > N_1 \log(\lambda) + N_{1,-} \log(\lambda_-)$. Then there exists $ t_0(g)+ N_1 \log(\lambda) < t < t_0(h)- N_{1,-} \log(\lambda_-)$. Consequently $\axe_{T_t}(g)$ contains an $L_0$-piece of $\Lambda^+$ and $[h (\axe_{{T_-}_t}(h))]$ contains an $L_0$-piece of $\Lambda^-$. This contradicts the fact that the pair $\{g,h\}$ is simple.
\end{proof}
 Here is a direct corollary:
\begin{cor}
 Let $s, s'$ be the constants from Lemmas \ref{lem:projection-g} and \ref{lem:distance-gh}. For a simple pair $\{g, h\}$, $\diam(\Theta(g) \cup \Theta(h)) < s + s'$.
\end{cor}

 In order to evaluate the distance $d_{\Lip}(X, \L_f)$ for some arbitrary $X \in \D$, we will use candidates of $X$. Lemma \ref{lem:distance-gh} applies in particular to candidates:
\begin{cor}\label{coro:candidats-proches}
 Suppose $b_1(G) \geq 3$.
 There exists $s'' > 0$ such that for every $X \in \D$, if $g, h$ are candidates in $X$, then for any $t_g \in \Theta(g)$ and $t_h \in \pi_f(h)$ we have $|t_g - t_h| < s''$.
\end{cor}
\begin{proof}
 By Lemma \ref{lem:troisieme-candidat-compatible} there exists $k \in G$ such that the pairs $\{g, k\}$ and $\{h, k\}$ are simple. Applying the previous corollary gives $|t_g- t_h| < |t_g - t_k| + |t_k - t_h| < 2s' + 2s$.
\end{proof}

For $X \in \D$, define $\Theta_X:= \{t \in \R / d_{\Lip}(X, T_t) \text{ minimal }\}$.

\begin{lem}
 For every $X \in \D$, the set $\Theta_X$ is non-empty.
 
 Moreover there exists $s >0$ such that for every $X \in \D$, $\diam( \Theta_X) < s$.
\end{lem}
\begin{proof}
 Let $X \in \D$. By Theorem \ref{theo:distance-lipschitz-candidats}, for all $t \in \R$, there exists a candidate such that $d_{\Lip}(X, T_t)= \log \frac{\|g\|_{T_t}}{\|g\|_X}$.
 
 Therefore we have 
 \[
  d_{\Lip}(X, T_t) = \max_{g \text{ candidate } } \log \frac{\|g\|_{T_t}}{\|g\|_X}
 \]

 Fix a candidate $g$. The function $t \mapsto \frac{\|g\|_{T_t}}{\|g\|_X}$  is minimal for $t \in \Theta_g$. We will prove that $d_{\Lip}(X, T_t)$ reaches its minimum in a $D$-neighbourhood of $\Theta_g$, where $D$ does not depend on $X$ nor on the number of candidates in $X$.

 Let $t_0 := t_0(g)$. 
 If $h$ is another candidate we have $|t_0(h) - t_0|< s''$ where $s''$ is the constant from Corollary \ref{coro:candidats-proches}. By Lemma \ref{lem:estimation-longueur-g} we have
 \[
  \|h\|_{T_{t_0}} \leq \|h\|_{T_{t_0(h)}} C \lambda^{\frac{s''}{\log(\lambda)}}
 \]
 Write $K = C \lambda^{\frac{s''}{\log(\lambda)}}$. For $t_* \geq \log(2CK)  + \log (\lambda)$ and $t > t_0 + s'' + t_*$ we have $t > t_0(h) + t_*$ and still by Lemma \ref{lem:estimation-longueur-g} we get
 \[
  \|h\|_{T_t} \geq 2K \|h\|_{T_{t_0(h)}} \geq 2 \|h\|_{T_{t_0}}
 \]
 
 Dividing both sides by $\|h\|_X$ does not change the inequality.
 
 \medskip
 
 Therefore we have for $t > t_0+ s'' + t_*$
 \[
  \max_{h \text{ candidate}} \frac{\|h\|_{T_t}}{\|h\|_X} \geq 2 \frac{\|h\|_{T_{t_0}}}{\|h\|_X}
 \]
 
 \medskip
 
 For $t < t_0$ we get a similar result. We deduce a constant $\Delta t$ such that for $t \in \R \setminus [t_0 - \Delta t, t_0 + \Delta t]$ we have $d_{\Lip}(X, T_t) > \log 2 d_{\Lip}(X, T_{t_0})$. Since $t \to T_t$ is continuous for the axes topology, $t \mapsto d_{\Lip}(X, T_t)$ reaches its minimum in a $\Delta t $-neighbourhood of $t_0$. 
\end{proof}

\begin{rema}
 The previous proof differs slightly from the proof of \cite[XX]{AlgomKfirStrongly} since there is no bound on the number of candidates in elements of $\D$, unlike in $\CV_N$. This comes from the fact that $\D$ is not finite dimensional so there is no bound on the number of orbits of edges in elements of $\D$.
\end{rema}

For $X \in \D$ we choose $t_X$ in  $\Theta_X$. Since $\Theta_X$ has bounded diameter and the bound does not depend on $X$ this will be well enough defined.

\section{Negative curvature properties of the projection} \label{sec:negative}

In this section we prove the analogues of Lemmas 5.7 and 5.8 of \cite{AlgomKfirStrongly} and deduce the strong contraction property.

The difference in the proof with \cite{AlgomKfirStrongly} is the proof of Lemma \ref{lem:passages-aretes}. The initial proof relies on special shapes of graphs such as roses (see \cite[Proposition 5.10]{AlgomKfirStrongly}). Here reduced graphs will take the role of roses. The other proofs are actually quite similar to the free group case.

\begin{lem}\label{lem:projections1} 
 There exist $s, c>0$ such that for any $X \in \D$, if $|t - t_X| > s$  then $d_{\Lip}(X, T_t) \geq d_{\Lip}(X, \pi(X)) + d_{\Lip}(\pi(X), T_t) - c$. 
\end{lem}
\begin{proof}
 Suppose $t \leq t_X$. Let $g$ be a candidate in $X$. The idea of the proof is that if $s$ is big enough, then $\LEG_f(g, T_{t_X + s})$ is also big and $g$ almost realizes $d_{\Lip}(T_{t_X+s}, T_t)$.
 
 There is a candidate $h$ of $X$ such that $\Lip(X, \pi(X)) = \frac{ \|h\|_{\pi(X)}}{\|h\|_X}$. Since $t_X \in \Theta(X)$ we have $t_X \in \Theta(h)$. By Lemmas \ref{lem:projection-g} and \ref{lem:distance-gh} there exists a constant $s$ such that for every candidate $g$ of $X$ we have $|t_X - t_0(g)| < s$. Thus for any candidate $g$ of $X$, for any $t_1 > t_X + s$, we have $\LEG_f(g, T_1) > \varepsilon_0$. 
 
 Let $Z := T_{t_1}$.
 Let $g$ be a candidate of $X$ such that $d_{\Lip}(X,Z) = \frac{\|g\|_{Z}}{\|g\|_X}  $.  
 
 Applying twice Lemma \ref{lem:estimation-longueur-g} to $g$ for $t$ and $t_1$ we obtain a constant $C$ such that 
 \[
 \frac{\|g\|_{T_t}}{\|g\|_Z} \geq C^{-2} \lambda^{\frac{t-t_0(g)}{\log(\lambda)}-  \frac{t_1 - t_0(g)}{\log(\lambda)}  - 1} = C^{-2}\lambda^{-1} e^{t-t_1}
 \]
 Remarking that $e^{t-t_1} = \Lip(Z, T_t)$ we have
 \[
  \Lip(Z, T_t) \leq \frac{\|g\|_{T_t}}{\|g\|_Z} \frac{1}{C^2 \lambda}
 \]
 with $C^2 \lambda > 1$.
 
 Since $\Lip(X, Z) =\frac{ \|g\|_{Z}} {\|g\|_X}$ we have
 \begin{align*}
  \Lip(X, T_t) & \geq \frac{ \|g\|_{T_t}} {\|g\|_X} \\
        &= \frac{ \|g\|_{T_t}} {\|g\|_Z} \frac{ \|g\|_{Z}} {\|g\|_X}\\
        & \geq \frac{1}{C^2 \lambda} \Lip(Z, T_t) \Lip(X,Z)
 \end{align*}
 Applying the logarithm we get a constant $K > 0$ such that $d_{\Lip}(X, T_t) \geq d_{\Lip}(X, Z) + d_{\Lip}(Z, T_t) - K$.

 \bigskip
 
 Finally by definition of the projection we have $d_{\Lip}(X,Z) \geq d_{\Lip}(X, \pi(X)$. If $t - t_X >s$ we have
 \begin{align*}
  d_{\Lip}(X, T_t) &\geq d_{\Lip}(X, \pi(X) + d_{\Lip}(Z, T_t)-K\\
                & \geq d_{\Lip}(X, \pi(X)) + d_{\Lip}(\pi(X), T_t) - s -K
 \end{align*}
\end{proof}

\begin{lem}\label{lem:projections2} 
 There exist $s, c>0$ such that for any $X, Y \in \D$, if $|t_X - t_Y|> s$ then $d_{\Lip}(Y, X) \geq d_{\Lip}(Y, \pi(X)) -c$.
\end{lem}
Before proving Lemma \ref{lem:projections2} we need some preliminary results.

\begin{lem} \label{lem:application-sympa}
 Let $X, T \in \D$ and $e_0 \in E(X)$.
 Suppose every edge orbit in $X\setminus G \cdot e_0$ is non-collapsible. There exists a $G$-equivariant map $\tau : X \to T$ such that every edge in $X \setminus G \cdot e_0$ is contained in a $\tau$-legal bi-infinite geodesic in $X \setminus G \cdot e_0$.
\end{lem}
\begin{proof}
 We will prove this by constructing the map $\tau : X \to T$ such that at every vertex $v \in V(X)$, at least two gates at $v$ for the gate structure induced by $\tau$ contain edges in $E(X) \setminus G \cdot e_0$. Then there exist bi-infinite $\tau$-legal geodesics with the desired property. 
 
 \bigskip

 There exists a $G$-equivariant map $\tau_0 : X \to T$. We may suppose that $\tau_0$ sends vertex to vertex and is linear on edges.
 
 Let $v_1, \dots, v_n$ be representatives of every vertex orbit of $V(X)$. In order to define a new map $\tau$, it suffices to choose the image of $v_i$ for every $i \in \{1, \dots, n\}$. The image of $v_i$ can be any vertex $w_i \in V(T)$ such that $G_{v_i} \subset G_{w_i}$.
 
 Suppose there exists $v \in X$ such that there is only one gate at $v$ which contains edges in $E(X) \setminus G \cdot e_0$. We include the case where there is only one gate at $v$.
 
 Let $w = \tau_0(v)$. There exists a vertex $w' \in T$ such that $\bigcap _{e \in E_v  \setminus G \cdot e_0} \tau_0(e) = [w, w']$. Since no edges in $\tau_0$ are collapsible except translates of $e_0$, the images $\tau_0(e)$ have non-zero length. Since  these edges are contained in a single gate, the intersection has non-zero length and we have $w \neq w'$. 
 
 Let us prove $G_v \subset G_{w'}$. Let $a \in G_v$. By contradiction suppose $aw' \neq w'$, then for any edge in $E_v \setminus G \cdot e_0$ we have $w' \notin \tau(ae)$ which contradicts the definition of $w'$, thus $G_v \subset G_{w'}$. Define $\tau_1$ by
 \begin{alignat*}{4}
	 \tau_1: ~ && x \in V(X) && ~\longmapsto &~ \left\{\begin{array}{ll}	\tau_0(x) & \mbox{if } x \notin G \cdot v \\ g w' & \mbox{if } x = g v	\end{array}\right.
 \end{alignat*}
 Note that if there exists $e \in E_v$ such that $\tau_0(e) = [w, w']$, then $e$ cannot be a loop in the quotient. If $e$ were a loop $[v, gv] \subset T$ then $g^{-1} \bar e = [v, g^{-1}v]$, which is also in $E_v$ and has same length, would also be sent to $[w,w']=[w, gw] = [w, g^{-1}w]$. This would imply that $g^2$ is elliptic, which is a contradiction. Thus the vertices of $e$ are in distinct orbits, and since the image of $e$ is a single point, $\tau_1$ factors through the collapse of $e$.  
 
 Since no edge orbit is collapsible except $e_0$, $e$ can be collapsed only if $e$ is a translate of $e_0$. We do not care about the image of $e_0$.

 By construction of $\tau_1$, there are at least two gates for $\tau_1$ at vertex $v$ which contain edges in $E_v \setminus G \cdot e_0$. Let us prove that gates at other vertices have not changed.
 
 Let $e \in E(T) \setminus G \cdot e_0$. Neither $\tau_0$ nor $\tau_1$ collapse $e$. If an endpoint $x$ of $e$ is in $G \cdot v$, then $\tau_1(x)$ is in the interior of $\tau_0(e)$. Otherwise $\tau_1(x) = \tau_0(x)$. Thus if $o(e) \notin G \cdot v$, the first edge of $\tau_1(e)$ and $\tau_0(e)$ are equal. Thus the gate structure at $o(e)$ is unchanged.

 Therefore the number of vertices with at most one gate containing edges not in $G \cdot e_0$ is smaller for $\tau_1$ than for $\tau_0$. 
 
 \bigskip
 
 We can iterate this procedure with $\tau_1$ instead of $\tau_0$ until we find a map $\tau$ such that all vertices in $X$ have at least two gates containing edges not in $G \cdot e_0$.
\end{proof}

\begin{lem} \label{lem:passages-aretes}
 Let $L$ be the constant from Proposition \ref{prop:candidats-lamination}.
 Suppose $r \in \R$ is such that for every candidate $u \in G$ for $X$, the axis $\axe_{T_r}(u)$ contains an $L$-piece of $\Lambda^-$. Let $g \in G$. Suppose that $\axe_{T_r}(g) / \langle g \rangle$ contains $k$ disjoint $2L$-pieces of $\Lambda^+$ for some $k \in \N$. Then for every edge $e_0 \in E(X)$, $\axe_X(g) / \langle g \rangle$ contains at least $k/2$ edges in the orbit of $e_0$.
\end{lem}
\begin{proof}
 Let $e_0$ be an orbit of edges in $E(X)$. There exists a collapse $X \to X'$ such that $e_0$ is not sent to a single point, and 
 every edge $e \neq e_0$ in $X/G$ is not collapsible.
 
 \bigskip
 
 By Lemma \ref{lem:application-sympa} there exists a map $\tau: X' \to T_r$ such that:
 \begin{itemize}
  \item at every vertex $v \in V(X')$, there are at least two gates for the train track structure induced by $\tau$
  \item at every vertex $v$, at least two gates contain edges which are not in $G \cdot e_0$.
 \end{itemize}
 
 Let $\sigma_1, \dots, \sigma_k$ be the $2L$-pieces of $\Lambda^+$ in $\axe_{T_r} (g)$. Let $\mu_j:=e_i \dots e_l$ be a minimal edge path in $\axe_{X'} (g)$ such that $[\tau(\mu_j)]$ contains $\sigma_j$ (see Figure \ref{fig:2Lpieces}). We claim that for all $j \in \{1, \dots, k\}$, the edge path $\mu_j$ contains one edge in $G \cdot e_0$. By contradiction, assume that for some $j \in \{1, \dots, k\}$ no edge in $\mu_j$ is in $G \cdot e_0$.
 
 The path $\mu_j$ can be completed into a bi-infinite line $\eta \subset X' \setminus G \cdot e_0$ such that every turn of $\eta$ is $\tau$-legal, apart from turns in the interior of $e_i \dots e_l$. By minimality of $\mu_j$, the image $\tau(e_i)$ (resp. $\tau(e_l)$) is not contained in $\tau(e_{i+1} \dots e_l)$ (resp. $\tau(e_i \dots e_{l-1})$). Therefore the legality property of $\eta$ implies that the image $[\tau(\eta)]$ contains the segment $[\tau(\mu_j)]$.
 
 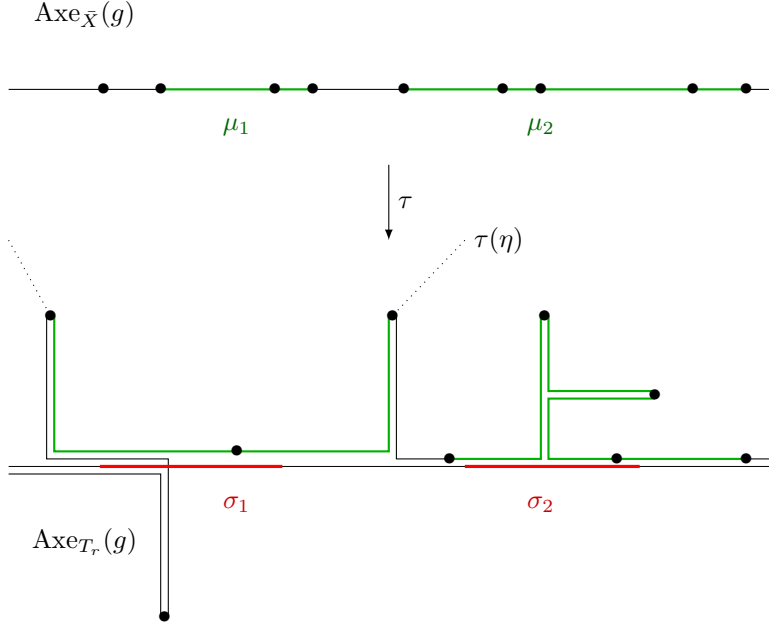
\begin{figure}
  \centering
    \begin{tikzpicture}[scale= 1]

\draw (1,6) node {$\axe_{\bar X}(g)$};
\draw (1,-1) node {$\axe_{T_r}(g)$};
\draw [-latex] (5,4) -- (5,3) node [midway, right]  {$\tau$};
\draw (0,5) -- (10,5);
\draw [thick, color = green!70!black] (2,5) -- (4,5);
\draw [thick, color = green!70!black] (5.2, 5) -- (9.7, 5);

\draw [green!40!black] (3,4.5) node {$\mu_1$};
\draw  [green!40!black] (7, 4.5) node {$\mu_2$};

\draw (1.25, 5) node {$\bullet$};
\draw (2, 5) node {$\bullet$};
\draw (3.5, 5) node {$\bullet$};
\draw (4, 5) node {$\bullet$};
\draw (5.2, 5) node {$\bullet$};
\draw (6.5, 5) node {$\bullet$};
\draw (7,5) node {$\bullet$};
\draw (9, 5) node {$\bullet$};
\draw (9.7, 5) node {$\bullet$};

\draw [very thin] (0,0) -- (10,0);

\draw [thick, color=green!70!black] (0.6,2) |- (3,0.2) ;
\draw [thick, color=green!70!black] (3, 0.2) node [color=black] {$\bullet$} -| (5, 2) -- (5.05,2) node [color=black] {$\bullet$};

\draw   [thick, color=green!70!black] (8.5, 0.9) -|  (7.1, 0.1) -- (8, 0.1) node [color=black] {$\bullet$} -- (9.7, 0.1) ;
\draw   [thick, color=green!70!black](8.5, 0.95) node [color=black]{$\bullet$} -- (8.5, 1) -|  (7.1, 2) ;
\draw  [thick, color=green!70!black] (5.8, 0.1) -|  (7,2) -- (7.05, 2) node [color=black] {$\bullet$};

\draw (0,-0.1) -| (2, -2) -- (2.05,-2) node {$\bullet$};
\draw (0.55,2) node {$\bullet$} --(0.5,2)   |- (1.2, 0.1)  -| (2.1, -2) ;

\draw (5.8, 0.1) node {$\bullet$} -|  (5.1,2);

\draw (9.7, 0.1) node {$\bullet$} -- (10, 0.1);

\draw[very thick, color = red] (1.2,0) -- (3.6,0);
\draw[very thick, color = red] (6,0) -- (8.3,0);

\draw [red!80!black] (3, -0.5) node {$\sigma_1$};
\draw [red!80!black] (7, -0.5) node {$\sigma_2$};

\draw[dotted] (0,3) -- (0.55, 2);
\draw[dotted] (6,3) node [right] {$\tau(\eta)$} -- (5.05, 2);

\end{tikzpicture}
  \caption{Axis of $g$ in $T_r$; the $2L$-pieces are in thick red, the image of $\axe_{X'}(g)$ in dotted line and the $\mu_j$ are highlighted in green. The line $\tau(\eta)$ is also represented. Due to minimality of $\mu_j$, $\eta$ contains $\sigma_j$.} \label{fig:2Lpieces}
 \end{figure}

 \bigskip
 
 Now we would like to find the axis of an element $h \in G$ such that $\axe_{T_r}(h)$ contains $[\tau(\mu_j)]$ and $\axe_{X'}(h) \cap G\cdot e_0 = \varnothing$. Suppose we find such an $h$. Then $\axe_{T_r}(h)$ contains an $L$-piece of $\Lambda^+$. However, in $X'$, there exists a candidate $u \in G$, possibly equal to $h$, whose axis does not cross $G \cdot e_0$. The assumptions of the lemma imply that $\axe_{T_r}(u)$ contains an $L$-piece of $\Lambda^-$. By Proposition \ref{prop:candidats-lamination} the pair $\{h, u\}$ is not simple. This contradicts the fact that their axes in $X'$ both avoid $G \cdot e_0$. The conclusion is that $e_0$ must appear somewhere in $\mu_j$.

 Let us explain how we construct $h$.
 In the special case where there exists an edge $e \in E(X')$ such that there are two translates $e, he$ in $\eta$ with same orientation, and one on each side of $\mu_j$ then $[e, he]$ contains a fundamental domain for the axis of $h$ and again by minimality of $\mu_j$, $\axe_{T_r}(h)$ contains $[\tau(\mu_j)]$. 
 
 In the general case, since no edge in $X' \setminus G \cdot e_0$ is collapsible, the connected component of $X' \setminus G \cdot \mathring{e_0}$ containing $\eta$ has no valence 1 vertex. Its stabilizer is a cyclic factor $H$ and this connected component is the minimal subtree $X'_H$. For every $g \in G$, $gX'_H \cap X'_H \neq \varnothing \Rightarrow g \in H$. The subtree $X'_H$ has infinite diameter because it contains $\eta$. As $X'/G$ is finite, there exists a vertex with unbounded $H$-orbit so $H$ is not elliptic. Thus it is not cyclic.

 If $H$ is not solvable, then the action of $H$ on $X'_H$ is irreducible: for every segment $I \subset X'_H$, there exists $h \in H$ whose axis contains $I$.
 
 Let $I \subset \eta$ be a segment containing a $2\BBT(\tau)/m$-neighbourhood of $\mu_j$, where $m = \min_{e \notin G \cdot e_0} \frac{\len(\tau(e))}{\len(e)}$. Let $h \in H$ be a loxodromic element whose axis contains $I$.
 
 The cancellation in $\tau(\axe_{X'}(h))$ does not reach $\mu_j$ so $\axe_{T_r}(h)$ contains $[\tau(\mu_j)]$.
 
 \bigskip
 
 Finally we must deal with the case where $H$ is isomorphic to $\BS(1,n)$. The subtree $X'_H$ is reduced, so $X'_H/H$ consists in a single edge. If $n = \pm 1$ then $X'_H$ is a line. If $h\in H$ is a loxodromic element then its axis contains $\eta$. Moreover since $X_H$ has only valence 2 vertices, they have to belong to different gates so all turns are $\tau$-legal. Therefore $\axe_{T_r}(h)$ contains $[\tau(\mu_j)]$. 
 
 If $|n|> 1$ then $X_H$ is not a line but there is a fixed point $\xi$ in $\partial X_H$ for the action of $H$. If the line $\eta$ has both endpoints different from $\xi$ then it might be impossible to find $h$ containing $\mu_j$ as a whole. However $\eta$ contains only one orbit of edge $e$. Up to reversing the orientation of $e$ we may assume $G_e = G_{t(e)}$, so every turn of the form $\{e, a \bar e\}$ with $a \in G_{t(e)}$ is degenerate. Therefore $\eta$ maps to $\dots \bar e \bar e \bar e \dots \bar e e \dots e e e \dots $ in $X'_H / H$. Since $e_i \dots e_l$ has length $> 2L$ there exists a subsegment $\eta_0$ with length $> L$ of the form $ e e e \dots$ or $\bar e \bar e \bar e \dots$ Once again such a segment is $\tau$-legal otherwise there would only be one gate at the vertices of $X_H$. There exists a loxodromic element $h \in H$ such that $\axe_{X'} (h)$ contains $\eta_0$. Therefore $\axe_{X'}(h)$ contains a $L$-piece of $\Lambda^+$ so once again we can apply the discussion above.
 
 \bigskip
 
 We proved that for any $j \in \{1, \dots, k\}$, there exists a translate of $e_0$ in $X'$ such that the minimal edge subpath $\mu_j$ contains a translate of $e_0$. 
 
\begin{figure}
 \centering
   \begin{tikzpicture}
    \draw[dotted] (0, 0) -- (5,0);
    
   \draw (0, 0.1) -| ( 1, 1)-- (1.05,1) node {$\bullet$} --(1.1,1) |- (1.5, 0.1);
   \draw (1.5, 0.1) -| ( 2, 1)-- (2.05,1) node {$\bullet$} --(2.1,1) |- (2.5, 0.1);
   \draw (2.5, 0.1) -| ( 3, 1)-- (3.05,1) node {$\bullet$} --(3.1,1) |- (3.5, 0.1);
   \draw (3.5, 0.1) -| ( 4, 1)-- (4.05,1) node {$\bullet$} --(4.1,1) |- (5, 0.1);

 \draw [very thick, red] (2.8, 0) -- (4.3, 0);
 \draw [very thick, red] (2.4, 0) -- (1.2, 0);

\draw (3,3) node {$\mu_i$s may overlap\textellipsis};
\draw (0,2) -- (5,2);
\draw (-0.5, 2) node {$\bar X$};
\draw (-0.5, 0) node {$T_r$};
\draw[green!70!black, very thick] (1,2.2)-- (3,2.2);
\draw[green!70!black, very thick] (2,2.1)-- (5,2.1);

\draw (1,2) -- +(0,0.5);
\draw (2,2) -- +(0,0.5);
\draw (3,2) -- +(0,0.5);
\draw (4,2) -- +(0,0.5);

\begin{scope}[xshift= 5cm]
     \draw[dotted] (1, 0) -- (5,0);

   \draw (1, 0.1) -| ( 1.5, 1)-- (1.55,1) node {$\bullet$} --(1.6,1) |- (2.5, 0.1);
   \draw (2.5, 0.1) -| ( 3, 1)-- (3.05,1) node {$\bullet$} --(3.1,1) |- (3.5, 0.1);
   \draw (3.5, 0.1) -| ( 4, 1)-- (4.05,1) node {$\bullet$} --(4.1,1) |- (5, 0.1);

 \draw [very thick, red] (2.8, 0) -- (4.3, 0);
 \draw [very thick, red] (1.8, 0) -- (1.2, 0);
 \draw [very thick, red] (2.1, 0) -- (2.6, 0);

\draw (2.8,3) node {but may not be one edge};
\draw (1,2) -- (5,2);
\draw[green!70!black, very thick] (1,2.2)-- (3,2.2);
\draw[green!70!black, very thick] (2,2.1)-- (5,2.1);
\draw[green!70!black, very thick] (2,2.3)-- (3,2.3);

\draw (2,2) -- +(0,0.5);
\draw (3,2) -- +(0,0.5);
\draw (4,2) -- +(0,0.5);

\end{scope}
\end{tikzpicture}
 \caption{The subsegments of $\axe_X(g)$ whose image contain a $2L$-piece of $\Lambda^+$ may overlap} \label{fig:cas-partage-e_0}
\end{figure}
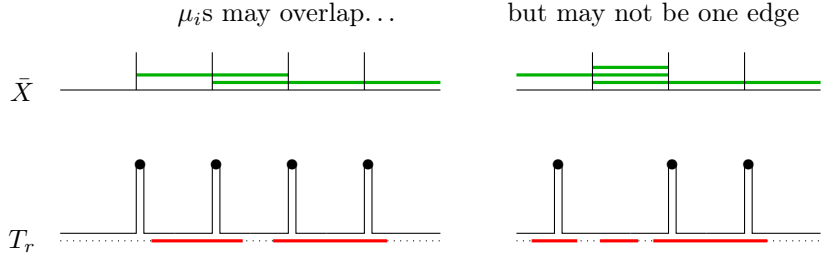

 If the segments $\mu_j$, $j \in \{1, \dots, k\}$ are disjoint, then we are done. This may fail though. We will see that $\mu_j \cap \mu_{j'} = \varnothing$ if $|j- j'|\geq 2$ and may be a single edge if $|j-j'|=1$. Thus when counting the translates of $e_0$ in the $\mu_j, j \in \{1, \dots, k\}$, a translate may be counted more than once, but it can be counted for at most twice (see Figure \ref{fig:cas-partage-e_0}). Therefore $\axe_{X'}(g)/ g$ contains at least $k/2$ translates of $e_0$. Since this lifts to $X$ we get the lemma.
 
 \bigskip
 
 Let us prove the fact about the intersection of the segments $\mu_j$.

 First we prove that a $2L$-piece of $\Lambda^+$ cannot be contained in the $\tau$-image of a single edge $e$. By contradiction, suppose otherwise: again we will construct a simple pair of elements containing long pieces of opposite laminations. The edge $e$ must be in the orbit of $e_0$. As above, one can find $\tau$-legal turns $\{e, e_1\}$ and $\{\bar e, e_2\}$ with $e_1, e_2 \notin G \cdot e_0$. Let $e'$ be an edge of $X' \setminus G \cdot e_0$. Since $e'$ is not collapsible, there exists a translate $h e' \neq e'$ such that $\{e', he'\}$ (if $e'$ is not a loop in the quotient) or $\{e', h \bar e'\}$ (if $e'$ is a loop) is a non-degenerate turn. Define 
 \[
  \rho_i:= e_i \cdot h_i \bar e_i \cdot h_i h_i' e_i \cdot h_i h_i' h_i \bar e_1 \dots
 \]
 if $e_i$ is not a loop, where $h_i, h_i'$ are such that $\{e_i, he_i\}$ and $\{\bar e_i, h' \bar e_i\}$ are non-degenerate. If $e_i$ is a loop define
 \[
  \rho_i := e_i \cdot h_i e_i \cdot h_i^2 e_i \dots
 \]

 Suppose $e_i$ is not a loop. Then $e_i$ cannot be identified with $h_i e_i$ by $\tau$, because the vertex groups at $t(e_i)$ and $h_i t(e_i)$ are not nested: it would imply that $T_r$ has an elliptic element which is not elliptic in $X'$. There exists a subdivision of $X'$ such that $e_i = a_i \cdot e_i' \cdot b_i$, with $\tau(b_i) = \tau(h_i b_i)$ and $\tau(a_i)= \tau(h_i'a_i)$. Then
 \[
  [\tau(\rho_i)] := \tau(a_i) \cdot \tau(e_i') \cdot h_i \tau(\bar e_i') \cdot h_i h_i' \tau(e_i') \dots
 \]
 and since the turns between $e$ and $e_i$ are legal, there is no simplification between $\tau(a_i)$ and $\tau(e)$.

 If $e_i$ is a loop, then we already proved in Lemma \ref{lem:application-sympa} that $e_i$ and $h_i\bar e_i$ cannot have the same image by $\tau$, or we would obtain new elliptic elements. Once again there exists a subdivision $e_i = a_i \cdot e_i' \cdot b_i$ with $\tau(\bar b_i) = h_i \tau(a_i)$. 
 Thus
 \[
  [\tau(\rho_i)] := \tau(a_i) \cdot \tau(e_i') \cdot h_i \tau(e_i') \cdot h_i^2 \tau(e_i') \dots
 \]
 and once again there is no cancellation between $\tau(a_i)$ and $\tau(e)$.
 
 Consequently we can construct a bi-infinite geodesic $\bar \rho_1 \cdot e \cdot \rho_2$ such that $\rho_1, \rho_2$ are rays which cross only one orbit of edges, and $\tau(e) \subset [\tau(\bar \rho_1 \cdot e \cdot \rho_2)]$. We proved above that the rays need not be legal, the point is that the cancellation which may occur at turns remains controlled.
 
 Let $l$ be such that the image of any segment of longer than $l$ by $\tau$ is longer than $2 \Lip(\tau)$. Let $\rho_1^0, \rho_2^0$ be prefixes of the rays longer than $l$. The path $\bar \rho_1^0 \cdot e \cdot \rho_2^0$ can be closed into a loop in the quotient, representing an element $h \in G$ such that $\axe_{T_r}(h)=[\tau(\axe_{X}(h))]$ contains $\tau(e)$, hence a $2L$-piece of $\Lambda^+$. The point is that this loop may be constructed such that it crosses only three orbits of edges in $X'$. Since $b_1(X'/G) \geq 3$, $\axe_{X'}(h)$ must avoid one orbit of edges.
 
 There exists a candidate $u$ of $X'$ whose axis avoids the same orbit of edge as $h$, thus $\{u,h\}$ is a simple pair. By assumption $\axe_{T_r}(u)$ contains an $L$-piece of $\Lambda^-$, which is a contradiction to Proposition \ref{prop:candidats-lamination}: thus $\tau(e)$ cannot contain a whole $2L$-piece of the lamination.
 
 \bigskip

 The second point is that the intersection of segments $\mu_i, \mu_j$ cannot be more than one edge:
 by minimality of $\mu_j$, the last endpoint of $\sigma_j$ lies in the image of the last edge of $\mu_j$ but not in the image of any other edge. Similarly for $j' > j$, the first point of $\sigma_{j'}$ lies in the image of the first edge of $\mu_{j'}$ and not in any other edge. Thus the intersection $\mu_i \cap \mu_j$ is at most one single edge.
 
 Besides, if $|j'-j|\geq 2$, then $\sigma_{j+1}$ must be contained between the last point of $\sigma_j$ and first point of $\sigma_{j'}$. This is not possible if both belong to the same edge, hence the fact.
\end{proof}

Now we have sufficient tools to prove Lemma \ref{lem:projections2}.

\begin{proof}[Proof of Lemma \ref{lem:projections2}]
 Let $X, Y$ be as in the statement of the lemma. Assume $t_Y < t_X$ in $\L_f$. The other case works similarly by exchanging the roles of $\phi$ and $\phi^{-1}$ and will give other constants $s, c$: we will take the greater constants.

 There exists $s_1$ such that for any $s > s_1$, for any $t \in \R$, for any candidate $g$ of $T_{t}$, $\LEG_f(g, T_{t+s}) > \varepsilon_0$. This is a consequence of the following facts: $T_t$ has an $f$-legal candidate $u$ (Lemma \ref{lem:existe-g-legal}) so $t_0(u) \leq t$, and there exists $s >0$ such that for any other candidate $v$ we have $t_0(v) \leq t_0(u) + s$  (Lemma \ref{lem:distance-gh}). 
 
 Let $L$ be the constant from Proposition \ref{prop:candidats-lamination}. By Lemma \ref{lem:lamination-ratio} there exists $N_1 \in \N$ such that for every candidate $g$ of $X$, $\LR(\phi^{N_1}(g), T, \Lambda^+, 2L) > \varepsilon_0/4$. 
 
 There exists $s_2>0$ such that for any $t \in \R$, the image of any candidate of $T_t$ in $T_{t-s_3}$ contains an $L$-piece of $\Lambda^-$.
 
 \bigskip
 
 Define $d=s_1 + N_1 \log(\lambda) + s_2$. Suppose $t_X - t_Y > d$. Let $r = t_X - s_2$. Let $g \in G$ be a candidate of $Y$ which realizes $\Lip(Y, \pi(X))$. Then the axis of $g$ in $T_r$ contains long leaf segments of $\Lambda^+$. Actually $\LR(g, T_r, \Lambda^+, 2L) \geq \varepsilon_0/4$. A given leaf segment of length longer than $2L$ can be at least half covered with disjoint $2L$-pieces of $\Lambda^+$. Thus a proportion of at least $\varepsilon_0/8$ of $\axe_{T_r}(g)$ can be covered by disjoint $2L$-pieces of $\Lambda^+$.

 Let $k(r)$ be the number of disjoint $2L$-pieces of $\Lambda^+$ which tile $\axe_{T_r}(g)/\langle g \rangle$. We have 
  \[
  2L k(r) > \frac{\varepsilon_0}{8} \|g\|_{T_r}
 \]

 Now $\axe_{T_r}(g)/\langle g \rangle$ contains $k(r)$ $2L$-pieces of $\Lambda^+$. By Lemma \ref{lem:passages-aretes} $\axe_{X} (g) / \langle g \rangle$ contains at least $k(r)/2$ edges in each orbit of $E(X)$. Since $\vol(X/G)= 1$ we have
 \[
  \vol( \axe_{X} (g) / \langle g \rangle) = \|g\|_X \geq k(r) / 2.
 \]
 Thus
 \[
  \|g\|_X \geq k(r)/2 \geq \|g\|_{T_r} \frac{ \varepsilon_0}{16L}
 \]
 Then 
 \begin{align*}
    \Lip(Y,X) &= \frac{\|g\|_X}{\|g\|_Y} \\
       &\geq \frac{\|g\|_{T_r}}{\|g\|_Y} \frac{ \varepsilon_0}{16L} \\
       &\geq \frac{\varepsilon_0}{16 L} \Lip(Y, T_r)
 \end{align*}
 By triangular inequality $d_{\Lip}(Y, T_r) \geq d_{\Lip}(Y, \pi(X))-r$ so
 \[
  d_{\Lip}(Y, X) \geq d_{\Lip}(Y, \pi(X)) - s_2 - \log(\frac{16L}{\varepsilon_0})
 \]
\end{proof}


\begin{defi}
 The ball of outward radius $r>0$ centered at $Y \in \D$  is 
 \[
  B_{\rightarrow}(Y, r):=\{X \in \D/ d_{\Lip}(Y, X) < r\}
 \]
\end{defi}

A \emph{closest point projection} to $\L_f$ is a map $p_f: \D \to \L_f$ such that for all $X \in \D$, the distance $d_{\Lip}(X, p_f(X))$ is minimal. The map $\pi_f$ constructed in Section \ref{sec:projection} is a closest point projection to $\L_f$.

Now we can state and prove the strong contraction property.
\begin{theo}
 Let $\phi$ be a fully irreducible automorphism such that $\phi, \phi^{-1}$ both admit train track representatives.
 
 Let $\L_f$ be an axis for $\phi$ in $\D$ and let $p_f$ be a closest point projection to $\L_f$. Then there exists $D > 0$ such that for any $Y \in \D$ and $r>0$ such that $B_{\rightarrow}(Y, r) \cap \L_f = \varnothing$
 \[
  \diam(p_f(B_{\rightarrow}(Y, r))) \leq D
 \]
\end{theo}
\begin{proof}
 Let $Y \in \D$, $r= d_{\Lip}(Y, p_f(Y))$. Let $B:=B_{\rightarrow}(Y, r)$. A ball centred at $Y$ intersects the axis if and only if its radius it greater than $r$. Balls with smaller radius are contained in $B$ so it suffices to bound $\diam(p_f(B))$ independently of $Y$. Let $X \in B$.
 Let $s, c$ be the constants from Lemma \ref{lem:projections2} and $s',c'$ be the constants from Lemma \ref{lem:projections1}. Suppose $d_{\Lip}(p_f(Y), p_f(X)) > \max\{s,s'\}$. Then Lemma \ref{lem:projections2} yields
 \[
  d_{\Lip}(Y,X) \geq d_{\Lip}(Y, p_f(X)) - c
 \]
 and using Lemma \ref{lem:projections1}: 
 \[
  d_{\Lip}(Y,X) \geq d_{\Lip}(Y, p_f(Y)) + d_{\Lip}(p_f(Y), p_f(X)) - c-c'
 \]
 Since $d_{\Lip}(Y,X) < r = d_{\Lip}(Y, p_f(Y))$ we have
 \[
  d_{\Lip}(p_f(Y), p_f(X)) \leq c+c'
 \]
 Therefore $\diam(p_f(B)) \leq 2 \max\{s, s', c+c'\}$. This bound is independent of $Y$.
\end{proof}

 \bibliography{./bibli_axes}

\begin{thebibliography}{{Coo}87}

\bibitem[{Alg}11]{AlgomKfirStrongly}
Yael {Algom-Kfir}.
\newblock {Strongly contracting geodesics in outer space}.
\newblock {\em {Geom. Topol.}}, 15(4):2181--2233, 2011.

\bibitem[BBF15]{BBF}
Mladen {Bestvina}, Ken {Bromberg}, and Koji {Fujiwara}.
\newblock {Constructing group actions on quasi-trees and applications to
  mapping class groups}.
\newblock {\em {Publ. Math., Inst. Hautes \'Etud. Sci.}}, 122:1--64, 2015.

\bibitem[{Bes}11]{BestvinaBersLike}
Mladen {Bestvina}.
\newblock {A Bers-like proof of the existence of train tracks for free group
  automorphisms.}
\newblock {\em {Fundam. Math.}}, 214(1):1--12, 2011.

\bibitem[BF91]{BestvinaFeighnBounding}
Mladen Bestvina and Mark Feighn.
\newblock Bounding the complexity of simplicial group actions on trees.
\newblock {\em Invent. Math.}, 103(3):449--469, 1991.

\bibitem[BFH97]{BestvinaFeighnHandelLaminations}
M.~Bestvina, M.~Feighn, and M.~Handel.
\newblock Laminations, trees, and irreducible automorphisms of free groups.
\newblock {\em Geom. Funct. Anal.}, 7(2):215--244, 1997.

\bibitem[BH92]{BestvinaHandelTrainTracks}
Mladen Bestvina and Michael Handel.
\newblock Train tracks and automorphisms of free groups.
\newblock {\em Ann. of Math. (2)}, 135(1):1--51, 1992.

\bibitem[Bou16]{BouetteThese}
Margot Bouette.
\newblock {\em {On the growth of the automorphisms of Baumslag-Solitar
  groups}}.
\newblock Theses, {Universit{\'e} Rennes 1}, December 2016.

\bibitem[{Coo}87]{CooperAutomorphisms}
Daryl {Cooper}.
\newblock {Automorphisms of free groups have finitely generated fixed point
  sets}.
\newblock {\em {J. Algebra}}, 111:453--456, 1987.

\bibitem[FM15]{FrancavigliaMartino}
Stefano {Francaviglia} and Armando {Martino}.
\newblock {Stretching factors, metrics and train tracks for free products}.
\newblock {\em {Ill. J. Math.}}, 59(4):859--899, 2015.

\bibitem[{For}02]{ForesterDeformation}
Max {Forester}.
\newblock {Deformation and rigidity of simplicial group actions on trees.}
\newblock {\em {Geom. Topol.}}, 6:219--267, 2002.

\bibitem[For06]{ForesterSplittings}
Max Forester.
\newblock Splittings of generalized {B}aumslag-{S}olitar groups.
\newblock {\em Geom. Dedicata}, 121:43--59, 2006.

\bibitem[GJLL98]{GJLL}
Damien {Gaboriau}, Andre {Jaeger}, Gilbert {Levitt}, and Martin {Lustig}.
\newblock {An index for counting fixed points of automorphisms of free groups}.
\newblock {\em {Duke Math. J.}}, 93(3):425--452, 1998.

\bibitem[GL07]{GuirardelLevitt07}
Vincent Guirardel and Gilbert Levitt.
\newblock Deformation spaces of trees.
\newblock {\em Groups Geom. Dyn.}, 1(2):135--181, 2007.

\bibitem[Mei15]{MeinertTheLipschitzMetric}
Sebastian Meinert.
\newblock The {L}ipschitz metric on deformation spaces of {$G$}-trees.
\newblock {\em Algebr. Geom. Topol.}, 15(2):987--1029, 2015.

\bibitem[Pap21]{PapinWhitehead}
Chloé Papin.
\newblock Whitehead algorithm for automorphisms of generalized baumslag-solitar
  groups, 2021.

\bibitem[Pap22]{PapinDetection}
Chloé Papin.
\newblock Detection of fully irreducible automorphisms in generalized
  baumslag-solitar groups, 2022.

\bibitem[{Sen}81]{Seneta}
E.~{Seneta}.
\newblock {\em {Non-negative matrices and Markov chains. 2nd ed}}.
\newblock Springer, New York, NY, 1981.

\bibitem[Ser77]{SerreArbresAmalgames}
Jean-Pierre Serre.
\newblock {\em Arbres, amalgames, {${\rm SL}_{2}$}}.
\newblock Soci\'{e}t\'{e} Math\'{e}matique de France, Paris, 1977.
\newblock Avec un sommaire anglais, R\'{e}dig\'{e} avec la collaboration de
  Hyman Bass, Ast\'{e}risque, No. 46.

\end{thebibliography}

 \end{document}